\documentclass{amsart}

\newtheorem{theorem}{Theorem}[section]
\newtheorem{lemma}[theorem]{Lemma}

\theoremstyle{definition}
\newtheorem{definition}[theorem]{Definition}

\theoremstyle{remark}

\numberwithin{equation}{section}
\numberwithin{figure}{section}

\usepackage{graphicx}
\usepackage{epsfig}
\usepackage{epsf}

\newtheorem{corollary}{Corollary}[section]

\newcommand{\cA}{{\cal A}}
\newcommand{\cP}{{\cal P}}
\newcommand{\conv}{\mbox{conv}}

\newcommand{\eps}{\epsilon}
\newcommand{\bitem}{\begin{itemize}}
\newcommand{\eitem}{\end{itemize}}
\newcommand{\goto}{\rightarrow}

\newcommand{\beqn}{\begin{equation}}
\newcommand{\eeqn}{\end{equation}}
\newcommand{\balign}{\begin{align}}
\newcommand{\ealign}{\end{align}}

\newcommand{\bR}{{\bf R}}

\newcommand{\cF}{{\cal F}}
\newcommand{\cX}{{\cal X}}

\def\cal{\mathcal}

\def\rhopW{\rho^{+}_W}
\def\rhopmW{\rho^{\pm}_W} 
\def\rhopS{\rho^{+}_S}
\def\rhopmS{\rho^{\pm}_S} 
\def\c{\tau}

\def\cs{\tau}
\def\cw{\tau}

\def\rstar{r^{\star}(\delta)}

\def\rstars{r^{\star}_S(\delta)}
\def\rstarsdn{r^{\star}_S(\delta_n)}
\def\rpmsdn{r^{\pm}_S(\delta_n)}
\def\rps{r^{+}_S(\delta)}
\def\rpms{r^{\pm}_S(\delta)}

\def\rw{r_W^{\star}(\delta)}
\def\rwp{r_W^{+}(\delta)}
\def\rwpm{r_W^{\pm}(\delta)}
\def\rwdn{r_W^{\star}(\delta_n)}

\def\rpsc{r^{+}_S(\delta;\cs)}
\def\rpmsc{r^{\pm}_S(\delta;\cs)}
\def\rstarsnc{r^{\star}_S(n/N_n;\cs)}

\def\rwcp{r_W^{+}(\delta;\cw)}
\def\rwcpm{r_W^{\pm}(\delta;\cw)}

\def\rwnc{r_W^{\star}(n/N_n;\cw)}

\def\psicomp{\Psi_{com}^+}
\def\psicompm{\Psi_{com}^{\pm}}
\def\psiint{\Psi_{int}}
\def\psiintp{\Psi_{int}^+}
\def\psiintpm{\Psi_{int}^{\pm}}
\def\psiextp{\Psi_{ext}^+}
\def\psiextpm{\Psi_{ext}^{\pm}}
\def\psinet{\Psi_{net}}
\def\psinetp{\Psi_{net}^+}
\def\psinetpm{\Psi_{net}^{\pm}}

\def\psiweakp{\Psi_{face}^+}
\def\psiweakpm{\Psi_{face}^{\pm}}

\def\xnu{x_{\nu}}
\def\xdelta{x_{\delta}}
\def\txnu{\tilde x_{\nu}}
\def\ynu{y_{\nu}}
\def\zp{z^+}
\def\zpd{z^+_{\delta}}
\def\zpm{z^{\pm}}
\def\zpmd{z^{\pm}_{\delta}}
\def\ygam{y_{\gamma}}
\def\sgam{s_{\gamma}}

\def\tsgam{\tilde s_{\gamma}}

\begin{document}

\title[Counting faces of randomly-projected polytopes]{Counting faces of randomly-projected polytopes
when the projection radically lowers dimension}



\author{David L. Donoho}
\address{Department of Statistics, Stanford University}
\curraddr{Department of Statistics, Stanford University}
\email{donoho@stanford.edu}
\thanks{DLD acknowledges partial support from NSF DMS 05-05303, and
01-40698 (FRG), and NIH}

\author{Jared Tanner}
\address{Department of Statistics, Stanford University}
\curraddr{Department of Statistics, Stanford University}
\email{tanner@math.utah.edu}
\thanks{JT acknowledges support from NSF fellowship DMS 04-03041.}

\subjclass[2000]{52A22, 52B05, 52B11, 52B12, 62E20, 68P30, 68P25, 68W20, 68W40, 94B20  94B35, 94B65, 94B70}


\date{November 2005}



\maketitle



\section{Introduction}
\setcounter{equation}{0}
\setcounter{table}{0}
\setcounter{figure}{0}

\subsection{Three surprises of high dimensions}

This paper develops asymptotic methods to
count faces of random high-dimensional polytopes;
a seemingly dry and unpromising pursuit. Yet our conclusions
have surprising implications \-- in  statistics, probability, information theory, and signal processing \-- with
potential impacts in practical subjects like medical imaging and digital
communications.  Before involving the
reader in our lengthy analysis of high-dimensional
face counting, we describe three implications of our results.

\subsubsection{Convex Hulls of Gaussian Point Clouds}\label{int:Hull}

Consider a random point cloud  of $n$ points $x_i$, $i=1,\dots,n$,
sampled independently and identically from a 
Gaussian distribution in $\bR^d$ with nonsingular covariance.
This is a standard model of multivariate data; its properties
are increasingly important in a wide range of applications.
At the same time, it is an attractive and in some sense
timeless object for theoretical study.

Properties of the convex hull of the random point cloud  $\cX = \{ x_i \}$
have attracted interest for several decades, increasingly so in recent years;
there is a now-voluminous literature on the subject. 
The results could be significant for understanding outlier detection,
or classification problems in machine learning. 

A classical asymptotic result, \cite{Hueter}, holds that if the 
dimension $d$ stays fixed, while the number
of points $n \goto \infty$, the convex hull has $\sim c_d \log^{(d-1)/2}(n)$
vertices, and the remaining points of $\cX$ are
all of course in the interior of the convex hull.

The modern trend in statistics and probability is
to consider the case where both the number of dimensions $d$
and the sample size $n$ are large \cite{HMN05,HW04}. 
In that case, the intuition
fostered by the classical fixed-dimension asymptotic
is wildly inaccurate.  Rather than the relatively 
few extreme points that we saw
in the fixed-dimension asymptotic, there are now
many extreme points, many edges, etc. -- in fact, the
maximal number conceivable. More precisely,
let $k_d^* = k_d^*(\cX)$ denote the largest
number $k$ such that
\bitem
\item Each point $x_j$ is a vertex of $\conv(\cX)$;
\item Each line segment $[x_j,x_i]$, $j \neq i$ is an edge
of $\conv(\cX)$;
\item ...
\item Every $k+1$ distinct points of $\cX$ span a $k$-face 
of $\conv(\cX)$.
\eitem
Then, not only is $k_d^*$ defined and positive,
it is rather large. Section \ref{sec:AppGaussianProof} below
gives a corollary of our main results
saying roughly that, for $\eps > 0$,
with overwhelming probability for large $d$,
\begin{equation} \label{result1}
   k_d^*  >  \frac{d}{2e \log(n/d)}(1-\eps).
\end{equation}
Not only are no points of $\cX$ `inside' $\conv(\cX)$,
it is also true that no edge between any pair of
points crosses the interior of $\conv(\cX)$, etc.
This is about as far from low-dimensional 
intuition as it is possible to get!

\subsubsection{Signal Recovery from Random Projections}\label{int:signal}

Suppose we are interested in a vector $x_0 \in \bR^N$
which, although unknown to us,  is known to be $k$-sparse \--
i.e. we know that it has at most $k$ nonzeros when represented
in the standard basis.  We are allowed to
ask some number $n$ of `questions' about $x_0$, each question
coming in the form of a projection $y_i = \langle a_i, x_0 \rangle$
on a vector $a_i \in \bR^N$.
How big should $n$ be so that we may recover $x_0$, 
i.e.: ``How many questions suffice to recover a $k-$sparse vector''?  

Obviously, $N$ suffice (simply ask for the $N$ coordinates
in the standard unit vector basis),
but in cases where $x_0$ is very
sparse, $k \ll N$, many fewer questions 
will do. Indeed, $n = 2k+1$ suffice;
simply take the $a_i$ as independent random vectors with iid Gaussian entries. 
(The matrix $A$ having $a_i$ for rows will then 
have its columns in general position, which implies that
there cannot be two $k$-sparse vectors $x_0$ and $x_1$ 
both answering the questions in the same way 
\cite{donoho1}.)  Although such a random set of questions
determines $x_0$ uniquely,
the task of actually recovering $x_0$
from such information is daunting; in general, one must enumerate
the $k$-subsets of columns of $A$ looking for a subset
which can be combined linearly to generate $y$.

A more useful question:
how many questions are needed in order to
permit {\it computationally tractable} recovery of $x_0$?
We will give precise and simply-stated
results for reconstruction using
standard linear programming.

Generate $n$ questions `at random' by simply taking for $A$ an
$n$ by $N$ matrix with iid Gaussian $N(0,1/n)$ entries. 
Obtain  a vector of $n$ measurements $y = Ax_0$ where $x_0$ has $k$ nonzeros.
Consider the convex optimization problem
\[
    (P_1)\quad\quad  \min \| x\|_1 \mbox{ subject to }  y = Ax.
\]
If $n$ is large enough relative to $k$, then
the solution $x_1$ to $(P_1)$ is very likely to be {\it exactly} $x_0$.
Section \ref{sec:AppHowManyProof} below gives
a corollary of this paper's main results showing that,
for $N$ much larger than $k$, and both large, this
exact equality happens as soon as 
\begin{equation} \label{result2}
        n \geq   2k \cdot \log(N/n)(1 + o_p(1)).
\end{equation}

Thus if we sample not $2k+1$ projections but instead roughly  $2 k \log(N/n)$
we can efficiently reconstruct the $k$-sparse vector; and this can be
far fewer than the $N$ samples superficially required.


\subsubsection{How many gross errors can we efficiently correct?}\label{int:ecc}

Consider a stylized problem of transmitting $m$ `pieces' of information
\-- i.e. $m$ numbers \-- with immunity to occasional
transmission errors.  A standard strategy is
encode the data to be transmitted as a block
of $N > m$ numbers, and to decode the received
block.  Let $B$ be an $m \times N$ matrix.
Given a vector $u \in \bR^m$ to be transmitted, encode it
as $v= B^T u \in \bR^N$ and transmit.  The receiver 
measures $w = v+ z$ where $w \in \bR^N$ and $z$ 
represents transmission errors. The receiver in some way decodes 
the $N$ numbers, hoping to produce the $m$ original entries in $u$.

The nonzeros in $z$ represent transmission errors; call
the number of nonzeros $k$. 
How many errors can such a scheme tolerate?  In principle,
if $B$ is well-chosen and $N - m = 2k+1$, it is possible to 
correct $k$ errors.  To do so, the receiver executes a
combinatorial search through all possible locations
of the $k$ nonzeros among the $N$ received values,
to see which error pattern is consistent with the underlying model.
Unfortunately, such a brute-force scheme is impractical for 
all but the smallest $N$.
More to the point is the question of how many errors
a practical decoding scheme can tolerate.
 
A simple decoding scheme based on $(P_1)$
can be used if the  encoding matrix $B$ is generated
in a specific way.
Let $U$ be a random orthogonal matrix, uniformly-distributed
 on $O(N)$, and partition it as $ 
U =  \begin{pmatrix}
     A  \\
     B 
\end{pmatrix}
 $
 where the encoding matrix $B$ is $m \times N$ and 
 the generalized checksum matrix $A$ is $n  \times N$,
 with $m  + n = N$.
 Given the received data $w$, form the generalized checksum
 $y = A w$. Solve the instance of $(P_1)$ given by $(y,A)$, obtaining
 $x_1$.
The generalized checksum is used to estimate
 the error pattern, and the optimization result $x_1$ is
 our estimate of $z$. Reconstruct by subtracting this estimate of the error
out of the received message, and projecting down from $\bR^N$
 to $\bR^m$:  $u_1 = B(w-x_1)$.
 
As $(P_1)$
is a standard convex optimization problem, this can be considered
computationally tractable.  How many errors can this scheme tolerate?

To answer this quantitatively, let us  call $R = n/N$ the {\em rate} of the
code, and consider the regime of high-rate coding, where $R$ is
nearly one. In this regime we don't want to expand the block length
by very much in our encoding, but we still 
want to gain some immunity to errors.

The results just stated in Section \ref{int:signal},
and a corollary in Section 6.4 below, together imply the following. 
Consider a sequence of problems $(n,N_n)$
with $R_n = n/N_n \goto 1$ sufficiently slowly.
Suppose the error vector $z$ contains
$k$ nonzeros and is stochastically independent of $(A,B)$,
so the sites and signs of the nonzeros are random and independent of $A$.
There is {\it perfect recovery} $u_1 = u$ provided $k  \leq k_n^W$,
where $k_n^W$ is a random variable dependent on $(A,B)$, and obeying
\begin{equation} \label{randomerrors}
    k_n^W  =   n / (2 \log(1/(1-R_n))) (1+ o_p(1)),  \qquad n \goto \infty.
\end{equation}
In short, if we use very long blocks, and stipulate a very small
loss in transmission rate $R_n = 1 - \eps_n$, with $\eps_n$
small, we can use linear programming to 
correct about $n/2\log(\eps_n)$ errors.

Results to be stated below -- see Section \ref{sec:AppAllECCProof} 
-- imply an even more impressive
result. Again, consider a sequence of problems $(n,N_n)$
with $R_n = n/N_n \goto 1$ sufficiently slowly.
Suppose the error vector $z$ contains
$k$ nonzeros {\it at arbitrary sites} and with {\it arbitrary nonzeros}.
There is {\em perfect} recovery $u_1 = u$ provided $k  \leq k_n^S$,
where $k_n^S$ is a random variable dependent on $(A,B)$ and obeying
\begin{equation} \label{arberrors}
    k_n^S  \geq  n/ (2e \log(\sqrt{\pi}/(1-R_n)) (1+ o_p(1))),  
    \qquad n \goto \infty.
\end{equation}
In short, if we use very long blocks, and stipulate a very small
loss in transmission rate $R_n = 1 - \eps_n$, with $\eps_n$
small, we can use linear programming to 
correct all possible patterns of about $n/2e\log(\eps_n)$ errors.

Note that the sites and values of the errors can here be arbitrary;
they can be chosen by a malicious opponent who knows $v$,$B$,$A$,
and $u$! The noise can thus be arbitrarily more energetic than
the signal, can be carefully chosen, and still it is completely suppressed.
In contrast to  (\ref{randomerrors}), which requires errors to
be in random positions, (\ref{arberrors}) allows them to occur in bursts
or in any other malicious patterns.

\subsection{Random Projections of Convex Polytopes}

The surprises  (\ref{result1}),(\ref{result2}), 
(\ref{randomerrors}),(\ref{arberrors})
are facets of a phenomenon which
makes appearances throughout the
mathematical sciences, in the fields of statistics, probability,
information theory, and signal processing. 
The phenomenon concern thresholds in the behavior  
of face counts of random high-dimensional polytopes.
We now develop the terminology and framework
for those results, only later explaining
how they imply (\ref{result1})-(\ref{arberrors}).

Let $T = T^{N-1}$ denote the standard simplex
$\{ x : \sum_i x_i = 1, x_i \geq 0 \}$ and let 
$C =C^N$ denote the standard cross-polytope 
in $\bR^N$, i.e. the  collection of vectors $ \{ x : \| x \|_1 \leq 1\}$.
(Also called the $\ell_1^N$-ball). 
Here and in what follows,
let $Q$ be either $T^{N-1}$ 
or $C^N$.

Let $A$ be
an $n \times N$ random matrix with Gaussian iid entries.
The image $AQ$ is a convex
subset of $\bR^n$; in fact, a convex polytope. 
We are interested in the case $n < N$,
so that multiplication by $A$ lowers the dimension.

It makes sense to count the number of $k$-dimensional faces of $Q$ and $AQ$,
$0 \leq k \leq n$.   In general $AQ$ will have fewer faces
than $Q$.   More precisely, if we enumerate
the $k$-faces $F$ of  $Q$, each $AF$ will either
be a face of $Q$ or will belong to the interior
of $AQ$. More picturesquely, some of the faces of $Q$
`survive projection', while some of the faces `do not survive'.

\subsubsection{Typical Faces of Random Polytopes}\label{int:typ}

The $k$-dimensional faces of $Q$ make a finite set, $\cF_{k}(Q)$ (say),
by placing uniform measure on this set, we may speak of {\em typical} 
{\em faces}, as follows.
\setcounter{theorem}{-1}
\begin{definition}
Consider a sequence of problem sizes $(n,N_n)$.  
Suppose that, for a given projector $A$, a property 
${\cal P} = {\cal P}(F;A)$ of the projected face $AF$ holds, 
at  a fraction $\pi_{k,n} = \pi_{k,n}(A)$ of $k$-faces $F \in \cF_k(Q)$. Suppose that the 
random variable $\pi_n\goto_p 1$ as $n\goto\infty$.  Then we say that 
(asymptotically) the {\em typical $k$-face $F \in \cF_k(Q)$ has property ${\cal P}$.}
\end{definition}
We now consider the fate of the typical $k$-face of  $A$ under the
projection $Q \mapsto AQ$.  In the following statements, 
fix $\epsilon>0$.
 \bitem
\item {\sl Let $F$ be a typical $k$-face of $T^{N-1}$.
Is $AF$ a face of $AT^{N-1}$?} 
The answer is {\em yes},
provided $N$ and $k$ are both large
and $n > 2k \log(N/n)(1 + \eps)$, and {\it no}
provided $n < 2k \log(N/n)(1 - \eps)$.
\item {\sl Let $F$ be a typical $k$-face of $C^N$.
Is $AF$ a face of $AC^N$?} The answer is {\it yes},
provided $N$ and $k$ are both large
and $n > 2k \log(N/n)(1 + \eps)$, and {\it no}
provided $n < 2k \log(N/n)(1 - \eps)$.
\eitem

In short, there are well-defined thresholds 
at which {\em typical}  $k$-faces of the 
simplex and the cross polytope begin to
get `swallowed up' under random lowering of dimension.

\subsubsection{All Faces of Random Polytopes}\label{int:arb}

We now consider the fate of the whole
collection of  $k$-faces simultaneously.
\bitem
\item {\sl For every $k$-face $F$  of $T^{N-1}$, is $AF$ also a $k$-face of 
$AT^{N-1}$?} The answer is overwhelmingly likely to be {\it yes},
provided $N$ and $k$ are both large
and $n > 2e k \log(N/(n \cdot 2\sqrt{\pi}))(1 + \eps)$. 
\item {\sl For every $k$-face $F$  of $C^N$, is $AF$ also a $k$-face of 
$AC^N$?} The answer is overwhelmingliy likely to be {\it yes},
provided $N$ and $k$ are both large
and $n > 2ek \log(N/(n\cdot \sqrt{\pi}))(1 + \eps)$.
\eitem
Below certain specific bounds on the
face dimension  $k$, no faces are lost in projection.
 
\subsection{Background: Proportional Growth Setting}\label{subsec:main}

 
Our promised applications, such as (\ref{result1}) and (\ref{result2}),
 were stated merely with $n$ and $N$ (respectively $k$ and $N$)
 both large. However, the backbone of our analysis (and the bulk
 of prior scholarly work) concerns a setting in which
 $(k,n,N)$ are large but also comparable in size. We consider this case
 first and later extend our results to a more general setting.
 
\begin{definition}
A sequence of triples $( (k_n,n,N_n) : n =n_0,n_0+1, \dots )$ 
will be said to {\bf grow proportionally} if there are $\delta \in (0,1)$
and $\rho \in (0,1)$ so that 
\begin{equation}\label{eq:triple_prop}
      k_n/n \goto \rho, \quad n/N_n \goto \delta , \qquad n \goto \infty.
\end{equation}
We omit subscripts $n$ on $k$ and $N$ unless they are absolutely necessary.
\end{definition}

There are several significant prior results concerning
thresholds for face counts in the proportional-growth setting.

\subsubsection{Weak Thresholds}
Consider first the question whether the {\em typical} face survives
projection. 
 \bitem 
  \item {\it Simplex.} 
There is a function $\rho_{W}^+ : [0,1] \mapsto [0,1]$
with the following property. In the proportional growth setting with
$\rho < \rho_W^+(\delta)$, we have 
\[
            Ef_k(AT^{N-1}) = f_k(T^{N-1})(1 - o(1)), \quad 0 \leq k < \rho n, \qquad n \goto \infty;
\]
while if $\rho > \rho_W^+(\delta)$ we have that for some $\eps > 0$ and some 
sequence $(k_n)$ with $k_n < \rho n$,
\[
            Ef_k(AT^{N-1}) < f_k(T^{N-1})(1-\eps), \quad n \goto \infty.
\]
Informally, the  fraction
of faces lost:
\[
     (f_k(T^{N-1}) - E f_k(AT^{N-1}))/f_k(T^{N-1}) ,
\]
is either negligible or non-negligible depending on
which side of $\rho_W^+(\delta)$ the fraction $k/n$ sits.
In words, for $k_n$ somewhat below the threshold $n \cdot \rho_W^+(\delta)$
the typical $k_n$-face of the simplex survives projection
into $n$ dimensions; but for $k_n$ somewhat 
above the threshold this is no longer true.
\item {\it Cross-Polytope.} 
There is a function $\rho_{W}^\pm : [0,1] \mapsto [0,1]$
with the following property. In the proportional growth setting with
$\rho < \rho_W^\pm(\delta)$, we have
\[
            Ef_k(AC^{N}) = f_k(C^{N})(1 - o(1)), \quad 0 \leq k < \rho n, \quad n \goto \infty;
\]
while if $\rho > \rho_W^\pm(\delta)$ we have for some $\eps > 0$ and some 
sequence $(k_n)$ with $k_n < \rho n$,
\[
            Ef_k(AC^{N}) < f_k(C^{N}))(1-\eps), \quad n \goto \infty.
\]
Again, for $k_n$ somewhat below the threshold $n \cdot \rho_W^\pm(\delta)$
the typical $k_n$-face of the cross-polytope survives projection
into $n$ dimensions; but for some $k_n$
at or above the threshold this is no longer true.

\eitem

In view of these results, the square $ 0 \leq \delta, \rho \leq 1$
may be decorated with a {\it phase diagram}.
The two $\rho_W$-functions mark {\it phase transitions}; there are two phases
for the property  ``the projected polytope has {\it approximately} as many faces
as the original''. Below the transitions, the property holds
asymptotically for large $n$, while above the transitions.
the property fails asymptotically for large $n$.
Both transitions are depicted in Figure \ref{fig:thresholds_one},
which displays a {\it phase diagram} in $(\delta,\rho)$ plane.
In the region below these curves, typical faces are not lost, in the region
above those curves, typical faces are lost. To interpret these curves,
note that if $\delta = 1/2$ so we are lowering dimension by 50\%, and if $n$ is large,
then the typical $k$-face of the simplex survives, for $k/n \leq .5581$,
while the typical $k$-face of the cross-polytope survives, for $k/n \leq .3848$.

Vershik and Sporyshev  \cite{VerSpor} pioneered study of the proportional growth setting, 
and proved the existence of what we call here the weak threshold for the Simplex case.
The weak threshold function $\rho_W^+$ was introduced using our notation and carefully
studied by the authors in \cite{DoTa05_polytope}, where numerical
methods were developed for its calculation and display. The weak threshold 
for the cross-polytope $\rho_W^\pm$ was introduced 
in \cite{Do05_polytope}, calculated, and displayed.

\subsubsection{Strong Thresholds}
We now ask when  the difference between $f_k(AQ)$ and $f_k(Q)$ is small in absolute,
not relative, terms.
 \bitem 
  \item {\it Simplex.} There is a function $\rho_{S}^+ : [0,1] \mapsto [0,1]$
with the following property. In the proportional growth setting with
$\rho < \rho_S^+(\delta)$, we have
\[
            Ef_k(AT^{N-1}) = f_k(T^{N-1}) - o(1) , \qquad 0 \leq k < \rho n, \quad n \goto \infty.
\]
Thus, for $k$ below $n \cdot \rho_S^+(\delta)$
there are on average as many $k$-faces of the 
projected simplex as the original simplex.  On the other hand, if $\rho > \rho_S^+(\delta)$,
then there is a sequence $(k_n)$ with $k_n < n \rho$ along which
\[
   f_k(T^{N-1}) - E f_k(AT^{N-1}) \goto \infty.
\]

\item {\it Cross-Polytope.} There is a function $\rho_{S}^\pm: [0,1] \mapsto [0,1]$
with the following property. In the proportional growth setting with
$\rho < \rho_S^{\pm}(\delta)$, 
\[
           E f_k(AC^{N}) = f_k(C^{N}) - o(1) , \qquad 0 \leq k < \rho n, \quad n \goto \infty.
\]
Thus, for $k$ somewhat below $n \cdot \rho_S^\pm(\delta)$
there are on average just as
many $k$-faces of the projected cross-polytope 
as the standard cross-polytope. On the other hand, if $\rho > \rho_S^\pm(\delta)$,
then there is $k_n < n \rho$ with
\[
   f_k(C^N) - E f_k(AC^N) \goto \infty.
\]
\eitem

The function $\rho_S^+$ was introduced and  carefully
studied by the authors in \cite{DoTa05_polytope}, and numerical
methods were developed for its calculation and display. The threshold function
 $\rho_S^\pm$ was introduced in \cite{Do05_polytope}, calculated, and displayed.

These strong thresholds have another interpretation.
Consider the event  ``all low-dimensional faces survive projection'',
i.e.
\[
       \Omega(k,n,N) =  \{ f_\ell(AT^{N-1}) = f_\ell(T^{N-1}), \ell = 0, \dots, k\} .
\]
Simple arguments as in \cite{Do05_polytope,DoTa05_polytope}
show that  if $\rho < \rho_S^+(\delta)$,
the probability 
\[
P(\Omega(k_n,n,N_n) ) \goto 1, \qquad n \goto \infty.
\]
Hence, below the strong phase transition,
all low-dimensional faces survive projection.
Parallel arguments
can be made in the cross-polytope case.
Thus in the region where $k/n$ is below the corresponding $\rho_S$ function not only
are very few $k$ faces lost on average; actually, there is overwhelming probability
that no faces are lost. 

\begin{figure}[h]
\begin{center}
\begin{tabular}{c} 
\includegraphics[width=3in]{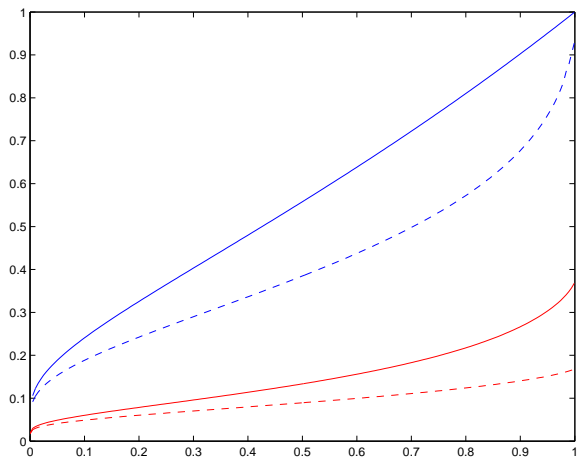}
\end{tabular}
\end{center}
\caption{Thresholds for $\delta\in(0,1)$ from top to bottom: $\rho_W^+$(blue - solid), $\rho_W^{\pm}$(blue - dashed), $\rho_S^+$(red - solid), and $\rho_S^{\pm}$(red - dashed).}
\label{fig:thresholds_one}
\end{figure}

These $\rho$-functions
are depicted in Figure \ref{fig:thresholds_one}. The strong 
thresholds $\rho_S^+$ and $\rho_S^{\pm}$
fall below the corresponding weak thresholds $\rho_W^+$, $\rho_W^{\pm}$;
indeed a property holding for every $k$-face
is less likely to hold than one holding for the typical $k$-face.
To interpret these curves,
note that if $\delta = 1/2$ so we are lowering dimension by 50\%,
then every $k$-face of the simplex survives, for $k \leq .1335$,
while every $k$-face of the cross-polytope survives, for $k \leq .0894$.

\subsection{Main Results: Proportional Growth Setting}

For applications,
the range where $\delta$ is small is very interesting;
it corresponds to:
\bitem
\item studying convex hulls of Gaussian point clouds where there are
many points relative to the number of dimensions -- Section \ref{int:Hull};
\item  recovering a sparse signal
from very few samples -- Section \ref{int:signal}; 
\item protecting against errors in digital transmission while sacrificing  
very little in the transmission rate -- Section \ref{int:ecc}. 
\eitem

Previous work by the authors \cite{DoTa05_polytope,Do05_polytope}
 considered the asymptotic
behavior of the several $\rho(\delta)$ 
functions just defined, and showed that $\rho(\delta) \geq c_\eps \log(1/\delta)^{-1-\eps}$ 
for each $\eps > 0$. Work by others \cite{CandesTao,Vershynin,LinialNovik} can be seen to imply that actually
$\rho(\delta) \geq c \log(1/\delta)^{-1}$. In this paper we determine the precise constants
in the asymptotic behavior as $\delta \goto 0$. These precise
constants are important in applications; they
can be used to plan how many samples to take in a digital
imaging system or how much transmission rate sacrifice
to make for a given error resistance. 

\begin{theorem}[Weak Threshold - Simplex]\label{thm:weakp}
\begin{equation} \label{eq:weakfrac}
\rhopW(\delta ) \sim  \left| 2\log(\delta)\right|^{-1}, \qquad \delta \goto 0.
\end{equation}
\end{theorem}

\begin{theorem}[Strong Threshold - Simplex]\label{thm:strongp}
\begin{equation}
\rhopS(\delta) \sim  \left| 2e\log(\delta 2\sqrt{\pi})\right|^{-1}, \qquad \delta \goto 0.
\end{equation}
\end{theorem}

Comparing these results:
\bitem\item Note the 
leading factor $e$ in $\rho_S$.  The
highest dimension $k$ where  the 
{\em vast majority} of $k$-faces survive projection is asymptotically $e$ 
times higher than the dimension where we can guarantee that {\em every} 
$k$-face survives.
\item An additional difference
is  the $2\sqrt{\pi}$ factor in the argument
of the logarithm.
\eitem

\begin{theorem}[Weak Threshold - Cross-polytope]\label{thm:weakpm}
\begin{equation} \label{eq:weakcrossasymp}
\rhopmW(\delta ) \sim  \left| 2\log(\delta)\right|^{-1}, \qquad \delta \goto 0.
\end{equation}
\end{theorem}

\begin{theorem}[Strong Threshold - Cross-polytope]\label{thm:strongpm}
\begin{equation}
\rhopmS(\delta) \sim \left| 2e\log(\delta \sqrt{\pi})\right|^{-1}, \qquad \delta \goto 0.
\end{equation}
\end{theorem}

Comparing the cross-polytope results to those for 
the simplex:
\bitem
\item Remarkably, to first order, the thresholds are the same for the 
simplex and cross-polytope.  This is surprising since at moderate values of $\delta$
the two functions are quite different; see Figure \ref{fig:thresholds_one}.
\item The bounds on strong 
thresholds agree, except for factors of $2$ in the argument of the logarithm.
\eitem

\begin{figure}[h]
\begin{center}
\begin{tabular}{c}
\includegraphics[width=3in]{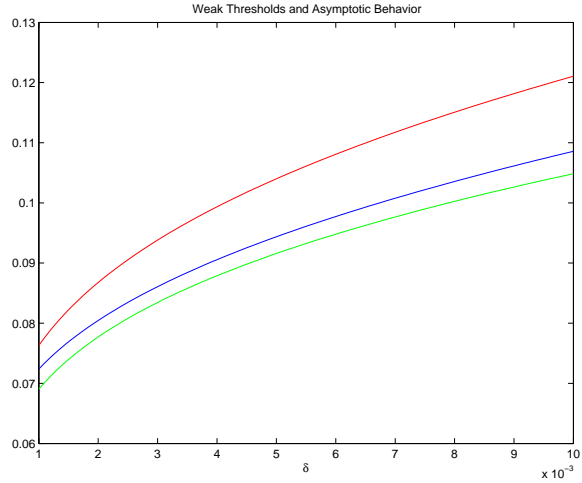}
\end{tabular}
\end{center}
\caption{Weak thresholds $\rho_W^+$ (red), $\rho_W^{\pm}$ (green), and their asymptotic behavior, $\left| 2\log(\delta)\right|^{-1}$ (blue), from Theorems \ref{thm:weakp} and \ref{thm:weakpm}, $\delta\in [10^{-3},10^{-2}]$.}
\label{fig:thresholds_asym_weak}
\end{figure}

\begin{figure}[h]
\begin{center}
\begin{tabular}{c}
\includegraphics[width=3in]{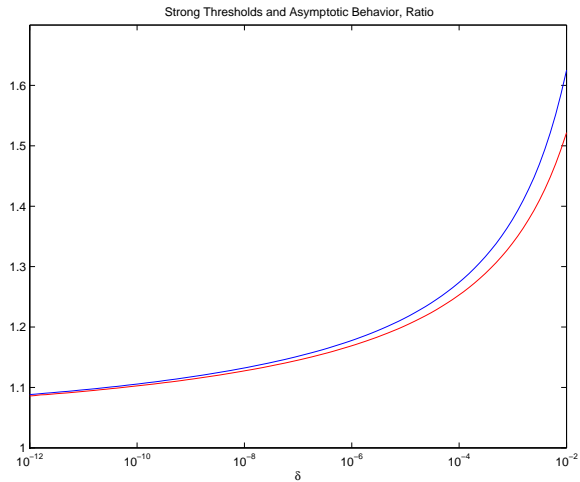} 
\end{tabular}
\end{center}
\caption{Ratio of the strong thresholds and their asymptotic behaviors, $\rho_S^+(\delta)$ in blue and $\rho_S^{\pm}(\delta)$ in red.  The asymptotic formulae approach slowly from above by the factors shown.}
\label{fig:thresholds_asym_strong}
\end{figure}

The weak-threshold asymptotic behavior (\ref{thm:weakp}) 
and (\ref{eq:weakcrossasymp})
closely matches $\rho_W^+$ and $\rho_W^{\pm}$ 
for modest values of $\delta$ -- see Figure \ref{fig:thresholds_asym_weak}. 
The strong-threshold asymptotic 
behavior, on the other hand,  slowly approaches $\rho_S^+$ and $\rho_S^{\pm}$ 
from above -- see Figure \ref{fig:thresholds_asym_strong}.

\subsection{Beyond Proportional Growth}\label{int:beyond_prop}

Having considered the Vershik-Sporyshev proportional 
growth scenario, we now generalize
to the case where $N$ can be dramatically larger than $n$.  
This  is important for applications where we want to
sample very few projections of a high dimensional object.
\cite{CS} exhibits stylized imaging problems where an $N$-pixel
image can be reconstructed by asking $n = O(N^a)$ questions, $a < 1$.
This of course lies outside the reach of proportional growth and
is dramatically smaller than $N$, 
underscoring the potential interest of the `how many questions' problem
of Section \ref{int:signal} where the number of questions $n\ll N$.

We would naively hope that the same threshold functions $\rho()$
``work'' even outside the proportional growth setting.
That is, in a setting where $n/N_n \goto 0$ , we would hope to get the `right answer'
for the behavior of face counts by simply `plugging in' a {\it varying $\delta$}
$= \delta_n = n/N_n \goto 0$ into the appropriate $\rho$-function. 
Happily, such naive hopes go unpunished. 

We say that $N$ {\it grows subexponentially} relative to $n$ if
\begin{equation} \label{defsubexpon}
     N_n/n \goto  \infty, \qquad  \frac{\log(N_n)}{n}  \goto 0, \qquad n \goto \infty.
\end{equation}

\begin{theorem}\label{thm:5}
Consider a sequence of problem sizes $(n,N_n)$ where $N_n$ grows subexponentially
relative to $n$. Let $\rho$ be one of the four functions $\rho_W^+$, $\rho_W^\pm$,
$\rho_S^+$,$\rho_S^\pm$.  Fix $\eps > 0$ and
consider  a sequence $(k_n)$
obeying $k_n/n  < \rho(n/N_n)(1-\eps)$ for $n > n_0$.  Then the same statement
that was made for that  $\rho$ in the proportional growth scenario
holds in this non-proportional growth scenario. 

Thus, for example,
$k_n < (1-\eps) \rho_W^+(n/N_n) \cdot n$ for $n = n_0, n_0+1, \dots$  implies
\[
            Ef_k(AT^{N-1}) = f_k(T^{N-1})(1 - o(1)), \quad n \goto \infty;
\]
similarly,
$k_n < (1-\eps) \rho_S^\pm(n/N_n) \cdot n $ for $n = n_0, n_0+1, \dots$  implies
that with overwhelming probability for large
$n$,
\[
            f_\ell(AC^{N}) = f_\ell(C^{N}), \quad \ell=0,\dots,k-1 .
\]
\end{theorem}

In short, the limit relations of Theorems \ref{thm:weakp}-\ref{thm:strongpm} 
are useful both in proportional and non-proportional growth settings.

\subsection{Contents}

Our paper proves Theorems 
\ref{thm:weakp}-\ref{thm:5}.  The development is organized 
as a branching tree, with initial sections mapping out the main concepts,
propositions and lemmas, and later sections dealing with  
detailed estimates and proofs.

Section 2 introduces the underlying machinery of face
counting and an analytic approach to studying asymptotic behavior.
Our starting point is  a beautiful and essential formula
for the expected number of faces of randomly projected
polytopes, due to Affentranger and Schneider and
Vershik and Sporyshev; it involves three factors, representing 
contributions from combinatorial aspects, 
from external angles and from internal angles.
We focus on the exponential
growth and/or decay of the factors by defining associated 
exponents $\Psi(\delta,\rho)$;
we work as if these factors behave exactly as $\exp\{ N \cdot \Psi \}$.
Each $\rho()$ function is defined as the smallest 
root $0=\Psi(\delta,\rho(\delta))$ of an associated exponent function $\Psi$, 
viewed as a function of $\rho$ with $\delta$ fixed. 

Section 3 gives the proofs for the lower bound half of 
Theorems \ref{thm:weakp}-\ref{thm:strongpm}.
The proofs are simple consequences of
 the asymptotic behavior of the net exponents as 
a function of $\delta$ and $\rho$ in the regime where $\delta \goto 0$.

Section 4 develops the basic asymptotic analysis
of the net exponents.
The exponents in question explicitly
involve tail probabilities of the Gaussian distribution;
our asymptotic analysis exploits detailed estimates for
the Mills' ratio of the standard normal density.

Section 5 turns to the proof of Theorem \ref{thm:5}, going outside the proportional
growth setting. Here we have to make careful estimates
of the errors incurred by treating the pieces in
the Affentranger-Schneider-Vershik-Sporyshev formula
as if they grow exactly like $\exp\{ N \cdot \Psi \}$. We refine our analysis
associated with Mills' ratio, getting remainder estimates assuming
$N_n$ is subexponential in  $n$.

Sections 2-5 are preoccupied largely with proving only half 
of Theorems \ref{thm:weakp}-\ref{thm:strongpm}; namely the bounds $\rho(\delta) \geq c_1 /\log(c_2/\delta)$.
Section 6 gives the arguments establishing inequalities in the 
other direction, in the process completing the proofs of Theorems \ref{thm:weakp}-\ref{thm:strongpm}.

Section 7 shows how our face-counting
results generate the applications mentioned in Section 1.1.
It also presents empirical results showing that
our asymptotic results work at moderate sample sizes,
and translates our asymptotic results into finite-sample
bounds. It  also considers extensions of this work,
and compares our results with other recent work.

\section{Definitions of $\rho_S^+$,$\rho_W^+$,$\rho_S^\pm$,$\rho_W^\pm$}\label{sec:geometry}

The various $\rho$ quantities referred
to in Theorems \ref{thm:weakp}-\ref{thm:strongpm} have so far 
been discussed behaviorally,
by their role in locating or bounding phase transitions in face counts.
In this section, we review an analytic definition
for these quantities given in  
\cite{Do05_polytope,DoTa05_polytope}.
The definition
unfortunately requires a considerable amount of machinery
associated with convex integral geometry.
Equipped with such machinery,
the claims made by Theorems \ref{thm:weakp}-\ref{thm:strongpm}
can be translated into sharply-defined 
questions about the leading-order asymptotics
of certain exponents. Sections 3 and 6 answers
those questions.

\subsection{Expected Face Counts of Projected Polytopes}
 
Let $Q$ be a polytope in $R^N$ and 
$A:R^N\mapsto R^n$ a random ortho-projection,
uniformly distributed on the Grassmann manifold
of all such projectors.
Affentranger  and Schneider \cite{AffenSchnei}
developed a useful identity
for the expected number of faces 
of $AQ$ \cite{AffenSchnei}:
\begin{equation}\label{eq:projection}
E f_k(AQ)=f_k(Q)-2{\sum_{\ell}}' \sum_{F\in\cF_k(Q)}\sum_{G\in\cF_{\ell}(Q)} 
\beta(F,G)\alpha(G,Q);
\end{equation}
here  $\cF_k(Q)$ denotes the set of $k$-faces of $Q$, each $F$ is a subface of $G$,
and $\sum'$ denotes the sum over $\ell = n+1, n+3, \dots$; $\ell < N$.
We are intensely interested in the discrepancy between the 
expected number of faces of the projected polytope $AQ$ and the 
necessarily larger number of faces of the original polytope $Q$;
i.e. in knowing {\it on average, how many faces are lost
in the projection from $\bR^N$ to $\bR^n$}.
The discrepancy in question is
\begin{equation}\label{eq:lost_faces}
\Delta(k,n,N;Q):=f_k(Q)-E f_k(AQ)=2{\sum_{\ell}}'  \sum_{F\in\cF_k(Q)}\sum_{G\in\cF_{\ell}(Q)} 
\beta(F,G)\alpha(G,Q).
\end{equation}
Here the sum covers the external angles between the 
original polytope $Q$ and its subfaces $G$, $\alpha(G,Q)$, multiplied 
by the sum of all internal angles between each particular subface $G$ and 
its faces $F$, $\beta(F,G)$.  For definitions of these angles see 
eg. Gr\"unbaum \cite[Chapter 14]{Gru67}, or Matousek \cite{Mat02}. 

\subsection{Analytic Definition of $\rho_S^+$,$\rho_S^\pm$}

In the remainder of the paper we are always interested
in just two choices of $Q$: the simplex, $Q = T^{N-1}$, and the
cross-polytope, $Q=C^N$. Various quantities associated with
the simplex case will be labeled with superscript $+$
(as the interior of the standard simplex consists of positive vectors)
and objects associated with 
the cross-polytope case will be labelled with 
superscript $\pm$ (as the standard cross-polytope contains vectors with
entries of both signs.) We frequently
use $\star$ as a superscript  in a statement which 
concerns either case, implying two different
statements, with obvious substitutions.

In the introduction, the functions $\rho_S^+$ and $\rho_S^\pm$ 
were partially characterized by 
the claim that, for $(k_n,n,N_n)$ growing proportionally
and limit ratios $(k_n/n , n/N_n) \goto (\rho,\delta)$ with
$\rho < \rho_S^\star(\delta)$, then 
\begin{equation}\label{discrep20}
\Delta(k_n,n,N_n;Q) \goto 0, \quad n \goto \infty.
\end{equation}
[Note: To make sure the reader follows our convention for $\star$,
the previous sentence is actually two sentences, one for the
symbol binding $(\star,Q) = (+, T^{N-1})$ and one for the
symbol binding $(\star,Q) = (\pm,C^N)$.]
It was also stated that if $\rho > \rho_S^\star$ then
for some sequence $(k_n)$ obeying $k_n < \rho n$, and some $\eps > 0$,
\begin{equation} \label{discreppos}
   \liminf_{n \goto \infty} \Delta(k_n,n,N_n;Q) \geq \eps > 0.
\end{equation}

The papers \cite{Do05_polytope,DoTa05_polytope} actually defined
 $\rho_S^\star(\delta)$
  with the following stronger property:
 if $\rho < \rho_S^\star(\delta)$, then,  in the proportional growth setting (\ref{eq:triple_prop})
for some $\epsilon>0$ and $n > n_0(\eps,\rho)$, we have
\begin{equation}\label{eq:decay_def_Q}
\Delta(k_n,n,N_n; Q)\le N_n\exp(-N_n\epsilon).
\end{equation}
Those papers implied/stated without proof that
if $\rho > \rho_S^\star$ then
for some sequence $k_n < \rho n$, some $\eps > 0$,
and $n_0$ we have
\begin{equation}\label{eq:growth_def_Q}
   \Delta(k_n,n,N_n;Q) \geq \exp(N \eps), \quad n > n_0.
\end{equation}

While conceptually, both (\ref{discrep20}) and (\ref{discreppos})
are equally important parts of the picture, in practice (\ref{discrep20})
is the more useful/surprising. Hence in Sections 3-5 of this paper
we focus on supporting assertions like (\ref{discrep20})
and (\ref{eq:decay_def_Q}) showing that the discrepancy is small,
rather than than assertions like (\ref{discreppos}) and (\ref{eq:growth_def_Q})
showing that the discrepancy is large. Section 6
will return to (\ref{discreppos}) and (\ref{eq:growth_def_Q}).



The analysis supporting the bound
(\ref{eq:decay_def_Q}) for the unit simplex and cross-polytope 
went by first rewriting (\ref{eq:lost_faces}) as 
a sum of contributions due to faces of different dimensions:
\[
\Delta(k_n,n,N_n;Q)={\sum_{\ell}}'  D_{\ell,n}^\star,
\]
where 
\[
D_{\ell,n}^\star :=2\sum_{F\in\cF_k(Q)}\sum_{G\in\cF_{\ell}(Q)} \beta(F,G)\alpha(G,Q).
\]
The papers
\cite{DoTa05_polytope,Do05_polytope} defined
functions $\Psi_{net}^\star(\nu,\gamma)$  for $\star  \in \{ +, \pm \}$
associated with our two choices for $Q$; these can be used to
bound $D_{\ell,n}^\star$ as follows. 

Put $\nu_{\ell,n} = \ell/N_n$
and $\gamma_{\ell,n} = k_n/\ell$, and note that $\nu_{\ell,n} \in [\delta,1]$
and $\gamma_{\ell,n} \in [0,\rho]$ over the relevant range 
$\ell = n+1, n+2, \dots $; $\ell < N$.
In the proportional growth setting \cite{DoTa05_polytope,Do05_polytope} showed that,
 for each $\eps > 0$, there is $n_0(\eps; \delta,\rho)$ so that
\begin{equation} \label{psinetprop}
     N_n^{-1} \log(D_{\ell,n}^\star) \leq \Psi_{net}^\star( \nu_{\ell,n}, \gamma_{\ell,n} ) + 3 \epsilon, 
     \quad \ell = n+1, n+3,\ldots, 
     \quad n \geq n_0,
 \end{equation}
Since our focus is the condition (\ref{eq:decay_def_Q}), 
we of course are interested in conditions guaranteeing
that the right side is negative, uniformly over the {\it admissible domain} of $(\nu,\gamma)$ 
pairs obeying $\nu \geq \delta$, $\gamma \leq \rho$. 

\begin{definition}
The {\bf maximal operator} $M[]$  associated to the family of
rectangles where
 $\nu \in [ \delta,1]$, $\gamma \in [0,\rho]$
 takes a function $\psi(\nu,\gamma)$,
and delivers the {\bf maximal function} $M[\psi](\delta,\rho)$ defined
by
\[
  M[\psi](\delta,\rho) = \sup \{ \psi(\nu,\gamma) : \nu  \in [\delta,1], \gamma \in [0,\rho] \}.
\]
\end{definition}

Applying this operator to each $\Psi_{net}^\star$
yields two maximal functions,
$M[\Psi_{net}^+]$ and $M[\Psi_{net}^\pm]$,
to be studied extensively below. 
Finally we can give
an analytic definition for the key quantities in Theorems \ref{thm:strongp} 
and \ref{thm:strongpm}:
 
 \begin{definition} \label{def:rhos}
For $\star \in \{ +, \pm\}$, define the
strong phase transition $\rho_S^\star(\delta)$ as the `first' zero of $M[\Psi_{net}^\star]$:
\[
      \rho_S^\star(\delta) = \inf \{ \rho : M[\Psi_{net}^\star](\delta,\rho) = 0, \rho \in [0,1] \} .
\]
\end{definition}

\begin{figure}[h]
\begin{center}
\begin{tabular}{c}
\includegraphics[width=3in]{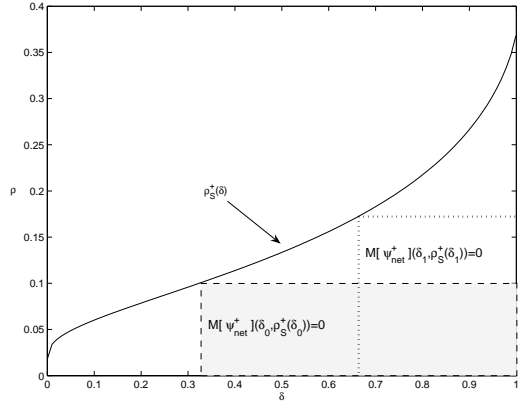}
\end{tabular}
\end{center}
\caption{Throughout the shaded region $\psinetp\le 0$.  Therefore the 
maximal function is $\le 0$ at the upper left corner 
$(\delta_0,\rho_S^+(\delta_0))$; in fact, $M[\psinetp](\delta_0,\rho)<0$ 
for $\rho<\rho_S^+(\delta_0)$, and 
$M[\psinetp](\delta_0,\rho_S^+(\delta_0))=0$; this is the ``first'' zero 
of $M[]$.  The family of such first zeros of $M[\psinetp](\delta,\rho)$ 
define $\rho_S^+(\delta)$.}
\label{fig:figure1}
\end{figure}

Definition \ref{def:rhos} is depicted in Figure \ref{fig:figure1}.

Several properties of the $\rho_S^\star$ are known
from \cite{Do05_polytope,DoTa05_polytope}. These functions are strictly positive
on $[0,1]$, strictly increasing, with limit $0$ as $\delta$ tends to 
$0$ and limits $\approx .3679$ and $.1685$ ($\star=+,\pm$ respectively) 
as $\delta$ tends to $1$.

The functions $M[\Psi_{net}^\star]$ are continuous. It follows that for
$\rho < \rho_S^\star(\delta)$,
\[
   M[\Psi_{net}^\star](\delta,\rho) < 0.
\]
Setting $\eps =  |M[\Psi_{net}^\star](\delta,\rho)|/4 $,
\begin{equation}\label{eq:prepsinet}
N^{-1}\log (D_{\ell,n}^\star )\le -\epsilon, \qquad \ell = n+1, n+3, \dots . 
\end{equation}
The result (\ref{eq:decay_def_Q}) follows. 

\subsection{Analytic Definition of $\rho_W^+$,$\rho_W^\pm$}

The papers \cite{DoTa05_polytope} and \cite{Do05_polytope}
also defined phase transitions 
 $\rho_W^+$ and $\rho_W^\pm$. Conceptually, these quantities
are defined by the notion that, for $(k_n,n,N_n)$ growing proportionally
with limit ratios $(k_n/n , n/N_n) \goto (\rho,\delta)$, then 
if $\rho < \rho_W^\star(\delta)$, the relative discrepancy
is negligible
\begin{equation} \label{eq:relneg}
\Delta(k,n,N;Q)/f_k(Q) = o(1), 
\quad k = 0, \dots, \lfloor \rho n \rfloor, \quad n \goto \infty, 
\end{equation}
while for $\rho > \rho_W^\star(\delta)$ the relative discrepancy
can be substantial; for some sequence $(k_n)$ obeying $k_n < \rho n$
and some $\eps > 0$ and $n_0$,
\begin{equation} \label{eq:relpos}
 \Delta(k,n,N;Q)/f_k(Q) \geq \eps > 0, \quad n > n_0. 
\end{equation}
Again while conceptually both  (\ref{eq:relneg}) and
(\ref{eq:relpos}) are equally important, practically
speaking the former is more useful/significant
than the latter, which mainly serves to show that
we cannot substantially improve on (\ref{eq:relneg}).
We will focus on (\ref{eq:relneg}) in Sections 3-5
and then return to discussion of (\ref{eq:relpos})
in Section 6.

Define 
\begin{equation} \label{eq:defacep}
  \Psi_{face}^+(\nu,\gamma) = H(\nu\gamma),
\end{equation}
with $H(\cdot)$ the type-$e$ Shannon entropy (\ref{ShannonDef}),
so that  under proportional growth
\[
        N^{-1} \log f_k(T^{N-1})  \goto  \Psi_{face}^+(\rho,\delta) .
\]
Also, put
\begin{equation} \label{eq:defacepm}
  \Psi_{face}^\pm(\nu,\gamma) = H(\nu\gamma) + \nu\gamma \log_e(2),
\end{equation}
so that under proportional growth
\[
        N^{-1} \log f_k(C^{N})  \goto  \Psi_{face}^\pm(\rho,\delta) .
\]
\begin{definition}
For $\star \in \{ +, \pm\}$, define 
$\rho_W^\star(\delta)$ as the `first' zero of 
$M[\Psi_{net}^\star-\Psi_{face}^\star]$:
\[
      \rho_W^\star(\delta) = \inf \{ \rho : M[\Psi_{net}^\star-\Psi_{face}^\star](\delta,\rho) = 0, \rho \in [0,1] \} .
\]
\end{definition}

If $\rho < \rho_W^\star(\delta)$, in the proportional growth setting, then
for some $\epsilon>0$ we have
\begin{equation}\label{eq:decay_def}
\Delta(k,n,N; Q)/ f_k(Q) \le N\exp(-N\epsilon), \qquad n > n_0.
\end{equation}
This establishes (\ref{eq:relneg}).

\subsection{Simplex exponent $\psinetp$}\label{subsec:lfaces_simplex}

We now give more details about the exponent
$\psinetp$ for the Simplex $T^{N-1}$.  We begin with 
observations by Affentranger and Schneider \cite{AffenSchnei} 
and Vershik and Sporyshev \cite{VerSpor}, that:
\bitem
\item There are ${N \choose k+1 }$ $k$-faces of $T^{N-1}$.
\item For $\ell > k$, there are ${N-k-1 \choose \ell-k }$ $\ell$-faces 
of $T^{N-1}$ containing a given $k$-face of $T^{N-1}$.
\item The faces of $T^{N-1}$ are all simplices, and the internal angle
$\beta(F,G) = \beta(T^k,T^\ell)$,
 where $T^d$ denotes the standard $d$-simplex.
 \eitem
Thus, for $\ell = n+1, n+3,\ldots$ we can write
\begin{eqnarray*}
      D_{\ell,n} ^+   &=&  2   {N \choose k+1 } {N-k-1 \choose \ell-k} \cdot \beta(T^k,T^\ell)  \cdot\alpha(T^\ell,T^{N-1}) \\
        &=& C^+_{\ell,n}\cdot \beta(T^k,T^\ell)\cdot  \alpha(T^\ell,T^{N-1}),
\end{eqnarray*}
with $C^+_{\ell,n}$ denoting the combinatorial prefactor.

Each of the factors in this product has either exponential growth or decay.
We will soon define associated 
exponents $\psicomp, \psiintp$, and $\psiextp$ so that, for any $\epsilon>0$ 
and $n > n_0(\delta,\rho)$,
\begin{equation} \label{BoundComb}
N^{-1} \log(C^+_{\ell,n}) \leq   \psicomp( \nu_{\ell,n} , \gamma_{\ell,n} ) 
+ \eps, 
\end{equation}
\begin{equation} \label{BoundInt}
N^{-1} \log( \beta(T^k,T^{\ell})) \leq  - \psiintp( \nu_{\ell,n} , \gamma_{\ell,n} )  + \eps,
\end{equation}
and
\begin{equation} \label{BoundExt}
N^{-1} \log( \alpha(T^{\ell},T^{N-1})) \leq  - \psiextp( \nu_{\ell,n} )  + \eps,
\end{equation}
uniformly in $\ell = n+1, n+3,\ldots$; $\ell < N$.  

The exponents were introduced in \cite{DoTa05_polytope},
which showed (\ref{BoundComb})-(\ref{BoundExt});
we repeat the definitions, although the reader
should not expect much insight 
at this point. The definitions are restated in Section 
\ref{sec:exponents}; equations (\ref{eq:psicom-restated}), 
(\ref{eq:psiint-restated}), (\ref{eq:xnu-restated}), and (\ref{eq:psiextpm-restated}),
where further details emerge. The combinatorial exponent
involves the base-$e$ Shannon entropy:
\begin{equation} \label{ShannonDef}
   H(p) = p \log(1/p) + (1-p) \log(1/(1-p)).
\end{equation}
Thus,
\begin{equation}\label{eq:psicom_related}
\psicomp(\nu,\gamma):=H(\nu)+\nu H(\gamma) .
\end{equation}
The internal exponent is 
\begin{equation}\label{eq:psiint}
\psiintp(\nu,\gamma):=\nu(1-\gamma)\left[
\log(\ygam/\gamma)+\frac{1}{2}\log(2\pi)+\frac{\gamma-1}{2\gamma}\ygam^2 .
\right]
\end{equation}
Here
$\ygam$ 
is defined implicitly by
\begin{equation}\label{eq:ygamma}
\frac{1-\gamma}{\gamma}\ygam=s_{\gamma}\quad\mbox{with}\quad
R(s_\gamma) = 1-\gamma,
\end{equation}
where
\begin{equation} \label{eq:defMillsR}
R(s):=s e^{s^2/2}\int_{s}^{\infty}e^{-y^2/2}dy.
\end{equation}
The function $R(s)$ is closely related to
a fundamental tool 
for studying tail probabilities of the
standard Normal distribution \-- the
so-called {\it Mills' ratio} of the Normal distribution,  \cite[Sec. 5.37]{KendallVol1}
about which more will be said in later sections.
The fact that $R(s) \goto 1$ as $s \goto \infty$
signifies that the tail probability under the normal
distribution is asymptotic to $s^{-1}$ times the normal
density. Details of this approximation will be
crucial for our work here.
Finally, the external exponent is:
\begin{equation}\label{eq:xnu}
\psiextp(\nu):=\nu \xnu^2-(1-\nu)\log Q(\xnu),
\end{equation}
with $\xnu$ the solution of
\begin{equation} \label{xnudef}
\frac{2xQ(x)}{q(x)}+1-\nu^{-1}=0 ;
\end{equation}
here $q(x):=\pi^{-1/2}e^{-x^2}$ and $Q(x)=\int_{-\infty}^x q(y)dy$.
Note that $Q()$ is the normal distribution with mean zero and standard 
deviation $1/\sqrt{2}$, and so $\xnu$ is again associated with the 
relationship between tail probabilities and density.
This definition seems at first very similar
to the definition of the internal angle;
however, note that $Q(x) \goto 1$ as $x \goto \infty$,
while $q(x) \goto 0$ rapidly. This difference is reflected in the behavior of the
$x_\nu$ as a function of $\nu$ which is very different than
the behavior of $\ygam$ as a function of $\gamma$.

These $\Psi$-functions are all smooth functions
of their arguments. For details on these exponents, see either the
original source \cite{DoTa05_polytope}, where graphical displays are provided,
or Section \ref{sec:exponents} below.

It follows from (\ref{BoundComb})-(\ref{BoundExt}) that 
for $\ell = n+1,n+3,  \dots$, 
\[
N^{-1}\log (D_{\ell,n}^+)\le 
\psicomp( \nu_{\ell,n} , \gamma_{\ell,n} )-\psiintp( \nu_{\ell,n} , \gamma_{\ell,n} )-\psiextp( \nu_{\ell,n} , \gamma_{\ell,n} )+3\epsilon.
\]
Defining now
\[
\psinetp(\nu,\gamma) :=\psicomp(\nu,\gamma)-\psiintp(\nu,\gamma)-\psiextp(\nu)
\]
provides us the desired property (\ref{psinetprop}) referred to earlier, in the simplex case.
Graphs were presented in \cite{DoTa05_polytope}  showing the exponent's
behavior for $\gamma=.5555$ over the range $\nu \in (.5555,1]$. 
Software is available to make similar graphs for other parameter choices.

\subsection{ $\psinetpm$, Cross-polytope case}

B\"or\"oczky and Henk \cite{BoroHenk} previously studied the expected number 
of faces for the randomly projected cross-polytope $C^N$, and although the 
analysis is quite different, we utilize a number of their observations.

\bitem
\item There are $2^{k+1}{N \choose k+1 }$ $k$-faces of $C^N$.
\item For $\ell > k$, there are $2^{\ell-k} {N-k-1 \choose \ell-k }$ 
$\ell$-faces of $C^N$ containing a given $k$-face of $C^N$.
\item The faces of $C^N$ are all simplices, and the internal angle
$\beta(F,G) = \beta(T^k,T^\ell)$.
\item The external angle $\alpha(G^{\ell},C^N)$ is the same for all 
$\ell$-faces of $C^N$,
 the closed form expression of which was originally 
given in \cite{BoroHenk}.  A version written in our notation 
was developed in
\cite{Do05_polytope}, and is spelled out below in
(\ref{eq:extp_def}).
\eitem
Thus, for $\ell = n+1, n+3,\ldots$ we can write
\begin{eqnarray*}
      D_{\ell,n}^\pm    &=&  2 \cdot 2^{\ell}\cdot {N \choose k+1 } {N-k-1 \choose \ell-k} \beta(T^k,T^\ell) \alpha(F^\ell,C^N) \\
        &=& C^{\pm}_{\ell,n} \cdot \beta(T^k,T^\ell) \cdot \alpha(F^\ell,C^N),
\end{eqnarray*}
with $C^{\pm}_{\ell,n}$ the combinatorial prefactor.

The factors in this product again have either exponential growth or decay.
We will soon define associated 
exponents $\psicompm, \psiintpm$, and $\psiextpm$ so that, for any $\epsilon>0$ 
and $n > n_0(\epsilon)$,
\begin{equation} \label{BoundCombPM}
N^{-1} \log(C^\pm_{\ell,n}) \leq   \psicompm( \nu_{\ell,n} , \gamma_{\ell,n} ) 
+ \eps, 
\end{equation}
\begin{equation} \label{BoundIntPM}
N^{-1} \log( \beta(T^k,T^{\ell})) \leq  - \psiintpm( \nu_{\ell,n} , \gamma_{\ell,n} )  + \eps,
\end{equation}
and
\begin{equation} \label{BoundExtPM}
N^{-1} \log( \alpha(F^{\ell},C^{N})) \leq  - \psiextpm( \nu_{\ell,n} )  + \eps,
\end{equation}
uniformly in $ \ell = n+1, n+3,\ldots$; $\ell < N$.  It follows that for $n > n_0$,
\[
N^{-1}\log (D_{\ell,n}^\pm)\le 
\psicompm( \nu_{\ell,n} , \gamma_{\ell,n} )-\psiintpm( \nu_{\ell,n} , \gamma_{\ell,n} )-\psiextpm( \nu_{\ell,n} , \gamma_{\ell,n} )+3\epsilon.
\]

The exponents were introduced in \cite{Do05_polytope},
which showed (\ref{BoundCombPM})-(\ref{BoundExtPM});
we rehearse the definitions, admitting
they yield little insight 
at this point. The definitions are restated in Section 
\ref{sec:exponents}; 
where further information can be obtained. 
The combinatorial exponent again
involves the base-$e$ Shannon entropy $H$:
\begin{equation}\label{eq:psicomp_related} 
\psicompm(\nu,\gamma):=H(\nu)+\nu H(\gamma) + \nu \log_e(2);
\end{equation}
thus $\psicompm = \psicomp + \nu \log_e(2)$.
The internal exponent is actually the same as 
in the simplex case: $\psiintp = \psiintpm$.
Finally, the external exponent is:
\begin{equation}\label{eq:psiextpm}
\psiextpm(\nu):=\nu \ynu^2-(1-\nu)\log G(\ynu),
\end{equation}
with $\ynu$ the solution to
\begin{equation} \label{eq:psiextpmb}
\frac{2yG(y)}{g(y)}+1-\nu^{-1}=0,
\end{equation}
and $g(y):=2\pi^{-1/2}e^{-y^2}$,  $G(y)= \mbox{erf}(y)=\int_{0}^y g(w)dw$.  
$G()$ is the Error function, also called the 
Half-Normal distribution $HN(0,\frac{1}{2})$.
Again the $\Psi$'s are smooth functions of their arguments. 

Defining now
\[
\psinetpm(\nu,\gamma) :=\psicompm(\nu,\gamma)-\psiintpm(\nu,\gamma)-\psiextpm(\nu)
\]
provides us, in the cross-polytope case, 
the property (\ref{psinetprop}) 
referred to earlier.
In \cite{Do05_polytope}  it was shown that this 
is a well-defined and in fact nicely behaved quantity
as a function of $\gamma$ for each fixed $\nu$. Graphs in  \cite{Do05_polytope}
portray its behavior over the range $\nu \in (.5555,1]$ for 
$\delta =.5555$; software is available
to compute similar graphs as other values for $\nu$.

\section{Asymptotics of $\Psi_{net}^\star$ as $\delta \goto 0$}
\setcounter{equation}{0}
\setcounter{table}{0}
\setcounter{figure}{0}

We now turn to the asymptotics at the heart of 
Theorems \ref{thm:weakp}-\ref{thm:strongpm}.
As indicated earlier, in Sections 3-5 we focus on
establishing lower bounds on $\rho$-functions,
practically most `important' or
`surprising' part of our results.

We introduce a parametrized 
family of simple comparison functions 
$\rstar$ of the form 
$|\c\log(c\cdot\delta)|^{-1}$ and control the behavior of $\rho^{\star}(\delta)$
by studying the maximal functions along the trajectories $(\delta,\rstar)$ 
as $\delta\goto 0$.  The central point will be that for $\c>\c_0$, 
each associated maximal functions is asymptotically negative along the 
trajectory $(\delta,\rstar)$. This forces $\rho^\star(\delta) > \rstar$.

We will also glean insights useful for establishing upper bounds on
$\rho$-functions.  It will emerge that fixing $\c>\c_0$ 
defines a trajectory along which
the net exponents are asymptotically positive and that fixing $\c = \c_0$
defines a trajectory such that 
the difference between net and face exponents
is vanishing; it will be explained in Section 6 how this implies
the upper bound half of Theorems \ref{thm:weakp}-\ref{thm:strongpm}.

It is convenient to develop the results in a permuted order.

\subsection{Theorem \ref{thm:strongp}}

Fix $\cs>  2e$ and define $\rps:=\rpsc:= |\cs\log(\delta 2\sqrt{\pi})|^{-1}$.
In what follows, $\tau$ is always held fixed throughout an argument,
while $\delta$ is sent towards 0.

We intend to show that   there is $\delta_S^+ = \delta_S^+(\cs) > 0$ so that
\begin{equation}\label{Mnegative}
  M[\Psi_{net}^+] ( \delta , \rps ) < 0, \quad 0 < \delta < \delta_S^+.
\end{equation}
This establishes the lower-bound half of Theorem \ref{thm:strongp},
i.e. that $\rho_S^+(\delta) \ge |2e \log(\delta 2 \sqrt{\pi})|^{-1} \cdot (1 + o(1))$.
The other half  of Theorem \ref{thm:strongp} can be inferred from the fact that 
if we instead have $\cs<2e$ there is $\delta_0 = \delta_0(\cs) > 0$ with
\[
  \Psi_{net}^+ ( \delta , \rps ) > 0, \quad 0 < \delta < \delta_0.
\]
See further discussion in Section 6 below.

We start the proof of (\ref{Mnegative})
by observing that
the maximal operator $M[]$ becomes `transparent'
in the limit $\delta \goto 0$ if we stay along the
trajectory $(\delta,\rpsc)$. Corollary \ref{cor:psinetp} below shows that,
if $ \cs > 2e $,
for some $\delta_1=\delta_1(\cs) > 0$,
\begin{equation} \label{eq:EvalMaximalStrongSimp}
    \Psi_{net}^+(\delta,\rps) = M[\Psi_{net}^+] (\delta,\rps),
    \quad 0 < \delta < \delta_1.
\end{equation}

The following limiting behavior of the individual exponents as 
$\delta\rightarrow 0$ and/or $\rho\rightarrow 0$ 
will be derived in Section \ref{sec:exponents}, see (\ref{eq:psicomp}), (\ref{eq:psiint_asym}),
(\ref{eq:psiextp}):
\begin{eqnarray}
\psicomp(\delta,\rho) & = & \delta\left[\log(1/\delta)+1
+{\cal O}(\delta \vee \rho\log\rho)\right], \quad\max(\rho,\delta)\rightarrow 0, \label{eq:33a}  \\
\psiintp(\delta,\rho) & = & -\frac{1}{2}\delta \left[\log\rho+\log(e/2\pi)  + {\cal O}(\rho\log\rho)\right], \quad\rho\rightarrow 0, \label{eq:33b} \\
\psiextp(\delta) & = & \delta\left[\log \zpd -\frac{1}{2}\log\log \zpd + 1 
+{\cal O}\left(\frac{\log\log \zpd}{\log \zpd}\right)\right], \quad\delta\rightarrow 0;
\label{eq:first_asym}
\end{eqnarray}
here $\zpd:=(\delta 2\sqrt{\pi})^{-1}$, ${\cal O}(x)$ denotes a term 
bounded by $Const\cdot |x|$ for all sufficiently small $|x|$, and 
$x\vee y$ is the maximum of $x$ and $y$.

>From  $\Psi_{net}^+ = \Psi_{com}^+ - \Psi_{int}^+ - \Psi_{ext}^+$ we have,
with $\rho = \rps$,
\begin{eqnarray}
M[\psinetp](\delta,\rho) & = & \delta\frac{1}{2}\left[
\log\rho  +\log\log\zpd+\log(2e)
+{\cal O}\left(\delta \vee \rho \log\rho \vee \frac{\log\log\zpd}{\log\zpd}\right)\right], 
\nonumber \\
 & = & \delta\frac{1}{2}\left[ \log\left(\frac{2e}{\c}\right)
+{\cal O} \left(\frac{\log\log\zpd}{\log\zpd}\right)\right], \quad \delta  \goto 0 . \label{eq:psinetpasym} 
\end{eqnarray}
The ${\cal O}()$ term tends to zero with $\delta$.
Now   $\cs>2e$ so $\log(2e/\cs) < 0$; for some $\delta_2(\cs) > 0 $ 
the bracketed term stays negative
on $0 < \delta<\delta_2(\cs)$. (\ref{Mnegative}) follows
with $\delta_S^+(\cs) = \min(\delta_1(\cs),\delta_2(\cs))$.
\qed

\subsection{Theorem \ref{thm:weakp}}

With (\ref{eq:psinetpasym}) in hand,
it is now convenient to prove 
the lower bound in Theorem \ref{thm:weakp}.

Fix $\cw > 2$ and define $\rwp:=\rwcp:=[\cw\log(1/\delta)]^{-1}$. 
We will show that there is $\delta_W^+=\delta_W^+(\cw) > 0$
so that 
\begin{equation} \label{eq:GoalThm1}
 M[\psinetp-\psiweakp] (\delta,\rwp)  < 0 , \quad 0 < \delta < \delta_W^+.
\end{equation}
Below, Corollary \ref{cor:psiweakp} shows that the maximal function 
becomes `transparent' \-- namely that,
fixing $\cw>2$, there is $\delta_1 = \delta_1(\cw) > 0$ so that
\begin{equation} \label{eq:EvalMaximalWeakSimp}
 M[\psinetp-\psiweakp] (\delta,\rwp) = (\psinetp-\psiweakp)(\delta,\rwp),
          \quad\quad\mbox{for}\quad\delta<\delta_1.
\end{equation}
  
Recall (\ref{eq:33a})-(\ref{eq:first_asym}) and (\ref{eq:defacep}). 
From (\ref{ShannonAsymp})
as the asymptotics for (\ref{eq:defacep}) 
we obtain the following display, in which  $\rho = \rwp$:
\begin{eqnarray}
M[\psinetp-\psiweakp](\delta,\rho) & = & 
\delta\frac{1}{2}\left[
\log\rho-2\rho\log(1/\delta)+\log\log\zpd+\log(2e)  \right.  \nonumber \\
&& \left. \qquad \qquad +{\cal O} \left(\delta \vee \rho\log\rho \vee \frac{\log\log\zpd}{\log\zpd}\right)\right], 
          \nonumber\\
& = & \delta\frac{1}{2} \left[ \log\left(\frac{2e}{\cw}\right)
-\frac{2}{\cw}+{\cal O} \left(\frac{\log\log 1/\delta}{\log1/\delta}\right)
\right] , \quad \delta \goto 0.  \label{eq:Masympweakp}
\end{eqnarray}
Since $\log(1+x) < x$ for $x \in (-1,\infty)$, by setting $1 + x = 2/\cw$
we see that  $\cw > 2$ implies $\log(2e/\cw) - 2/\cw < 0$.
Hence there is $\delta_2(\cw) > 0$ so that 
the term in brackets  is  negative for $\delta$ sufficiently 
small. Define now $\delta_W^+(\cw) = \min(\delta_1(\cw), \delta_2(\cw))$,
 establishing (\ref{eq:GoalThm1}).
 
Looking further ahead to proving the upper bound half of the theorem,  
we record the following remark.
Fix $\cw = 2$. Then as $\delta \goto 0$,
\begin{equation} \label{eq:Goal2Thm1}
 (\psinetp-\psiweakp) (\delta,\rwp)   \goto 0 .
\end{equation}
The implications will emerge in Section 6. 
 \qed

\subsection{Theorem \ref{thm:strongpm}}

The proof of this lower bound  is structurally analogous to the proof
of the lower bound in Theorem \ref{thm:strongp}.

Fix $\cs > 2e$, and define $\rpms:=\rpmsc:= |\cs\log(\delta \sqrt{\pi})|^{-1}$.
As in the proof of Theorem \ref{thm:strongp},
we will show there is $\delta_S^\pm = \delta_S^\pm(\cs) > 0$ so that
\begin{equation}\label{Mnegativepm}
  M[\Psi_{net}^\pm] ( \delta , \rpms ) < 0, \quad 0 < \delta < 
\delta_S^\pm.
\end{equation}
This establishes half of Theorem \ref{thm:strongpm}.
Again, the other half  can be inferred from the fact that for
$\cs<2e$ there is $\delta_0 = \delta_0(\cs) > 0$ with
\[
  \Psi_{net}^\pm]( \delta , \rpms ) > 0, \quad 0 < \delta < \delta_0.
\]
Section 6 will give the details.

 Corollary \ref{cor:psinetpm} below shows that
 the maximal operator $M[]$ becomes transparent
in the limit $\delta \goto 0$;
for some $\delta_1 = \delta_1(\cs) > 0$,
\begin{equation} \label{eq:EvalMaximalStrongCross}
    \Psi_{net}^\pm(\delta,\rpms) = M[\Psi_{net}^\pm] (\delta,\rpms),
    \quad 0 < \delta < \delta_1 .
\end{equation}

The following limiting behavior of the individual exponents as 
$\delta\rightarrow 0$ and/or $\rho\rightarrow 0$ 
will be derived in Section \ref{sec:exponents};
see (\ref{eq:psicomp}), (\ref{eq:psiint_asym}), and (\ref{eq:psiextpm_asym}):
\begin{align}
\psicompm(\delta,\rho) = & \delta\left[\log(1/\delta)+1 + \log_e(2) 
+{\cal O}(\delta \vee \rho\log\rho)\right], \quad\max(\rho,\delta)\rightarrow 0, \label{eq:38a} \\
\psiintpm(\delta,\rho) = & -\frac{1}{2}\delta \left[\log\rho+\log(e/2\pi)  + {\cal O}(\rho\log\rho)\right], \quad\rho\rightarrow 0, \label{eq:38b}\\
\psiextpm(\delta) = & \delta\left[\log \zpmd -\frac{1}{2}\log\log \zpmd + 1 
+{\cal O}\left(\frac{\log\log \zpmd}{\log \zpmd}\right)\right], \quad\delta\rightarrow 0;
\label{eq:38c}
\end{align}
where $\zpmd:=(\delta \sqrt{\pi})^{-1}$.

Combining asymptotics using $\Psi_{net}^\pm = \Psi_{com}^\pm - \Psi_{int}^\pm - \Psi_{ext}^\pm$ yields,
with $\rho = \rpms$,
\begin{eqnarray}
M[\psinetpm](\delta,\rho) & = & \delta\frac{1}{2}\left[
\log\rho  +\log\log\zpmd+\log(2e)
+{\cal O}\left(\delta \vee \rho \log\rho \vee \frac{\log\log\zpmd}{\log\zpmd}\right)\right], 
\nonumber \\
 & = & \delta\frac{1}{2}\left[ \log\left(\frac{2e}{\cs}\right)
+{\cal O}\left(\frac{\log\log \zpmd}{\log \zpmd}\right)   
 \right ], \quad \delta  \goto 0 . \label{eq:psinetpmasym} 
\end{eqnarray}
As $\log(2e/\cs) < 0$, there is $\delta_2(\cs) > 0 $ 
so the term in brackets is negative
for $0 < \delta<\delta_2(\cs)$. Setting $\delta_S^\pm(\cs) = \min(\delta_1(\cs),\delta_2(\cs))$, (\ref{Mnegativepm}) follows. \qed

\subsection{Theorem \ref{thm:weakpm}}

Structurally, the argument for this lower bound resembles 
that in the proof of Theorem \ref{thm:weakp},
in the same way as the proof of  the lower bound in Theorem
\ref{thm:strongpm} resembles that in  Theorem \ref{thm:strongp}.

Fix $\cw > 2$ and define $\rwpm:=\rwcpm:=[\cw\log(1/\delta)]^{-1}$. 
Note that $\rwp=\rwpm$, unlike the strong threshold comparison functions
$\rps$ and $\rpms$, which are not equal. 
We will show that for $\delta_W^\pm= \delta_W^\pm(\cw) > 0$,
\begin{equation}  \label{eq:GoalThm3}
 M[\psinetpm-\psiweakpm] (\delta,\rwpm)  < 0 , 
\quad 0 < \delta < \delta_W^\pm.
\end{equation}
Corollary \ref{cor:psiweakpm} shows that the maximal function
machinery again simplifies for small $\delta$. Thus 
for $\tau > 2$ and for $\delta_1=\delta_1(\cw)> 0$,
\begin{equation} \label{eq:EvalMaximalWeakSimpm}
 M[\psinetpm-\psiweakpm] (\delta,\rwpm) = (\psinetpm-\psiweakpm)(\delta,\rwpm),
         \quad 0 < \delta<\delta_1.
\end{equation}
  
Recall (\ref{eq:psinetpmasym}) and (\ref{eq:defacepm})
with asymptotic behavior following from (\ref{ShannonAsymp}).
We have the following display, in which  $\rho = \rw$, 
\begin{eqnarray}
M[\psinetpm-\psiweakpm](\delta,\rho) & = & 
\delta\frac{1}{2}\left[
\log\rho-2\rho\log(1/\delta)+\log\log\zpd+\log(2e)  \right. \nonumber  \\
&& \left. \qquad \qquad +{\cal O} \left(\delta \vee \rho\log\rho \vee \frac{\log\log\zpmd}{\log\zpd}\right)\right], \nonumber  \\
& = & \delta\frac{1}{2} \left[ \log\left(\frac{2e}{\cw}\right)
-\frac{2}{\cw}+o(1) \right] , \quad \delta \goto 0.  \label{eq:Masympweakpm}
\end{eqnarray}
As in the proof of the lower bound for Theorem \ref{thm:weakp},
for each $\cw>2$ there is $\delta_2(\cw) > 0$ so that 
the term in brackets is negative for all $\delta \in (0, \delta_2)$. 
Setting $\delta_W^\pm(\cw) = \min(\delta_1(\cw),\delta_2(\cw))$,
(\ref{eq:GoalThm3}) follows. \qed

\section{Analysis of the Exponents}\label{sec:exponents}
\setcounter{equation}{0}
\setcounter{table}{0}
\setcounter{figure}{0}

We now verify earlier claims about 
the asymptotic behavior of the exponents.

\subsection{Combinatorial exponents, $\Psi_{com}^\star$}\label{subsec:com}
\setcounter{equation}{0}
\setcounter{table}{0}
\setcounter{figure}{0}

The combinatorial exponents for the simplex and cross-polytope were
defined in (\ref{eq:psicom_related}) 
and (\ref{eq:psicomp_related}) respectively; they obey
\begin{equation}\label{eq:psicom-restated}
\psicomp(\nu,\gamma)=\psicompm(\nu,\gamma)-\nu\log_e(2):=H(\nu)+\nu H(\gamma),
\end{equation}
where again $H(p)=p\log(1/p)+(1-p)\log(1/(1-p))$ is the (base-$e$) Shannon entropy.
Both identities in (\ref{eq:psicom-restated}) derive from the limit
\beqn \label{shannasymp}
    n^{-1} \log{n  \choose {\lfloor p n \rfloor}}  \goto H(p), \quad n \goto  \infty, \quad p\in [0,1],
\eeqn
which of course is fundamental in asymptotic analysis and has proven useful
in earlier research concerning polytopes  \cite{VerSpor,donoho1,Do05_polytope}.
The asymptotic behavior
\begin{equation}\label{eq:psicomp}
\psicomp(\nu,\gamma)=\nu\left[\log(1/\nu)+1+{\cal O}(\nu \vee \gamma\log\gamma)\right],\quad \nu \vee \gamma\rightarrow 0,
\end{equation}
follows directly from that of the Shannon entropy,
\begin{equation} \label{ShannonAsymp}
H(p)=p\log(1/p)+p+{\cal O}(p^2),\quad\quad p\rightarrow 0.
\end{equation}

\newcommand{\psiints}{\Psi_{int}^\star}

\subsection{Internal exponents, $\Psi_{int}^\star$}\label{subsec:int}

The internal-angle exponent is the same for both $\star = +$ and $\star = \pm$;
it was defined in (\ref{eq:psiint}) by 
\begin{equation}\label{eq:psiint-restated}
\psiints(\nu,\gamma) :=\nu(1-\gamma)\left[
\log(\ygam/\gamma)+\frac{1}{2}\log(2\pi)+\frac{\gamma-1}{2\gamma}\ygam^2,
\right]
\end{equation}
where $\ygam$ 
was defined implicitly by
\begin{equation}
\frac{1-\gamma}{\gamma}\ygam=s_{\gamma}, \quad\mbox{and $s_\gamma$ solves}\quad
R(s_\gamma)  = 1 -\gamma;
\end{equation}
here  $R(s)$ -- defined at (\ref{eq:defMillsR}) -- 
is closely associated to a famous quantity in probability theory,
the Mills' ratio of the standard Normal distribution
\cite[Sec 5.38]{KendallVol1}.   The asymptotic properties of $s_\gamma$
as $\gamma \goto 0$ (and hence also of $y_\gamma$)
 were studied in \cite{Do05_polytope}
using properties of Laplace's asymptotic series for $R$.  In the Appendix, we 
refine that  approach, obtaining the following error bounds.

\begin{lemma}\label{lem:ygam}
\begin{equation}\label{eq:lemygam}
\ygam=\frac{\gamma^{1/2}}{1-\gamma}+r_2(\gamma),\quad\quad\mbox{with}\quad 
|r_2(\gamma)|\le 4\gamma^{3/2} \quad\mbox{ for} \quad \gamma\le 1/30.
\end{equation}
\end{lemma}

The behavior  (\ref{eq:33b})
of the internal exponent as $\gamma\rightarrow 0$ 
follows from this lemma directly.
Indeed, substitute the behavior of $y_{\gamma}$ given by 
Lemma \ref{lem:ygam}, and rearrange terms:

\begin{eqnarray}
\psiints(\nu,\gamma) & = &
\nu(1-\gamma)\left[
\log(\ygam/\gamma)+\frac{1}{2}\log(2\pi)+\frac{\gamma-1}{2\gamma}\ygam^2
\right] \nonumber \\
& = & \nu(1-\gamma)\frac{1}{2}\left[
-\log\gamma+\log(2\pi)-\frac{1}{1-\gamma} \right.\nonumber \\
& - & 2\left.\left(\log(1-\gamma)+\gamma^{-1/2}r_2(\gamma)
+\frac{1-\gamma}{2\gamma}r_2^2(\gamma)
+\log\left(1+\frac{r_2(\gamma)}{\gamma(1-\gamma)}\right)\right)\right]
\nonumber \\
& = & \nu(1-\gamma)\frac{1}{2}\left[
-\log\gamma+\log(2\pi/e)+{\cal O}(\gamma)\right],\quad\gamma\rightarrow 0 
\quad\quad(\mbox{by (\ref{eq:lemygam})}) \nonumber \\
& = & -\frac{1}{2}\nu \left[\log\gamma+\log(e/2\pi)  + {\cal O}(\gamma\log\gamma)\right],\quad\gamma\rightarrow 0; \label{eq:psiint_asym}
\end{eqnarray}
this is (\ref{eq:33b}). \qed

\subsection{External exponents, $\Psi_{ext}^\star$}\label{subsec:ext}

Each external exponent
$\Psi_{ext}^\star$ is defined  implicitly through a relation qualitatively resembling 
\begin{equation}\label{eq:fund}
f(x(z),z)=xe^{x^2}-z=0;
\end{equation}
that is to say, we will soon be interested in quantities resembling the solution
$x(z)$. We briefly sketch an
analysis technique for such quantities.  

Our approach approximates
the asymptotic behavior of $x(z)$ for $z$ large 
by
\[
x_2(z) = \sqrt{\log z -\frac{1}{2}\log\log z};
\]
the approximation error obeys
\begin{equation}
|x(z)-x_2(z)|\le \frac{\log\log z}{2\log z} \quad\quad\mbox{as}\quad z\goto\infty.
\nonumber
\end{equation}

The subscript $2$ signals that  $x_2(z)$ is the second
in a {\it sequence} of approximations. The sequence
starts from a very crude approximation,
 $x_1(z)$, and then improves with each stage.  
The initial approximation, $x_1:=\sqrt{\log z}$,
is obtained by treating the 
factor $x$ in (\ref{eq:fund}) as if it were constant, 
so that instead of solving (\ref{eq:fund}), we simply solve 
\[ 
e^{x_1^2}=z.
\]
This approximation, substituted into equation (\ref{eq:fund}), yields an error 
\begin{equation}\label{eq:err1}
f(x_1,z)=z((\log z)^{1/2}-1).
\end{equation}
The next approximation, $x_2$, comes from attempting 
to cancel the $(\log z)^{1/2}$ factor in the above error.  This is done by solving 
\[ 
e^{x_2^2}=z(\log z)^{-1/2},
\]
which indeed yields $x_2(z)$.
This sequence continues on to increasingly accurate approximations, 
but we stop here because the second term
is sufficiently accurate for our purposes.   

\subsubsection{Simplex case $\psiextp$}\label{subsubsec:extp}
\label{subsec:ExtSimp}

Recall the definition given in (\ref{eq:xnu}):
\begin{equation}\label{eq:xnu-restated}
\psiextp(\nu):=\nu \xnu^2-(1-\nu)\log Q(\xnu),
\end{equation}
where $\xnu$ solves
\begin{equation}\label{eq:xnu_again}
\frac{2xQ(x)}{q(x)}+1-\nu^{-1}=0 ;
\end{equation}
here $q(x):=\pi^{-1/2}e^{-x^2}$ and $Q(x)=\int_{-\infty}^x q(y)dy$ 
is related to the Error function by $Q(x)=2(1+$erf$(x))$.
Since there is no closed form solution to $Q(x) = c$ as a function of $c$,
to analyze the implicitly defined $\xnu$, we develop an asymptotic 
approximation using the technique just sketched. Define
\begin{equation} \label{eq:zp}
\zp = \zp (\nu)  :=  (\nu 2\sqrt{\pi})^{-1} , \qquad \tilde{x}_\nu := \left[ \log\zp-\frac{1}{2}\log\log\zp\right]^{1/2}.
\end{equation}
In the Appendix, we prove the approximation result:
\begin{lemma}\label{lem:xnu}
There is $\nu_0 > 0$ so that 
\begin{eqnarray}
\xnu & = & \tilde{x}_\nu +r_3(\nu), \quad 
\tilde{x}_\nu := \left[ \log\zp-\frac{1}{2}\log\log\zp\right]^{1/2}, 
\nonumber \\
|r_3(\nu)| & \le & \frac{1}{2} \left( \log\zp-\frac{1}{2}\log\log\zp\right)^{-1/2}\frac{\log\log\zp}{\log\zp}, \quad 0 < \nu < \nu_0 .\nonumber 
\end{eqnarray}
\end{lemma}

We now plug this approximation  into (\ref{eq:xnu-restated}), 
and derive the asymptotic behavior. As the cumulative distribution of normal $Q()$ 
famously has no known
closed form expression, we
approximate $Q(x)$ for large $x$ using 
the asymptotic series \cite[Sec. 5.38]{KendallVol1},
\[
Q(x)=1-\frac{e^{-x^2}}{2\sqrt{\pi} x}\sum_{r=0}^{\infty}\frac{(r-1/2)!}{(-x^2)^r} .
\]
Keeping the first two terms, and applying bounds from
\cite[eq (5.109)]{KendallVol1}, we have
\[
Q(x)=1-\frac{1}{2\sqrt{\pi}x}e^{-x^2}+{\cal O}(x^{-3}e^{-x^2}), \quad x \goto \infty.
\]
Recalling (\ref{eq:xnu_again}), we now substitute the approximation to $\xnu$ 
from Lemma \ref{lem:xnu}; note that $\xnu^2=\txnu^2+r_4(\nu)$ with $|r_4(\nu)|\le 2(\log\log\zp)/\log\zp$,
and $\zp$ as in (\ref{eq:zp}). Hence,
\begin{eqnarray}
Q(\xnu) & = & Q(\txnu+r_3(\nu)) \nonumber \\
&=&  1-\nu\left[1-\frac{\log\log\zp}{2\log\zp}\right]^{-1/2}\cdot
\left[1+r_3(\nu)/\txnu\right]^{-1}
e^{-r_4(\nu)}+{\cal O}(\nu/\log\zp),\quad\nu\rightarrow 0 ,\nonumber \\
& = & 1-\nu+{\cal O}\left(\nu\frac{\log\log\zp}{\log\zp}\right),
\quad\nu\rightarrow 0, \nonumber
\end{eqnarray}
from which follows 
\begin{equation}
\log Q(\xnu)
=\nu\left[1+{\cal O}\left(\frac{\log\log\zp}{\log\zp}\right)\right], \quad\nu\rightarrow 0.
\end{equation}

We obtain, finally,
\begin{equation}\label{eq:psiextp}
\psiextp(\nu)=
\nu\left[\log \zp -\frac{1}{2}\log\log \zp + 1 
+{\cal O}\left(\frac{\log\log \zp}{\log \zp}\right)\right], 
\quad\nu\rightarrow 0.
\end{equation}
This is (\ref{eq:first_asym}). \qed

\subsubsection{Cross-polytope case: $\psiextpm$}\label{subsubsec:extpm}

The definition given in (\ref{eq:psiextpm}) was
\begin{equation}\label{eq:psiextpm-restated}
\psiextpm(\nu):=\nu \ynu^2-(1-\nu)\log G(\ynu),
\end{equation}
with $\ynu$ the solution of
\begin{equation}
\frac{2yG(y)}{g(y)}+1-\nu^{-1}=0,
\end{equation}
where we recall from before $g(y)=2\pi^{-1/2}e^{-y^2}$ on $y \geq 0$,  and $G(y)= \mbox{erf}(y)=\int_{0}^y g(w)dw$ is the Error function.  
The procedure 
just used in Section \ref{subsec:ExtSimp} also works here.  
We merely state results, omitting proofs. 

Let $\zpm=\zpm(\nu):=(\nu\sqrt{\pi})^{-1}$, and set $\tilde{y}_\nu =\left[ \log\zpm-\frac{1}{2}\log\log\zpm\right]^{1/2}$.
\begin{lemma}\label{lem:ynu}
There is $\nu_0 > 0$ so that 
\begin{equation}
\ynu= \tilde{y}_\nu +r_5(\nu), \quad\quad
|r_5(\nu)|\le \frac{1}{2}\txnu^{-1}\frac{\log\log\zpm}{\log\zpm}, 
\quad 0 < \nu< \nu_0 \nonumber .
\end{equation}
\end{lemma}

This approximation is motivated by the asymptotic series of 
$2yG(y)/g(y)$, giving
\[
\ynu e^{\ynu^2}-\pi^{-1/2}\nu^{-1}={\cal O}(\ynu^{-2}).
\]
The series is identical to the series motivating $\xnu$ in $\psiextp$ 
but now $\zpm(\nu)=\zpm:=\pi^{-1/2}\nu^{-1}$.  
The precise bound on the remainder, $r_5(\nu)$, 
can be recovered by following the same steps as 
in the proof of Lemma \ref{lem:xnu}, replacing $J(x,\nu)$ 
in that proof by $2ye^{-y^2}\int_0^y e^{-w^2}dw+1-\nu$. 

\qed

The asymptotic behavior of the external exponent,
\begin{equation}\label{eq:psiextpm_asym}
\psiextpm(\nu)=
\nu\left[\log \zpm -\frac{1}{2}\log\log \zpm + 1 
+{\cal O}\left(\frac{\log\log \zpm}{\log \zpm}\right)\right], \quad\nu\rightarrow 0
\end{equation}
follows by substituting  $\tilde{y}_\nu$, as justified by Lemma \ref{lem:ynu}.

\subsection{Maximal Function for  $\Psi_{net}^\star$ }
\label{subsec:net_bound}

We now support our earlier claim (\ref{eq:EvalMaximalStrongSimp})  that
$M[\Psi_{net}^\star] = \Psi_{net}^\star$
along the trajectory $(\delta,\rps)$  for $\delta $ small enough. 

\begin{corollary}\label{cor:psinetp}
Fix $\cs>2e$, and recall the definition  
$\rps:=\rpsc:=|\cs\log(\delta 2\sqrt{\pi})|^{-1}$.  
There is $\delta_1(\cs) > 0$ so that
\[
\psinetp(\delta,\rps)= M[\psinetp](\delta,\rps) ,
\quad 0 < \delta < \delta_1(\cs).
\]
\end{corollary}

This follows from two lemmas, proved in the Appendix,
which clarify how $\psinet$ changes with $\nu$ in the regime of interest.

\begin{lemma}\label{lem:psinetp_nu}
Fix $\cs>2e$.
There is $\delta_1 = \delta_1(\cs) > 0$ so that for $0 < \delta < \delta_1$, 
and $0<\gamma\le \rps$,
$\psinetp(\nu,\gamma)$ is a decreasing function of 
$\nu$ for $\nu\in[\delta,1)$.
\end{lemma}

\begin{lemma}\label{lem:psinetpgamma}
For $0 < \gamma< \gamma_0$, $\psinetp(\nu,\gamma)$ 
is an increasing function of $\gamma$.
\end{lemma}

Similar results hold for the cross-polytope 
[note the slight difference in definition
between $\rps$ and $\rpms$].

\begin{corollary}\label{cor:psinetpm}
Pick $\cs>2e$ and again set
$\rpms:=\rpmsc:=|\cs\log(\delta \sqrt{\pi})|^{-1}$.   
For $\delta < \delta_0(\cs)$, 
$\psinetpm(\nu,\gamma)$ obtains its maximum value over $\nu\in[\delta,1)$ and $\gamma\le\rpms$ at 
$(\nu,\gamma)=(\delta,\rpms)$: 
\[
\psinetpm(\delta,\rpms)= M[\psinetpm](\delta,\rpms).
\]
\end{corollary}

We omit the proof, whose arguments parallel those
for Lemmas \ref{lem:psinetp_nu} and \ref{lem:psinetpgamma}.

\subsection{Maximal Function for $\Psi_{net}^\star - \Psi_{face}^\star$ }
\label{subsec:net_weak_bound}

We now consider the maximal function associated with the weak exponent, 
establishing the earlier claim (\ref{eq:EvalMaximalWeakSimp}).

\begin{corollary}\label{cor:psiweakp}
Fix $\cw>2$.  There is $\delta_1=\delta_1(\cw) > 0$ so that 
\[
(\psinetp-\psiweakp)(\delta,\rwp)=M[\psinetp-\psiweakp] (\delta,\rwp) ,
\quad 0 < \delta < \delta_1.
\]
\end{corollary}

This  follows immediately from 
the next lemmas, which are proven in the Appendix.

\begin{lemma}\label{lem:psiweakp_nu}
Fix $\cw>2$.  For $0 < \delta < \delta_1(\cw)$,
$0 < \gamma \leq \rwp$,

$(\psinetp-\psiweakp)(\nu,\gamma)$ is a decreasing function of $\nu$ over $\nu\in [\delta,1)$.  
\end{lemma}

\begin{lemma}\label{lem:psiweakp_gamma}
Fix $\cw>2$.  For $\delta \in(0,\delta_1(\cw))$, 
$\rho \in (0,\rwp)$ and
$\nu\in [\delta,1)$,
$(\psinetp-\psiweakp)(\nu,\gamma)$ is an increasing function of $\gamma$ ,
$0 \leq \gamma \leq r_W^+(\delta)$.
\end{lemma}

Similar results for the cross-polytope are obtained by following the same 
arguments line-by-line with appropriate substitutions.
One obtains the following, though we omit the argument.

\begin{corollary}\label{cor:psiweakpm}
Fix $\cw>2$.  There is $\delta_1 = \delta_1(\cw) > 0$ so that
\[
(\psinetpm-\psiweakpm)(\delta,\rwpm)= M[\psinetpm-\psiweakpm](\delta,\rwpm),
 \quad 0 <  \delta < \delta_1.
\]
\end{corollary}

\section{Beyond Proportional Growth}
\setcounter{equation}{0}
\setcounter{table}{0}
\setcounter{figure}{0}

Theorem \ref{thm:5} can be reformulated as follows.
\begin{theorem} \label{thm:NonPropGrowthRestated}
Let $N_n$ grow subexponentially with $n$.
\bitem
\item {\it Strong Exponents.}
Fix $\cs > 2e$ and consider a sequence
$(k_n)$ with  $k_n \leq  n \cdot \rstarsnc$.
There is a sequence $(\eps_n)$ with $N_n \eps_n \goto \infty$ and
\begin{equation} \label{eq:StrongGoal}
    N_n^{-1} \log( D_{\ell,n}^\star) \leq  - \eps_n , \qquad \ell = n+1, n+3,\ldots .
\end{equation}
\item {\it Weak Exponents.}
Fix $\cw > 2$ and consider 
a sequence $(k_n)$ with $k_n \leq  n \cdot \rwnc$.
There is a sequence $(\eps_n)$ with $N_n \eps_n \goto \infty$ and
\begin{equation} \label{eq:WeakGoal}
    N_n^{-1} (\log( D_{\ell,n}^\star) - \log f_k(Q) ) \leq  - \eps_n , \qquad \ell = n+1, n+3,\ldots.
\end{equation}
\eitem
\end{theorem}

To venture outside the proportional growth setting
requires to strengthen all previous arguments. First,
we have to show not just that each maximal function is negative before its
first zero, but that it is sufficiently negative in a quantitative sense.
Fortunately, the hard work has already been done;
summarizing the implications
of (\ref{eq:psinetpasym}), 
(\ref{eq:Masympweakp}), 
(\ref{eq:psinetpmasym}), and 
(\ref{eq:Masympweakpm}), we have:
\begin{lemma} \label{lem:SuffNegative}
Let $\delta_n = n/N_n$.
\bitem
\item {\it Strong Exponents.}
Fix $\cs > 2e$. There are $\zeta_S^\star(\cs) > 0$ so that for $n > n_0$
\begin{equation} \label{Mstrongp}
     M [ \Psi_{net}^\star ] ( \delta_n , \rstarsdn )  < - \zeta_S^\star(\cs) \delta_n.
\end{equation}
\item {\it Weak Exponents.}
Fix $\cw > 2$. There are $\zeta_W^\star(\cw) > 0$ so that for $n > n_0$
\[
     M [ \Psi_{net}^\star - \Psi_{face}^\star ] ( \delta_n , \rwdn)  < -  \zeta_W^\star(\cw) \delta_n.
\]
\eitem
\end{lemma}

We must also strengthen the previously-discussed 
inequalities (\ref{BoundComb}),(\ref{BoundInt}),
(\ref{BoundExt}),(\ref{BoundCombPM}),(\ref{BoundIntPM}),
and (\ref{BoundExtPM}), giving precise information about the remainders.
We start with the combinatorial exponent.
\begin{lemma} \label{lem:ComBound}
\begin{equation}
N^{-1}\log C_{\ell,n}^\star  
\leq  \Psi_{com}^\star(\nu_{\ell,n},\gamma_{\ell,n}) 
+O(N^{-1} \log(N)),
\label{ComBound1} 
\end{equation}
where the $O()$ term is uniform in $\ell = n+1, n+3, \dots$.
\end{lemma}

The proof is given in Section \ref{sec:ReanalCombin}.
We next consider the external angles.
\begin{lemma}
\begin{equation} \label{ExtBound1}
N^{-1}\log\alpha(T^{\ell},T^{N-1})\le-\Psi_{ext}^+(\nu_{\ell,n})+
O(N^{-1} \log(N)),
\end{equation}
where the $O()$ is uniform in $\ell = n+1, n+3, \dots$.
Similarly,
\begin{equation} \label{GoalExtExplogNoN}
N^{-1}\log\alpha(F^{\ell},C^{N})\le -\Psi_{ext}^{\pm}(\nu_{\ell,n})+
O(N^{-1}\log N ),
\end{equation}
where the $O()$ is uniform in $\ell = n+1, n+3, \dots$.
\end{lemma}
For the proof see Section \ref{sec:ReanalExtern}.
We finally consider the internal angles.
\begin{lemma} \label{lem:IntBound}
Fix  $\cs > 2e$. 
\begin{equation} \label{IntBound1}
N^{-1}\log\beta(T^k,T^{\ell})\le -\psiint^\star(\nu_{\ell,n},\gamma_{\ell,n})  + o(1)\Psi_{net}^\star(\nu_{\ell,n},\gamma_{\ell,n}) + O(N^{-1} \log N ),
\end{equation}
where the $o()$ is uniform in $\ell = n+1, n+3, \dots$
and in $k = 1,\dots, \lfloor n \cdot \rstarsdn \rfloor$.
Fix $\cw > 2$. 
\begin{eqnarray} \label{IntBound2}
N^{-1}\log\beta(T^k,T^{\ell}) &\le& -\psiint^\star(\nu_{\ell,n},\gamma_{\ell,n})  \nonumber \\
   && + o(1)(\Psi_{net}^\star(\nu_{\ell,n},\gamma_{\ell,n}) - \Psi_{face}^\star(\nu_{\ell,n},\gamma_{\ell,n})) + O(N^{-1} \log N ), 
\end{eqnarray}
where the $o()$ is uniform in $\ell = n+1, n+3, \dots$
and in $k = 1,\dots, \lfloor n \cdot \rwdn \rfloor$.
\end{lemma}
We also need analogous results for the number of faces of $T^{N-1}$ and $C^N$.
\begin{lemma}\label{lem:face_order}
\[
-N^{-1}\log f_k (T^{N-1})\le \Psi_{face}^+(\nu,\gamma) +O(N^{-1}\log(N)), 
\]
where the $O()$ is uniform in $k=1,2,\ldots,n$.  Similarly, 
\[
-N^{-1} \log f_k (C^{N})\le \Psi_{face}^{\pm}(\nu,\gamma) +O(N^{-1}\log(N)), 
\]
where the $O()$ is uniform in $k=1,2,\ldots,n$.
\end{lemma}
For the proof see Section \ref{subsec:precise_face}.

These Lemmas easily combine to finish the
argument for Theorem \ref{thm:NonPropGrowthRestated}.
Under the subexponential growth assumption $\log(N_n) = o(n)$,
the remainder terms 
\[
O(N_n^{-1} \log(N_n)) = o(n/N_n) = o( \delta_n).
\]
Hence
the remainders are much smaller than the bounds on $M[]$
terms associated with (\ref{eq:WeakGoal}), (\ref{Mstrongp}). 
  Consider the case of the strong exponent for
the cross-polytope. Uniformly in $\ell = n+1, n+3, \dots$,
\begin{eqnarray*}
  N^{-1} \log(D_{\ell,n}^\pm) &\leq& 
       M[\Psi_{net}^\pm](\delta_n, \rpmsdn) \cdot (1+o(1))  
       + O(N_n^{-1} \log(N_n)) \\
      &\leq &  -\zeta_S^\pm(\cs) \cdot \delta_n \cdot (1 + o(1))
      +  O(N_n^{-1} \log(N_n)) \\
      &<& -(\zeta_S^\pm(\cs)/2) \cdot  \delta_n, \quad n > n_0 .
\end{eqnarray*}
Hence (\ref{eq:StrongGoal}) follows, with $\star = \pm$
and  $\eps_n = (\zeta_S^{\pm}(\cs)/2) \delta_n$.
The rest of Theorem  \ref{thm:NonPropGrowthRestated} follows similarly. \qed

It remains to prove Lemmas \ref{lem:ComBound} - \ref{lem:face_order}.
This we do in the coming subsections.
\subsection{Combinatorial exponents}
\label{sec:ReanalCombin}

 Stirling's formula provides error bounds for the 
combinatorial exponents.

\begin{lemma}[Stirling's inequality, \cite{Cheney}]\label{lem:ShannonBounds}
\[
(2\pi n)^{1/2}\left(\frac{n}{e}\right)^n\le n!
\le \frac{5}{4}(2\pi n)^{1/2}\left(\frac{n}{e}\right)^n
\quad\mbox{for}\quad n\ge 1.
\]
\end{lemma}

To verify (\ref{ComBound1}), recall the combinatorial factors  
\[
C_{\ell,n}^+=2\left(\begin{array}{c}N \\ \ell \end{array}\right)
\left(\begin{array}{c}\ell \\ k+1 \end{array}\right) \quad\mbox{and}\quad
C_{\ell,n}^{\pm}=2^{\ell+1}\left(\begin{array}{c}N \\ \ell+1 \end{array}\right)
\left(\begin{array}{c}\ell+1 \\ k+1 \end{array}\right)
\]
and that $\nu_{\ell,n}=\ell/N$ and $\gamma_{\ell,n}=k/\ell$.

Using Lemma \ref{lem:ShannonBounds} we arrive at  
\begin{equation}\label{eq:psicom_precise}
C_{\ell,n}^{\star}\le \frac{5}{8\pi}\cdot 
N e^{N\Psi_{com}^{\star}(\nu,\gamma)},
\end{equation}
establishing (\ref{ComBound1}).



\subsection{External Angle}
\label{sec:ReanalExtern}

\subsubsection{Simplex case $\Psi_{ext}^+$}\label{sec:extp}

It is enough to show that uniformly over $\ell\ = n+1, n+3, \dots$,
\begin{equation} \label{mainExtSimp}
N^{-1}\log\alpha(T^{\ell},T^{N-1})\le-\Psi_{ext}^+( \nu_{\ell,n})+
N^{-1}\log\left(N\right);
\end{equation}
of course, the remainder term is $O( \log(N)/N)$. 
The simplex part of Lemma \ref{lem:IntBound} follows.

The external angle for the simplex is given by
\begin{equation}\label{eq:extp_def}
\alpha(T^{\ell},T^{N-1})=\sqrt{\frac{\ell+1}{\pi}}\int_0^{\infty}
e^{-(\ell+1)x^2}\left(\frac{1}{\sqrt{\pi}}\int_{-\infty}^x e^{-y^2}dy
\right)^{N-\ell-1}dx.
\end{equation}
As before, $Q(x):=\pi^{-1/2}\int_{-\infty}^x e^{-y^2}dy$.
Recall that $\nu_{\ell,n} = \ell/N$ and rewrite the simplex external angle as
\begin{equation}\label{eq:extpgoal}
\alpha(T^{\ell},T^{N-1})=\sqrt{\frac{\ell+1}{\pi}}\int_0^{\infty}
\exp(-N[\nu_{\ell,n} x^2+(\nu_{\ell,n}-1)\log Q(x)])\frac{e^{-x^2}}{Q(x)}dx.
\end{equation}

The factor $N$ in the integral might
suggest the use of Laplace's method as
in \cite{Do05_polytope}.
A simpler, direct approach is possible.  The following is obvious but very useful.

\begin{lemma} \label{lem:BeatLaplace}
Let $\psi: [0,\infty) \mapsto \bR$ achieve its global minimum at $x^*$ and
let $\varphi:[0,\infty) \mapsto [0,\infty)$ be integrable.
Then
\begin{equation} \label{EllOneEllInfty}
   \int_0^\infty \exp(-N \psi(x)) \varphi(x) dx \leq \exp(-N\psi(x^*))  \int_0^\infty \varphi(x) dx.
\end{equation}
\end{lemma}

\noindent
Recall that $\xnu$ is the minimizer of 
$[\nu x^2+(\nu-1)\log Q(x)]$, and 
\[
\Psi_{ext}^+(\nu):=\nu \xnu^2+(\nu-1)\log Q(\xnu).
\]
Apply Lemma \ref{lem:BeatLaplace}
to the integral (\ref{eq:extpgoal}); set $\psi(x) =  [\nu x^2+(\nu-1)\log Q(x)]$
and $\varphi(x) = e^{-x^2}/Q(x)$.
Because $\psi(x^*) = \Psi_{ext}^+(\nu)$
and $\int \varphi = (3\pi)/8$,
(\ref{EllOneEllInfty}) yields 
\begin{equation}\label{eq:psiext_precise_p}
\alpha(T^{\ell},T^{N-1})\le e^{-N\Psi_{ext}^+(\nu)}\cdot \sqrt{\frac{\ell+1}{\pi}}
\cdot \frac{3\pi}{8} \leq \sqrt{N+1} \cdot e^{-N\Psi_{ext}^+(\nu)} .
\end{equation}
(\ref{mainExtSimp}) follows. \qed

\subsubsection{Cross-Polytope case $\Psi_{ext}^\pm$}\label{sec:extpm}

Our goal is to prove (\ref{GoalExtExplogNoN}).
We introduce a perturbed version of $\nu_{\ell,n}$; 
\[
\widehat{\nu}_{\ell,n} = \nu_{\ell,n} + \frac{1}{2N}.
\]
Note that $\widehat{\nu}_{\ell,n}\in [0,1)$ as is 
the unperturbed $\nu_{\ell,n}$.
This is used in our first step, where we find that 
it appears naturally in the bound
\begin{equation} \label{extCrossGoal1}
N^{-1}\log\alpha(F^{\ell},C^N)\le -\Psi_{ext}^{\pm}(\widehat{\nu}_{\ell,n})+
N^{-1}\log\left(N\right), \quad N > 3.
\end{equation}
Note that the remainder is $O(N^{-1} \log(N))$ 
uniformly over $\ell \in [n+1,N-1]$,
as our goal requires. 
Indeed (\ref{GoalExtExplogNoN}) is an inequality like (\ref{extCrossGoal1})
but with  $\nu_{\ell,n}$ rather than $\widehat{\nu}_{\ell,n}$.
In our second step, we verify that 
the perturbation of the argument is unimportant:
\begin{equation} \label{extCrossGoal2}
  \Psi_{ext}^{\pm}(\widehat{\nu}_{\ell,n}) = \Psi_{ext}^{\pm}({\nu}_{\ell,n}) + O( N^{-1} \log(N)),
\end{equation}
uniformly over $\ell = n+1,n+3,\ldots$; $\ell<N$. The cross-polytope half 
of Lemma \ref{lem:IntBound}
then follows. It remains to show (\ref{extCrossGoal1})-(\ref{extCrossGoal2}).

The external angle for the cross-polytope is given by 
\begin{equation}\label{eq:cp_ext_angle_def}
\alpha(F^{\ell},C^N)=\sqrt{\frac{\ell+1}{\pi}}\int_0^{\infty}
e^{-(\ell+1)y^2}\left(\frac{2}{\sqrt{\pi}}\int_{0}^y e^{-w^2}dw
\right)^{N-\ell-1}dy.
\end{equation}
Following the same approach as for the simplex, recall 
$G(y):=2\pi^{-1/2}\int_{0}^y e^{-w^2}dw$ and rewrite the cross-polytope external angle as 
\begin{equation} \label{eq:ExtCrossInt}
\alpha(F^{\ell},C^N)=\sqrt{\frac{\ell+1}{\pi}}\int_0^{\infty}
\exp(-N[\widehat{\nu}_{\ell,n}y^2+(\widehat{\nu}_{\ell,n}-1)\log G(y)])\left(\frac{e^{-y^2}}{G(y)}\right)^{1/2}dy.
\end{equation}
Let $\psi(y;\nu) = \nu y^2+(\nu-1)\log G(y)]$, and  $\widehat{y}_{\nu}$ be
the minimizer of $\psi(\cdot; \nu)$. Set
\[
\Psi_{ext}^\pm(\widehat{\nu}):=
\nu \widehat{y}_\nu^2+(\nu-1)\log G(\widehat{y}_\nu).
\]
Apply Lemma \ref{lem:BeatLaplace} to the integral (\ref{eq:ExtCrossInt});
set $\psi = \psi(\cdot; \widehat{\nu}_{\ell,n})$, and
$\varphi(y) = \left(\frac{e^{-y^2}}{G(y)}\right)^{1/2}$.
The factor $\exp(-N\psi(x^*)) = \exp(-N\psiextpm(\widehat{\nu}_{\ell,n}))$, 
while $\int \varphi \leq 2.175$; we obtain 
\begin{equation}\label{eq:psiext_precise_pm}
\alpha(F^{\ell},C^N)\le \frac{5}{4}\sqrt{{\ell+1}}\exp(-N\Psi_{ext}^{\pm}(\widehat{\nu}_{\ell,n})),
\end{equation}
Hence
\[
N^{-1}\log\alpha(F^{\ell},C^N)\le -\Psi_{ext}^{\pm}(\widehat{\nu}_{\ell,n})+N^{-1}\log\left(\frac{5}{4}\sqrt{{\ell+1}}\right),
\]
from which (\ref{extCrossGoal1}) follows.  

We earlier studied the asymptotic behavior of 
$\psiextpm(\nu)$; see (\ref{eq:psiextpm_asym}). The effect of the 
perturbation $1/2N$ in $\widehat{\nu}_{\ell,n}$ can be bounded simply. 
Put  $z_n^\pm = z^\pm(\widehat{\nu}_{\ell,n})$. Then
\begin{eqnarray}
\psiextpm(\widehat{\nu}_{\ell,n}) & = &
(\nu_{\ell,n} +1/2N) \left[\log z_n^\pm -\frac{1}{2}\log\log z_n^\pm + 1 
+{\cal O}\left(\frac{\log\log z_n^\pm}{\log z_n^\pm}\right)\right], \quad n \goto \infty, \nonumber \\
& = & \nu_{\ell,n} \left[\log z_n^\pm -\frac{1}{2}\log\log z_n^\pm + 1  \right]
+{\cal O}\left(\log z_n^\pm \vee \frac{\log\log z_n^\pm}{\log z_n^\pm}\right)/N, \quad n \goto \infty. \nonumber 
\end{eqnarray}
Our goal (\ref{extCrossGoal2}) follows. 
Combined with (\ref{extCrossGoal1}) we obtain (\ref{GoalExtExplogNoN}). \qed

\subsection{Internal angle}\label{sec:int}

We aim to demonstrate (\ref{IntBound1}).
We again introduce perturbed  variables:
\[
\tilde{\nu} = \tilde{\nu}_{\ell,n} = \frac{\ell+2}{N},  \qquad \tilde{\gamma} =  \tilde{\gamma}_{\ell,n} := \frac{k+1}{\ell+2}.
\]
Our plan is to first show that for $n > n_0$ 
\begin{equation} \label{IntGoal1}
N^{-1}\log\beta(T^k,T^{\ell})\le -\psiint(\tilde{\nu},\tilde{\gamma})
+N^{-1}\log\left[\frac{2}{\pi}(N+3)^{5/2}\right].
\end{equation}
The remainder here is $O(N^{-1} \log(N))$.
We then show that the perturbation of variables has a negligible impact:
\begin{equation} \label{IntGoal2}
  \psiint(\tilde{\nu},\tilde{\gamma}) -\psiint(\nu,{\gamma})  = o( \Psi_{net}^\star )
\end{equation}
uniformly in $0 \leq k \leq n\cdot \rstarsdn$.
Our goal (\ref{IntBound1}) follows. It remains to
prove (\ref{IntGoal1}), (\ref{IntGoal2}).

An expression for the internal angle was developed in \cite{Do05_polytope}: 
\begin{equation} \label{intdefbeta}
\beta(T^k,T^{\ell})
=\left(\pi\frac{\ell+2}{k+1}\right)^{1/2} 2^{k-\ell} g_{T+W_m}(0);
\end{equation}
here $g_{T+W_m}(0)$ denotes the probability density of a certain
random variable expressible as a sum of $m+1$ independent random variables;
here $m = \ell - k+1$.
\cite{Do05_polytope} used large deviations analysis to bound this term
using a certain nonnegative convex rate function $\Lambda^*: [0,\infty) \mapsto [0,\infty)$;
the bound was:
\begin{eqnarray*}
g_{T+W_m}(0) &\le& \frac{2}{\sqrt{\pi}}\frac{m^2}{2\theta}\int_0^{\sqrt{2/\pi}} 
y\exp\left(-m\left[\left(\frac{m}{2\theta}\right)y^2+\Lambda^*(y)\right]\right)dy + \frac{2}{\sqrt{\pi}}\exp\left(-\frac{m^2}{\pi\theta}\right)\\
&=:& I_m+II_{m},
\end{eqnarray*}
say, with $m=\ell-k+1$ and $\theta=k+1$.  
The second term was argued to be negligible in the proportional growth setting by soft
analysis; later below we will check that it is still negligible in the current
non-proportional growth setting.

Focusing on the supposedly dominant term $I_m$, substitute in the values for $m$ and $\theta$, 
and recall that $\tilde{\gamma}=(k+1)/(\ell+2)$:
\[
I_m=\frac{1}{\sqrt{\pi}}\frac{(\ell-k+1)^2}{(k+1)}\int_0^{\sqrt{2/\pi}}
y\exp\left(-(\ell-k+1)\left[\left(\frac{1-\tilde{\gamma}}{\tilde{\gamma}}\right)\frac{y^2}{2}+\Lambda^*(y)\right]\right)dy.
\]
The integral here can be rewritten as:
\[
  J_m := \int_0^{\sqrt{2/\pi}} y \exp(- N\tilde{\nu}(1 - \tilde{\gamma})\xi_{\tilde{\gamma}}(y)) dy .
\]
where, consistent with earlier definitions,
\[
\xi_{\tilde \gamma}(y) =
\left[\left(\frac{1-\tilde{\gamma}}{\tilde{\gamma}}\right)\frac{y^2}{2}+\Lambda^*(y)\right].
\]
Note that
 $y_{\tilde{\gamma}}$  is the minimum of $\xi_{\tilde{\gamma}} (y)$.
Again apply Lemma \ref{lem:BeatLaplace}   to bound $J_m$;
setting $\psi = \tilde{\nu}(1 - \tilde{\gamma})\xi_{\tilde{\gamma}}(y)$ and 
$\varphi = y 1_{[0,\sqrt{2/\pi}]}$, (\ref{EllOneEllInfty}) gives
\[
   J_m \leq \exp(-N  \tilde{\nu}(1 - \tilde{\gamma})\xi_{\tilde{\gamma}}(y_{\tilde{\gamma}}))  / \pi.
\]
Note that
\[
\psiint(\tilde{\nu},\tilde{\gamma}):=
\tilde{\nu}(1 - \tilde{\gamma})[\xi_{\tilde{\gamma}}(y_{\tilde{\gamma}})+\log 2],
\]
and so
\[
    2^{k-\ell-1} \exp(-N \tilde{\nu}(1 - \tilde{\gamma})\xi_{\tilde{\gamma}}(y_{\tilde{\gamma}})) = \exp(-N\Psi_{int}(\tilde{\nu},\tilde{\gamma})).
\]
Noting the presence of a factor $2^{k-\ell}$ in (\ref{intdefbeta}) and
noting that $\ell+1-k \leq N$, we obtain 
\begin{eqnarray}
\beta(T^k,T^{\ell}) & \le & 2\left(\frac{\ell+2}{k+1}\right)^{1/2}
\frac{(\ell-k+1)^2}{(k+1)}
\exp(-N\Psi_{int})/\pi   + \sqrt{\pi (N+2)} \cdot 2^{k-\ell} II_{m} \nonumber \\
& \le & 2\left( N+2 \right)^{5/2}\cdot \exp(-N\Psi_{int}(\tilde{\nu},\tilde{\gamma})) +\sqrt{\pi (N+2)} \cdot 2^{k-\ell} II_{m} .
\end{eqnarray}
This essentially verifies (\ref{IntGoal1}).

However, it remains to verify that $II_{m} \ll I_m$.
 Put $\mu = \sqrt{2/\pi}$ and recall from \cite{Do05_polytope} that
 $\mu = E (T + W_m)$. We focus on
$y = \mu =  \sqrt{2/\pi}$ and use the fact that the 
large deviations rate function always vanishes at the underlying mean,
i.e. $\Lambda^*(\mu) = 0$ essentially
by definition. Then
\[
-m\left[\left(\frac{m}{2\theta}\right)\mu^2+\Lambda^*(\mu)\right] =   \frac{-m^2}{\pi\theta}.
\]
It follows that
\[
   II_{m}  =  \exp( - N  \cdot \tilde{\nu}(1 - \tilde{\gamma}) \xi_{\tilde{\gamma}}(\mu) ) \cdot \sqrt{\frac{2}{\pi}}.
\]
But by definition of $y_{\tilde{\gamma}}$ as the minimizer of $\xi_{\tilde{\gamma}}$,
and the asymptotic $y_{\tilde{\gamma}} \goto 0$,
\[
 \xi_{\tilde{\gamma}}(\mu) >  \xi_{\tilde{\gamma}}(y_{\tilde{\gamma}});
\]
in fact $ \xi_{\tilde{\gamma}}(\mu)  \sim  \tilde{\gamma}^{-1} \mu^2 \gg \log(\tilde{\gamma}^{-1}) \sim \xi_{\tilde{\gamma}}(y_{\tilde{\gamma}})$
as $\tilde{\gamma} \leq \rwdn( 1+ o(1)) \goto 0$.
Hence $II_{m}$ is exponentially smaller than $J_m$, and (\ref{IntGoal1}) is fully proven.

As for  (\ref{IntGoal2}), recall that 
\begin{equation}\label{eq:psiint_asym_k}
\psiint(\tilde{\nu},\tilde{\gamma}) = 
 -\frac{1}{2}\tilde{\nu} \left[\log\gamma+\log(e(1+k^{-1})/2\pi) + o(1) \right], 
\end{equation}
while, if $\rho = \rstars$
\[
\Psi_{net}^\star(\nu,\eta \rho)  = \nu\frac{1}{2}\left[ \log\left(\frac{2e}{\cs}\right) + \log(\eta) 
+ o(1) \right] .
\]
Look now in the vicinity of $k = \gamma n$, where 
$\gamma = \eta \cdot \rstarsdn$.
\[
\frac{ |\psiint^\star(\tilde{\nu},\tilde{\gamma}) - \psiint^\star({\nu},{\gamma}) | }{ |\Psi_{net}^\star(\nu,\gamma) |}  \leq  
\frac{\min(1, 1/ \eta \cdot \frac{1}{n\cdot\rstarsdn} ) + o(1)  }{| \frac{2e}{\cs} + \log(\eta) + o(1)| }  = o(1) . 
\]
Here all the $o(1)$'s are uniform in $0 \leq \eta \leq 1$.

The argument for (\ref{IntBound2}) is similar to that of 
(\ref{IntBound1}) detailed above, replacing $\Psi_{net}^\star$ by 
$\Psi_{net}^\star - \Psi_{face}^\star$. \qed

\subsection{Face Counts of $T^{N-1}$ and $C^N$}\label{subsec:precise_face}
The number of $k-$faces for the simplex and cross-polytope are
\[
f_k (T^{N-1}) = \left(\begin{array}{c}N \\ k+1 \end{array}\right)
\quad\mbox{and}\quad
f_k (C^N) = 2^{k+1} \left(\begin{array}{c}N \\ k+1 \end{array}\right).
\]

Invoking Lemma \ref{lem:ShannonBounds} and recalling that 
$\nu_{\ell,n}=\ell/N$ and $\gamma_{\ell,n}=k/\ell$, we arrive at  
\begin{equation}\label{eq:psiface_precise_p_upper}
f_k (T^{N-1}) \ge \frac{8}{25}\sqrt{\frac{2}{\pi}}
N^{-1} e^{N\Psi_{face}^+(\nu,\gamma)},
\end{equation}
\quad\mbox{and}
\begin{equation}\label{eq:psiface_precise_pm}
f_k (C^N) \ge \frac{16}{25}\sqrt{\frac{2}{\pi}}
N^{-1} e^{N\Psi_{face}^{\pm}(\nu,\gamma)},
\end{equation}
establishing Lemma \ref{lem:face_order}.

\newcommand{\Psinet}{\Psi_{net}}
\newcommand{\Psiint}{\Psi_{int}}
\newcommand{\Psiext}{\Psi_{ext}}
\newcommand{\Psicom}{\Psi_{com}}
\newcommand{\Psiface}{\Psi_{face}}

\section{Upper Bounds on Phase Transitions}
\setcounter{equation}{0}
\setcounter{table}{0}
\setcounter{figure}{0}

Until this point, we have focused on establishing lower bounds on the several $\rho$-functions
introduced in Section 1. Our work so far has given the lower-bound ``half'' of Theorems \ref{thm:weakp}-\ref{thm:strongpm}; 
we now give upper bounds on the $\rho$-functions and complete the proof of 
Theorems \ref{thm:weakp}-\ref{thm:strongpm}.

We remark, parenthetically, that the ``half'' already proven is the more surprising/interesting
part of the result, in view of applications. However, the remaining part settles any question
about whether the lower bounds have slack, i.e. whether they
 actually agree with the precise phase transitions.
 
For establishing {\it lower} bounds on the $\rho$'s, we have been applying
{\it upper} bounds on the combinatorial factor, and on the internal and external angles.
Now that we want {\it upper} bounds on the $\rho$'s, we will turn to
{\it lower} bounds on the combinatorial factor and the angles.

The required lower bounds will be developed in later subsections of
this section, effectively we will be using standard
ideas such as Stirling's inequality,
 Laplace's method and the Saddlepoint method.
 
Before turning to those lower bounds,
 we give the arguments completing
 the proofs of Theorems \ref{thm:weakp}-\ref{thm:strongpm}.
 
\subsection{Upper Bounds on Strong Phase Transition}
 
 The key to tying down the strong phase transitions
 $\rho_S^\star$ is to use the fact that $\Psinet^\star(\delta,\cdot)$ 
 makes a sign change at $\rho_S^\star$.
 Indeed, by definition, $\Psinet^\star(\delta,\rho)$ has a zero
 at $\rho = \rho^\star(\delta)$; but actually it is strictly increasing
 in the vicinity of this zero. For sufficiently small $\eps > 0$,
 we can find $\rho = \rho_\eps > \rho_S^\star(\delta)$ so that
 \[
    \Psinet^\star(\delta,\rho) >  2 \eps.
 \]
 Set now $k_n = \lfloor \rho_\eps n \rfloor$;  for
 all sufficiently large $n$,
\begin{equation} \label{epsthresh}
      \Psinet^\star\left( \frac{n+2}{N}, \frac{k+1}{n+2}\right) > \eps .
\end{equation}
We now invoke lemmas placing lower bounds
on the combinatorial, internal and external angle factors.

\begin{lemma}
There is an absolute constant $c_1 > 0$ so that
\[
  C_{\ell,n}^\star \geq c_1 \cdot \frac{\ell^{1/2}}{N^{3/2}} \cdot 
     \exp\left( N 
       \Psicom^\star\left( \frac{\ell+1}{N}, \frac{k+1}{\ell+1} 
      \right) \right).
\]
\end{lemma}

The next lemma is more than we really need at this stage; the extra generality 
will be useful in discussion of the weak phase transition in the next subsection.

\begin{lemma} \label{lem:IntegralBounds}
In the proportional growth setting, we have constants $c_2$, $c_3$, and $c_4$ depending
at most on $\delta$,
so that, for $\ell = n+1, n+3, \dots$, $ \ell \leq n + \sqrt{N}$, $k = k_n = \lfloor \rho_\eps n \rfloor$,
and $n > n_0$:
\begin{equation}\label{eq:int_lower}
   \beta(T^k,T^\ell) \geq c_2 \cdot \exp\left( 
-N \Psiint^{\star}\left( \frac{\ell+1}{N}, \frac{k+1}{\ell+1} \right) \right)
\end{equation}
\begin{equation} \label{tltn}
   \alpha(T^\ell,T^{N-1}) \geq c_3 \cdot \exp\left( -N \Psiext^+\left( \frac{\ell+1}{N}, \frac{k+1}{\ell+1} \right) \right)
\end{equation}
\begin{equation} \label{extFlCN}
   \alpha(F^\ell,C^N) \geq c_4 \cdot \exp\left( -N \Psiext^{\pm} \left( \frac{\ell+1}{N}, \frac{k+1}{\ell+1} \right) \right)
\end{equation}
\end{lemma}
Combining the last two lemmas, we get -- specializing to the case $\ell = n+1$\--
\begin{eqnarray*}
 D_{n+1,n} & \geq & c_5 \ell^{1/2} N^{-3/2} \exp\left( N \Psinet^\star\left( \frac{n+2}{N}, \frac{k+1}{n+2} \right) \right) \\
        & \geq & c_5  \ell^{1/2} N^{-3/2} \exp(N \eps )  \goto \infty, 
        \quad N_n \goto\infty .
\end{eqnarray*}
As $f_k(Q) - E f_k(AQ) > D_{n+1,n}^\star$, we conclude that $\Delta(k_n,n,N_n) \goto \infty$
as $n \goto \infty$;
this completes the upper bound
for the strong phase transition $\rho_S^\star$. \qed

\subsection{Upper Bounds on the Weak Phase Transition}

We now aim to show that, in the proportional growth setting with 
$n/N_n \goto \delta > 0$, and $k = \lfloor \rho_W^\star(\delta) n \rfloor$,
\begin{equation} \label{weaktransition}
(f_k(Q) - E f_k(AQ))/f_k(Q) > \eps> 0, \quad n > n_0.
\end{equation}
In words, `above $\rho_W^\star(\delta)$ a nonvanishing fraction of faces 
get lost under projection'.  

In fact we will show that for all large enough $n$, and
for all $\ell$ in the range $n+1, n+3, \dots$, $\ell \leq n + \sqrt{N}$,
\begin{equation} \label{dlbound}
   D_{\ell,n}^\star / f_k(Q) \geq c n^{-1/2}, \quad \ell = n+1, n+3, \dots; \quad \ell \leq n+ \sqrt{N}.
\end{equation}
Since this inequality holds for at least $\frac{1}{2} \sqrt N$ terms from
the sum $\sum'$, we have
\[
    \Delta(k_n,n,N_n) = \Sigma' D_{\ell,n}^\star \geq \eps f_k(Q),
\]
for $\eps = c/2$, which implies (\ref{weaktransition}).

The different structure of our argument in the weak transition case
can be traced to the fact that $\Psinet^\star - \Psiface^\star$
does not change sign at $\rho = \rho_W^\star$. Instead, it achieves
its global maximum $0$. This means that
\[
   \frac{\partial}{\partial \nu} ( \Psinet^* - \Psiface^\star )( \nu , \delta \rho/\nu) = 0
\]
from which it follows that, for $\nu \in [\delta, \delta + 1/\sqrt{N}]$,
and some $c > 0$,
\[
( \Psinet^* - \Psiface^\star )( \nu , \delta \rho/\nu)  \geq c/{N}.
\]

The combinatorial identity
\[
   { {n}\choose{k+1} }  { {n-k+1}\choose{\ell-k} } ={ {n}\choose{\ell} }  { {\ell}\choose{k+1} } .
\]
implies
\[
 D_{\ell,n}^+ / f_k(T^{N-1})  = 2 \cdot { {N-k+1}\choose{\ell-k} }  \cdot \beta(T^k,T^\ell) \alpha(T^\ell,T^{N-1})
\]
and
\[
 D_{\ell,n}^\pm / f_k(C^N)  = 2 \cdot { {N-k+1}\choose{\ell-k} }  \cdot \beta(T^k,T^\ell) 
      \alpha(T^\ell,C^{N}).
\]

We need the following combinatorial result; it follows from Stirling's
inequalities (Lemma \ref{lem:ShannonBounds}) and we omit the proof.

\begin{lemma}
\[
  { {N-k-1}\choose{\ell-k} } \geq \frac{1}{3} (N - k - 1)^{-1/2}
       \exp\left((N-k-1) H\left(\frac{\ell-k}{N-k-1}\right) \right) .
       \]
\end{lemma}

We combine this with Lemma \ref{lem:IntegralBounds}
and get that, under proportional growth
\[
 D_{\ell,n}^\star / f_k(Q) \geq c N^{-1/2} \cdot \exp\left(N ( \Psinet^* - \Psiface^\star ) \left(\frac{\ell+1}{N},\frac{k+1}{\ell+1}\right) \right).
\]
This implies (\ref{dlbound}) and (\ref{weaktransition}) follows.

\subsection{Analysis of the External Angle}\label{sec:ext_lower}

{\it Simplex case.}
We recall the exact formula
\begin{equation} \label{ExtExact}
 \alpha( T^{\ell}, T^{N-1}) = \sqrt{\frac{\ell+1}{\pi}} \int_0^\infty \exp(-N \psi_\nu(x) ) dx ,
\end{equation}
where $\psi_\nu(x) = \nu x^2 - (1-\nu) \log Q(x)$ and $\nu = \frac{\ell+1}{N}$.
Note that throughout this Section, \ref{sec:ext_lower}, we use the 
convention $\nu = \frac{\ell+1}{N}$ consistent with 
Lemma \ref{lem:IntegralBounds}.  In the next subsection,
we use Laplace's method to obtain lower bounds on general integrals of this type.
That lemma requires estimates which are, in turn supplied by Lemma
\ref{lem-x} below.
Here the exponent $\psi_\nu$ should not be viewed as constant in $N$; it depends
on the variable $\nu = \frac{\ell+1}{n}$ which varies slightly; also
$\ell$ is a variable which ranges in the vicinity of $n$.
Lemma \ref{lem-x} gives lower bounds on 
Laplace integrals with uniform multiplicative remainders;
this yields that for a fixed subinterval $0 < \nu_0 < \nu_1 < 1$, and for each 
$\eps > 0$ there is $N_0(\eps)$
so that for all $\nu =  (\ell+1)/N$ in $(\nu_0,\nu_1)$,
\[
 \sqrt{\frac{\ell+1}{\pi}}
  \int_0^\infty \exp(-N \psi_\nu(x) ) dx  \geq \sqrt{\frac{2\pi}{N \psi''_\nu(x_\nu)}}
     \exp(-N \Psi_\nu(x_\nu)) \cdot (1 -\eps).
\]
We conclude from (\ref{ExtExact}) and (\ref{EvalPsipp}) that
for each $\eps > 0$ and all sufficiently large $N > N_0(\nu,\eps)$,
\[
 \alpha( T^{\ell}, T^{N-1}) \geq \sqrt{(1-\nu)(1+2x_\nu^2)} \cdot \exp(-N \Psiext^+(\nu))( 1 - \eps),
\]
where the threshold $N_0$ may be taken 
locally uniform in $\nu \in [0,1)$. Here, again,
$x_\nu$ is the minimizer of $\psi_\nu(x)$.
It follows that there is  
a constant $c > 0$ so that for all sufficiently large $N$,
and all $\nu \in I_N$ , $I_N = [\delta, \delta + 1/\sqrt{N}]$,
\[
 \alpha( T^{\ell}, T^{N-1}) \geq c \cdot \exp(-N \Psiext^+(\nu)).
\]
Equation (\ref{extFlCN}) follows. 

\begin{lemma} \label{lem-x}
Let $\psi_\nu(x) = \nu x^2 - (1-\nu) \log Q(x)$. Then $\psi_\nu(x)$ is $C^4(0,\infty)$,
\begin{equation} \label{EvalPsipp}
  \psi''_\nu(x_\nu) = \frac{2 \nu}{1-\nu} ( 1+ 2 x_\nu^2) .
\end{equation}
\[
  \psi'''_\nu(x_\nu) =  -(1-\nu) \left[\frac{4 \nu}{1-\nu} ( 2 x_\nu^3 - x_\nu)
       + \frac{ 24\nu}{1-\nu} x_\nu^2
       - \frac{16 \nu}{1-\nu} x_\nu^3 \right] .
\]
For $\eps > 0$, set
\[
   C(\nu,\eps) = \sup_{ |x - x_\nu| < \eps} \frac{| \psi'''_\nu(x)|}{ \psi''_\nu(x_\nu) }.
\]
Then for  small $\eps > 0$, $C(\nu,\eps) < \infty$, and as $\nu \goto 0$,
$C(\nu,\eps) \sim 2 x_\nu$.
 \end{lemma}
 
 {\it Cross-polytope case.}
We recall the exact formula
\[
 \alpha( F^{\ell}, C^N) = \sqrt{\frac{\ell+1}{\pi}}  \int_0^\infty \exp(-N \psi_\nu(x)) dx ,
\]
where $\psi_\nu(x) = \nu x^2 - (1-\nu) \log G(x)$ and $\nu = (\ell+1)/N$ 
We apply Lemma \ref{lem-x} bounding Laplace integrals with multiplicative remainder
to conclude
\[
 \sqrt{\frac{\ell+1}{\pi}} \int_0^\infty \exp(-N \psi_\nu(x) ) dx  \geq \sqrt{\frac{2\pi}{N \psi''_\nu(x_\nu)}}
     \exp(-N \psi_\nu(x_\nu) ) \cdot (1 + o(1)).
\]
Here the $o(1)$-term is locally uniform over $\nu \in [0,1)$. We conclude that
for each subinterval $(\nu_0,\nu_1)$ with $0 < \nu_0 < \nu_1 < 1$
and for $\eps > 0$ we have for  $N > N_0(\nu,\eps)$,
\[
 \alpha(F^{\ell}, C^N) \geq (1+\frac{4\nu}{1-\nu} x_\nu^2)^{-1/2} 
    \cdot \exp(-N \Psiext^\pm(\nu))( 1 - \eps).
\]
Here, again,
$x_\nu$ is the minimizer of $\psi_\nu(x)$.
It follows that there is a constant $c > 0$ so that for all sufficiently large $N$
and all $\nu \in I_N$, $I_N = [\delta, \delta + 1/\sqrt{N}]$,

\[
 \alpha(F^{\ell}, C^N) \geq c \cdot \exp(-N \Psiext^\pm(\nu)).
\]
Equation (\ref{tltn}) follows. 

\begin{lemma}
Let $\psi_\nu(x) = \nu x^2 - (1-\nu) \log G(x)$. Then $\psi_\nu(x)$ is $C^4(0,\infty)$,
\[
  \psi''_\nu(x_\nu) = {2 \nu} \cdot ( 1+  x_\nu^2\frac{4 \nu}{1-\nu} ) .
\]
\[
  \psi'''_\nu(x_\nu) =  (1-\nu) \left[\frac{4 \nu}{1-\nu} ( 2  - 4x_\nu^2)
       +  6x_\nu z^2
       + 2 x_\nu z^3 \right],
\]
where $z_\nu = \frac{2\nu x_\nu}{1-\nu}$. For $\eps > 0$, set
\[
   C(\nu,\eps) = sup_{ |x - x_\nu| < \eps} \frac{| \psi'''_\nu(x)|}{ \psi''_\nu(x_\nu) }.
\]
Then for  small $\eps > 0$, $C(\nu,\eps) < \infty$, and as $\nu \goto 0$,
$C(\nu,\eps) \sim 4 x_\nu^3$.
 \end{lemma}
 
\subsection{Uniform Laplace's Method}
 
 We use a uniform variant of Laplace's method, suitable for bounding
 a collection of integrals uniformly. The approach is similar to
\cite{Do05_polytope}.

\begin{lemma} \label{lem:UnifLaplace} 
Let $I = [-\eps,\eps]$ and suppose that $f$ attains its minimum on $I$ at $0$.
Let 
\[
C = \sup_I \frac{ |f'''(x)| }{f''(0)}.
\]
Then
\[
  \int_I exp(-Nf(x)) dx \geq \sqrt{\frac{2\pi}{N f''(0)}} 
    \cdot \exp(-N\psi(0)) \cdot  R(\eps,N) .
\]
\[
   R(\eps,N) =  ( 1 - \sqrt{\frac{2}{\pi}} \exp(-N\eps^2 f''(0)))
                \cdot \exp( - N f''(0) C \eps^3/16) .
\]
\end{lemma}
The derivation of the lemma is similar to that of Lemma A.4
in  \cite{Do05_polytope} (although with all inequalities reversed).

\begin{lemma} \label{lem:Jbound}
Consider the collection of integrals
\[
   J(N,\lambda) = \int_0^\infty \exp(-N f_\lambda (x) ) dx
\]
and suppose either that $\lambda = \lambda_0$ independent of $N$
or that $\lambda = \lambda_N \goto \lambda_0$ as $N \goto \infty$.
Suppose that $f_\lambda$ has a unique minimizer $x_\lambda$ interior to $(0,\infty)$
and suppose that $f_\lambda$ is $C^4(0,\infty)$. Let 
\[
   C(\lambda,\eps) = \sup_{ |x - x_\lambda| < \eps} \frac{| f'''_\lambda(x)|}{ f''_\lambda(x_\lambda) }.
\]
Suppose that
\[
   N f_{\lambda_N} (x_{\lambda_N}) \goto \infty,
\]
and 
\[
   \frac{C(\lambda_N,\eps_N)}{\sqrt{N\psi''(x_{\lambda_N} })} \goto 0.
\]
Then
\[
  J(N,\lambda_N) \geq \sqrt{\frac{2\pi}{Nf_\lambda''(x_\lambda)} }\cdot \exp(-N f_\lambda(x_\lambda)) ( 1+ o(1))
\]
\end{lemma}

To prove Lemma \ref{lem:Jbound}, simply translate coordinates so that
$x_\lambda = 0$, pick $\eps_N = N^{-2/5}$ and set $I = [-\eps,\eps]$,
then apply Lemma \ref{lem:UnifLaplace}.
 
\subsection{Analysis of the Internal Angle}
 
 Our earlier analysis of the internal angle employed an upper bound
 derived in  \cite{Do05_polytope} from large-deviations theory.  We now
 develop a lower bound using complex analysis techniques;
our analysis is related to the approach of Vershik and Sporyshev \cite{VerSpor}.

Let $X \sim HN(0,1)$ be a real half-normal random variable,
i.e. $X = |Z|$ where $Z$ is standard normal. The moment generating function
$M(t) = E e^{tX}$ can be continued to the complex plane.
We have the explicit formula $M(t) = e^{t^2/2} \cdot 2 \Phi(t)$,
where $\Phi$ denotes the standard $N(0,1)$ cumulative distribution
function. Operations with Taylor series show that for $\omega$ real,
$\Phi(i \omega)$ has real part $1/2$ along the imaginary axis
 and so the cumulant generating function
 $\log(2\Phi(z))$ can be consistently defined in a neighborhood
 of both the real and imaginary axes. Define
 \[
    \psi_\gamma(z) =  z^2/2 + (1-\gamma)  \log( 2\Phi(z)).
 \]
 We begin by justifying our interest in the complex domain:
 \begin{lemma} \label{lem:AngleIdent}
For $\gamma = \frac{k+1}{\ell+2}$,
 \[
   \beta(T^k,T^\ell) = \sqrt{\ell+3} \cdot 2^{-\ell-k+1} \cdot \frac{1}{\sqrt{2\pi}}
    \int_{-i\infty}^{i \infty} e^{(\ell+2) \psi_\gamma(z)} dz .
 \] 
 \end{lemma}
 
\begin{figure}[h]
\begin{center}
\begin{tabular}{cc} 
\includegraphics[width=2.35in]{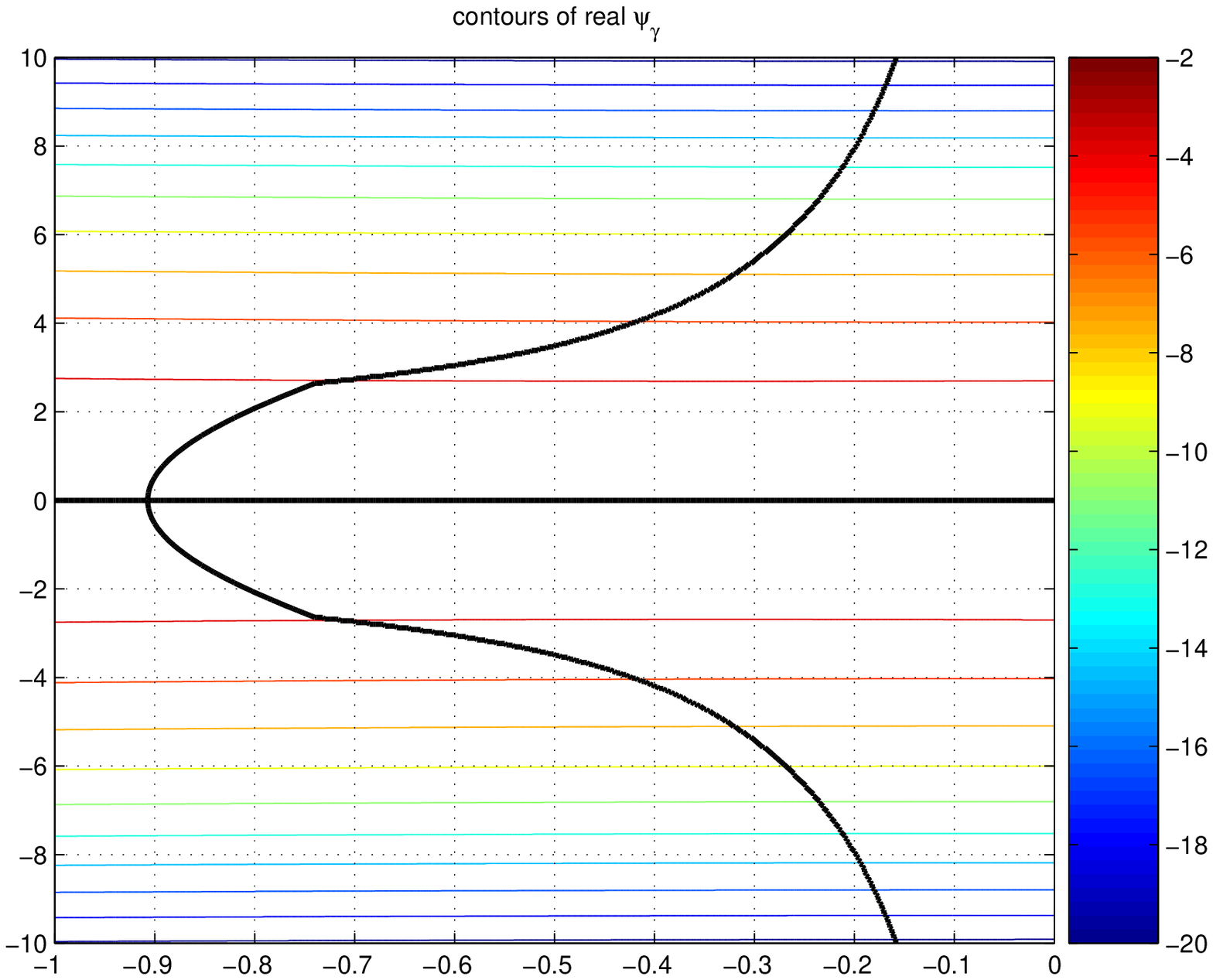} &
\includegraphics[width=2.35in]{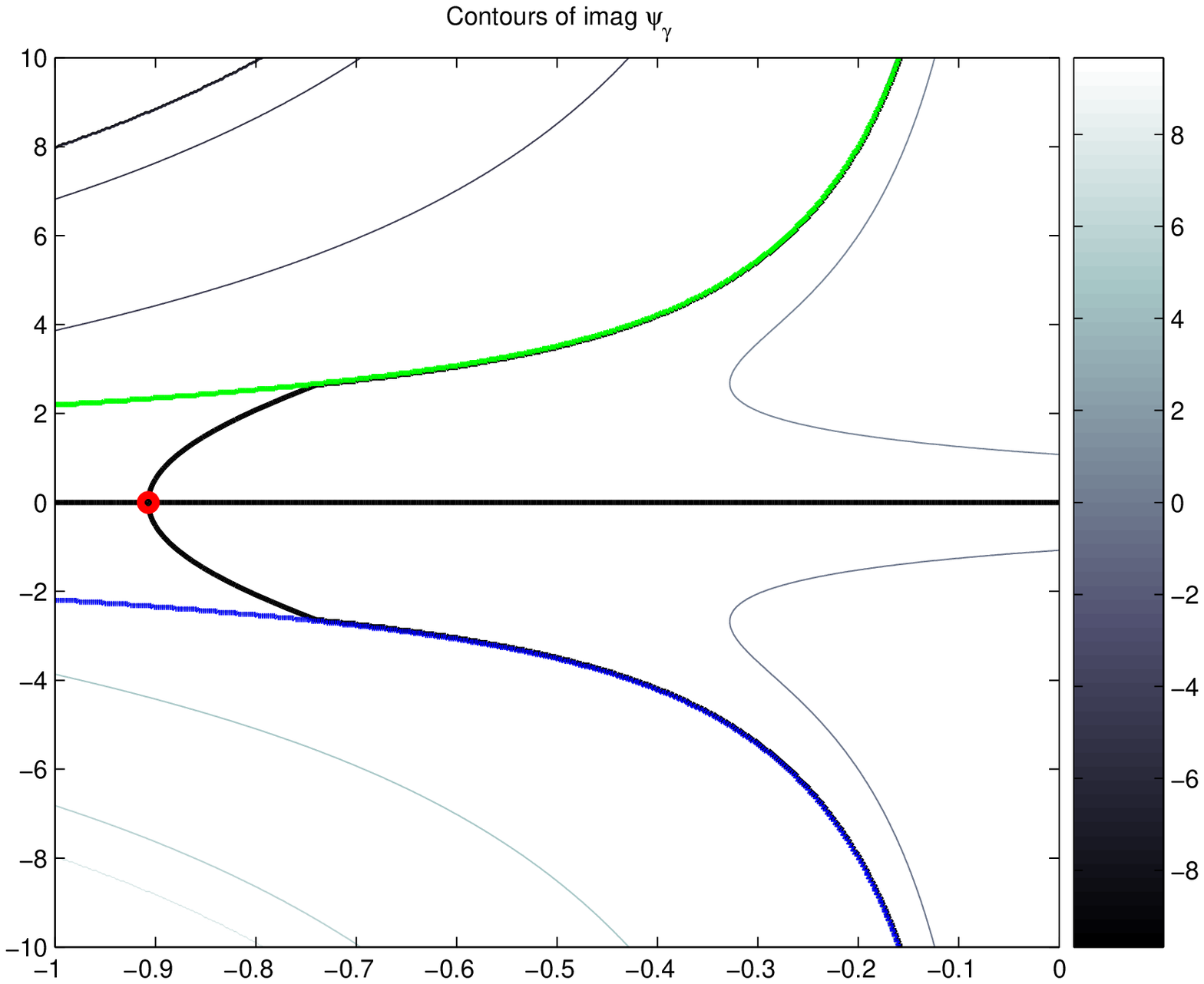} \\
(a) & (b) 
\end{tabular}
\end{center}
\caption{Level curves for the exponent $\psi_\gamma(z)$ with $\gamma=3/8$; real (a) and imaginary (b) components.  The path $C_\gamma$ along which $\psi_\gamma(z)$ is purely real is overlaid as the wider black line.  Panel (b) additionally overlays the level curves with the imaginary part of $\psi_\gamma(z)$ being equal to $-\pi$ (green) and $\pi$ (blue).  The path $C_\gamma$ lies between these hyperbolae and the imaginary axis, allowing the deformation in (\ref{eq:deform}) without necessitating branch cuts.  The {\em saddlepoint}, $z_{3/8}\approx - 0.907 + i0$, is indicated by the red circle in Panel (b).}
\label{fig:Saddle}
\end{figure}

Contour integration was previously used in the analysis of the internal angle
by Vershik and Sporyshev, without 
making the connection to the cumulant generating function.
The contour integral and the form of the integrand
suggests to use the method of steepest descents  \cite{BleisteinHandelsman}.
An analysis of $\psi_\gamma(z)$ is easily performed computationally.
One learns that  there is a path $C_\gamma$ along which $\psi_\gamma(z)$ is 
purely real and which is asymptotic, for large $|z|$, to the imaginary axis; 
see Figure \ref{fig:Saddle}.  
This path crosses the real axis at a point $z_\gamma$. Because $\psi_\gamma$
is real for real $z$, $z_\gamma$ is necessarily a saddlepoint of $\psi_\gamma$.
Within the region bounded by the imaginary axis and $C_\gamma$,
$\psi_\gamma$ is analytic, and so we have the identity
\begin{equation}\label{eq:deform}
 \int_{-i\infty}^{i \infty} e^{(\ell+2) \psi_\gamma(z)} dz  =  \int_{C_\gamma} e^{(\ell+2) \psi_\gamma(z)} dz ,
\end{equation}
provided the orientation of the path $C_\gamma$ is chosen properly.
Parametrizing by arclength, the contour integral can be rewritten
purely in terms of real variables:
\[
   \int_{-\infty}^\infty  e^{(\ell+2) \tilde{\psi}_\gamma(t)}  dt
\]
where $\tilde{\psi}(t) = \psi_\gamma(z(t))$; this of course is
in the form of a Laplace integral. Taking into account that
\[
  \tilde{\psi}_\gamma(0) = \psi_\gamma(z_\gamma), \qquad \tilde{\psi}_\gamma''(0) = \psi_\gamma''(z_\gamma),
\]
and that $\tilde{\psi}_\gamma(t)$ is $C^4(-\infty,\infty)$,
we immediately have:

\begin{lemma}
Let $\gamma$ be fixed in $(0,1)$. $\psi_\gamma(z)$ has a saddlepoint
$z_\gamma$ on the negative real axis and
\[
 \int_{-i\infty}^{i \infty} e^{(\ell+2) \psi_\gamma(z)} dz  = \sqrt{\frac{2\pi}{(\ell+2) \psi''_\gamma(z_\gamma)}}
   \cdot \exp\{(\ell+2) \psi_\gamma(z_\gamma)\} \cdot (1 + o(1)), \quad \ell \goto \infty.
\]
\end{lemma}

Actually, however, we are interested in the case
where $\gamma$ is changing slightly with $n$, i.e. $\gamma = \gamma_n = \frac{k_n + 1}{\ell_n+2}$,
and need a stronger result. We note that  the third and fourth derivatives
of $\tilde{\psi}_\gamma(t)$ near $t=0$ are bounded locally uniformly in $\gamma$. 
We conclude:

\begin{lemma}
Fix $c > 0$.
Let $\gamma_n = \frac{k_n + 1}{\ell_n+2}$.  In the proportional growth setting, we have
\[
 \int_{-i\infty}^{i \infty} e^{(\ell_n+2) \psi_{\gamma_n}(z)} dz  = \sqrt{\frac{2\pi}{(\ell+2) \psi''_{\gamma_n}(z_{\gamma_n})}}
   \exp\{(\ell_n+2) \psi_{\gamma_n}(z_{\gamma_n})\} \cdot (1 + o(1)), \quad n \goto \infty,
\]
with the term $o(1)$ uniform in $n \leq \ell_n \leq n + c \sqrt{n}$.
\end{lemma}

To complete the evaluation of the asymptotics of the internal angle
we need
\begin{lemma} \label{lem:psicalc}
\[
   \psi''_\gamma(z_\gamma) = 1 - z_\gamma^2 \cdot \frac{\gamma}{1-\gamma}
\]
Let $\xi_\gamma(y)$ denote the function introduced earlier in connection
with the internal angle. Then
\[
    \psi_\gamma(z_\gamma) = -(1-\gamma) \cdot \xi_{\gamma}(y_{\gamma}).
\]
\end{lemma}
We conclude that 
\begin{equation} 
   \beta(T^k,T^\ell) \geq c_2 \cdot \exp\left( -N \Psiint^{\star}\left( \frac{\ell+1}{N}, \frac{k+1}{\ell+1} \right) \right).
\end{equation}
The result (\ref{eq:int_lower}) follows.

\subsection{Proof of Lemmas \ref{lem:AngleIdent} and \ref{lem:psicalc}}

\subsubsection{Proof of Lemma \ref{lem:AngleIdent}}

B\"or\"ozcky and Henk gave the formula
\begin{equation} \label{star}
  \beta(T^k,T^\ell) = \theta^{(m-1)/2} \cdot \sqrt{(m-1)\alpha + 1} \cdot \pi^{-m/2} \cdot \alpha^{-1/2} \cdot J(m,\theta),
\end{equation}
where 
\begin{equation} \label{starstar}
\theta = k+1, \quad \alpha = 1/(k+2), \quad m = \ell -k +1,
\end{equation}
and
\[
J(m,\theta) = \frac{1}{\sqrt{\pi}} \int_{-\infty}^\infty e^{-\lambda^2} 
   \left(\int_0^\infty \exp(-\theta v^2 + 2i v\lambda) dv \right)^m d\lambda .
\]
Note that
\[
  \int_0^\infty \exp(-\theta v^2 + 2i v\lambda) dv  
   = \frac{\sqrt{\pi}}{2 \sqrt{\theta}} \cdot E e^{i \sqrt{\frac{2}{\theta}} \lambda X},
\]
where $X$ is standard half normal, $X = |Z|$, $Z \sim N(0,1)$.
Using now the cumulant generating function of the half-normal,
\[
  \Lambda(z) = \log E e^{zX} ,
\]
we write
\[
  J(m,\theta) = \frac{\pi^{m/2-1/2}}{2^m \theta^{m/2}} \cdot  
   \int_{-\infty}^\infty e^{(i\lambda)^2} e^{m \Lambda(i \sqrt{\frac{2}{\theta}} \lambda)} d \lambda .
\]
Now change variables $\omega = \sqrt{\frac{2}{\theta}} \lambda$, and write
\[
  J(m,\theta) = \frac{\pi^{m/2-1/2}}{2^m \theta^{m/2}} \cdot  \sqrt{\frac{\theta}{2}}
   \int_{-\infty}^\infty e^{\theta  (i\omega)^2/2} e^{m \Lambda(i \omega)} d \omega .
\]
Recalling (\ref{star})-(\ref{starstar}) and noting that
\[
  \theta = k+1, \qquad ((m+1)\alpha + 1)/\alpha = \ell+3,
\]
we have 
\[
 \beta(T^k,T^\ell) = \sqrt{\ell+3} \cdot 2^{-m} \cdot \frac{1}{\sqrt{2\pi}} \cdot 
    \int_{-\infty}^\infty e^{-(k+1) \omega^2/2 + (\ell-k+1) \Lambda(i\omega)} d \omega.
\] 
The cumulant generating function of the half-normal
obeys $\Lambda(s) = e^{s^2/2} \cdot 2 \Phi(s)$. Setting
$\gamma = \frac{k+1}{\ell+2}$ the exponent
can be rewritten as $(\ell+2) \psi_{\gamma}(z)$. \qed

\subsubsection{Proof of Lemma \ref{lem:psicalc}}

Note that
\[
  \Lambda'(z) = \phi(z)/\Phi(z),
\]
where $\psi$ is the standard normal density and $\Phi$ is the
standard normal cumulative distribution function.
Hence from $\psi'_\gamma(z) = z + (1-\gamma) \phi(z)/\Phi(z)$
and $\psi_\gamma(z_\gamma) = 0$ we have
\begin{equation} \label{zgammaident}
  -z_\gamma/(1-\gamma) = \phi(z_\gamma)/\Phi(z_\gamma).
\end{equation}

We also have
\[
  \psi''_\gamma(z) = 1 + (1-\gamma) [ \frac{\phi'}{\Phi} - \frac{\phi^2}{\Phi^2}]
\]
and $\phi' = (-z)\phi$. Hence
\[
 \psi''_\gamma(z_\gamma) = 1 - z_\gamma^2 \cdot \frac{\gamma}{1-\gamma} .
\]
This proves half the lemma.

For the other half of the lemma, we need to establish a connection between
the values of $\psi_\gamma(z_\gamma) =  z_\gamma^2/2 + (1-\gamma) \Lambda(z_\gamma)$ and
$\xi_\gamma(y_\gamma)$, where
\[
   \xi_\gamma(y) = \frac{1-\gamma}{\gamma} y^2/2 + \Lambda^*(y) .
\]
Here $\Lambda^*(x) = \max_{s} sx - \Lambda(s)$ is the 
the classical Fenchel-Legendre transform
of cumulant generating function on the real axis.
It is worth reviewing Sections 6.4 and 6.5 of \cite{Do05_polytope}.
The definition of $\Lambda^*$ sets up a one-one relationship between
variables $(y,s)$, where $y = y(s)$ and $s = s(y)$, where
\[
    \Lambda^*(y) = s(y)y - \Lambda(s).
\]
Hence 
\[
   \xi_\gamma(y_\gamma) = s(y_\gamma)y_\gamma - \Lambda(s_\gamma) +  \frac{1-\gamma}{\gamma} y_\gamma^2/2 .
   \]
Formula (6.12) in \cite{Do05_polytope} reads
\[
 \frac{1-\gamma}{\gamma} y_\gamma = - s_\gamma;
\]
this implies 
\begin{eqnarray*}
 \xi_\gamma(y_\gamma) &=&  - \frac{\gamma}{1-\gamma} s_\gamma^2/2 - \Lambda(s_\gamma) \\
           &=& - \frac{\gamma}{1-\gamma} s_\gamma^2/2 -  s_\gamma^2/2 - \log(2 \Phi(s_\gamma))\\
           &=& - \frac{1}{1-\gamma} s_\gamma^2/2  - \log(2 \Phi(s_\gamma)).
\end{eqnarray*}
We note \-- parenthetically \-- that the variable $s$ is in this subsection the argument to a cumulant generating function, and elsewhere in the paper, the same symbol denotes the negative of this
same quantity.
Moreover the dual relationship between $s,y$ variables
is expressed through $\Lambda'(s_\gamma) = y_\gamma$. 
We compute that $\psi_\gamma(s_\gamma + i0) = 0$, i.e. $z_\gamma = s_\gamma + i0$.
In words, the saddlepoint value $z_\gamma$ is identical to the dual variable $s_\gamma$.
Finally we have 
\[
-\Psi_{int}(\nu,\gamma) =  -( \xi_{\gamma}(y_{\gamma})  + \log_e(2))\cdot  \nu \cdot (1-\gamma) = 
                                             (\psi_\gamma(z_\gamma) + \log_e(2)(1-\gamma)) \cdot \nu .
\]
Compare also section 6.5 of \cite{Do05_polytope}. \qed

\section{Discussion}
\setcounter{equation}{0}
\setcounter{table}{0}
\setcounter{figure}{0}

In this section, we first show how the applications  (1.1)-(1.4)
follow from Theorems \ref{thm:weakp} - \ref{thm:5}.  We next consider the 
performance of these rules at finite $n$.  Finally we discuss
extensions, open questions, and relations to other work.

\subsection{Convex Hulls of Gaussian Point Clouds. Proof of (1.1)}
\label{sec:AppGaussianProof}

In the 1950's, David Gale \cite{Gale}
introduced an important extremal property
of polytopes; the following is now classical:
\begin{definition} \cite[Chapter 7]{Gru67}
A convex polytope $P$ 
is called {\bf $k$-neighborly} if every subset of $k+1$ vertices
spans a $k$-face of $P$.
\end{definition}

By mere face counting, we can determine
whether a polytope is $k$-neighborly.
In this section, put for short $T = T^{N-1}$.

\begin{lemma} \cite[Chapter 7]{Gru67} 
Let $P = AT$. Suppose that
\bitem
 \item $P$ has $N$ vertices.
 \item $P$ has $ { N \choose k+1 } $ $k$-faces.
\eitem
Then $P$ is $k$-neighborly.
\end{lemma}

Combining Theorems \ref{thm:strongp} and \ref{thm:5} we have:
\begin{corollary}
Let $(k_n,n,N_n)$ be a sequence of triples with $n$ tending to $\infty$,
and $N_n$ growing subexponentially with $n$. Fix $\eps > 0$
and suppose that
\[
   k_n < (1-\eps) \cdot  \frac{n}{2e \cdot \log(N_n/(n\cdot 2\sqrt{\pi}))}, 
\quad n > n_0 .
\]
Let $A$ be a random $n \times N$
matrix with iid $N(0,1/n)$ entries. 
Define the event
\[
\Omega(k,n,N)=\{P=AT \mbox{ is $k-$neighborly} \}.
\]
Then
\[
P(\Omega(k_n,n,N_n))\goto 1,\quad \mbox{as} \quad n\goto\infty.
\]
\end{corollary}

In words, with overwhelming probability for large $n$, 
$P = AT$ is at least $k_n$-neighborly.

This is simply (1.1) in another language.
To see why,  note that, for each $k > 1$, a $k$-neighborly polytope
is also $k-1$-neighborly. 
If $a_1$, \dots , $a_N$ are
vertices of $\cA = \conv(a_1, \dots a_N)$, then
$k$-neighborliness
of $\cA$  
is equivalent to the following $k$ 
simultaneous properties:
\bitem
 \item  every pair $(a_i, a_j)$ spans an edge of $\cA$
 \item ...
 \item every $k+1$-tuple $(a_{i_1}, \dots, a_{i_k})$
 spans a $k$-face of $\cA$.
\eitem
This is precisely the condition mentioned in Section 1.1.1 
with the substitutions: $x_i\leftrightarrow a_i$, 
$\cA\leftrightarrow\cX$, $n\leftrightarrow d$, and $N\leftrightarrow n$.

To conclude, we note that $P = AT$ has
$N$ vertices with probability 1, those vertices are
simply the columns of $A$, and so $P = \conv(a_1,\dots, a_N)$.
Invoking now the above corollary we obtain the conclusion (1.1). \qed

\subsection{Correcting all patterns of $k$ or fewer errors. Proof of (1.4)}
\label{sec:AppAllECCProof}

A convex polytope is centrosymmetric
if it has $2N$ vertices made of $N$ antipodal pairs.
Neighborliness {\it per se} does not apply to centrosymmetric polytopes,
instead one needs the following notion:
see eg \cite[Chapter 8]{Gru67}.
\begin{definition}
A centrosymmetric convex  polytope $P$ with vertices $\pm a_1, \dots, \pm a_N$
is called {\bf centrally $k$-neighborly} if every subset of $k+1$ vertices
{\it not including an antipodal pair} spans a $k$-face of $P$.
\end{definition}
By face counting, we can determine
whether a polytope is centrally $k$-neighborly.
In this section, put for short $C = C^N$.

\begin{lemma}  \cite[Lemma 1]{Do05_polytope}
Let $P = AC$. Suppose that
\bitem
 \item $P$ has $2N$ vertices; and
 \item $P$ has $2^{k+1} \cdot  { N \choose k+1 } $ $k$-faces.
\eitem
Then $P$ is centrally $k$-neighborly.
\end{lemma}

Combining Theorems \ref{thm:strongpm} and \ref{thm:5} we have
\begin{corollary}
Let $(k_n,n,N_n)$ be a sequence of triples with $n$ tending to $\infty$,
and $N_n$ growing subexponentially with $n$. Fix $\eps > 0$
and suppose that
\[
   k_n < (1-\eps) \cdot  \frac{n}{2e \cdot \log(N_n/(n\cdot \sqrt{\pi}))} , \quad n > n_0 .
\]
Let $A$ be a random $n \times N$
matrix with iid $N(0,1/n)$ entries. 
Then $P=AC$ is a random centrosymmetric polytope.  
Define the event
\[
\Omega(k,n,N)=\{P=AC \mbox{ is centrally $k-$neighborly} \}.
\]
Then
\[
P(\Omega(k_n,n,N_n))\goto 1,\quad \mbox{as}\quad n\goto\infty.
\]
\end{corollary}

In words, with overwhelming probability for large $n$, 
$P = AC$ is at least $k_n$-centrally neighborly.

We now relate central $k$-neighborliness
to (1.4).
Recall the optimization problem
\[
 (P_1)  \qquad \min_x \| x \|_1 \mbox{ subject to } y = Ax.
\]
Call the solution $x_1$; it obviously depends on $y$ and $A$.

\begin{theorem}   \cite{Do05_signal}
\label{thm:l1l0equiv}
The following statements about an $n \times N$ matrix $A$
are equivalent.
\bitem
\item The polytope $AC$ has $2N$ vertices and is centrally $k$-neighborly.
\item For every problem instance $y = A x_0$ where
$x_0 \in \bR^N $ has at most $k$ nonzeros, the solution
$x_1$ to the corresponding instance of $(P_1)$ is unique
and is equal to $x_0$.
\eitem
\end{theorem}

To apply this, recall the setting of Section \ref{int:ecc}.
The encoding matrix $B$ mentioned there was
obtained as follows: a random orthogonal matrix $U$
is generated, and $B$ makes up $N-n$ rows
of this matrix. The checksum matrix $A$ makes up the other $n$
rows of $U$.

Given received data $w \in \bR^N$,
form the generalized checksum $y = A w \in \bR^{n}$.
Then solve the instance of $(P_1)$
defined by $(y,A)$.  Define the  reconstruction
$u_1 = B(w-x_1)$.  (1.4) now follows from the above,
and the following:

\noindent
{\bf Claim.}  If $AC$ is centrally $k$-neighborly,
and if the error vector $z$ has
at most $k$ nonzeros,  one has
perfect error-correction:
\[
u = u_1.
\]

\noindent
{\bf Proof.} The received message $w = B^T u + z$ 
where the error vector $z$ has, by hypothesis,
nonzeros in at most $k$ positions. Since $A B^T = 0$,
$y = Az$.  Invoking Theorem \ref{thm:l1l0equiv}, we have
$x_1 = z$. Hence $B(w-x_1) = B(w-z) = BB^T u + z -z  = u$. \qed


%
%
  
\subsection{How Many Projections? Proof of (1.2)}
\label{sec:AppHowManyProof}

We first transform the ``how many questions''  problem into
face counting. 

\begin{definition}
The random $n \times N$ matrix $A$ will be called {\it orthant-symmetric}
if, for every signed permutation $\Pi$, and for 
every measurable $\Omega \subset \bR^{n\times N}$,
\[
   P \{ A \in \Omega \} = P \{ A\Pi \in \Omega \}.
\]
\end{definition}

\begin{theorem}  \cite{Do05_signal}
\label{l1l0equiv}
Let $A$ be an orthant symmetric random $n \times N$ matrix.
Let $x_0$ be a fixed vector with $k$ nonzeros.
Form a random problem instance $(y,A)$ of $(P_1)$,
where $y = Ax_0$.  Let $x_1$ denote the solution of this
instance of $(P_1)$.
\[
 P \{ x_1 = x_0 \} \geq \frac{ E f_{k-1}(AC)}{f_{k-1}(C)} .
\]
\end{theorem}

Theorems \ref{thm:weakpm} and \ref{thm:5} 
imply the following precise version of (1.2).

\begin{corollary} \label{lem:expfaces}
Let $(k_n,n,N_n)$ be a sequence of triples with $n$ tending to $\infty$,
and $N_n$ growing subexponentially with $n$. Fix $\eps > 0$
and suppose that
\begin{equation}\label{indicasymp}
   k_n < (1-\eps) \cdot  \frac{n}{2 \cdot \log( N_n/n)} , \quad n > n_0 .
\end{equation}
Then 
\[
    \frac{ E f_{k-1}(AC)}{f_{k-1}(C)} \goto 1 , \quad n \goto \infty.
\]
\end{corollary}

In words, for $(k,n,N)$ obeying the asymptotics (\ref{indicasymp}),
an overwhelming fraction of the $k-1$ faces $F$ of
$C$ induce $k-1$ faces $AF$ of $AC$.

$P = AT$ is at least $k_n$-neighborly.

\subsection{Correcting random patterns of $k$ errors or fewer. Proof of (1.3)}
\label{sec:AppECCProof}

Let $\| z \|_0$ count the number of nonzeros in $z$.
\begin{definition}
The random vector $z$ is a {\it symmetric $k$-sparse random vector} if 
\bitem
\item $P \{ z \in \Omega \} = P \{ -z \in \Omega \} $ for all measurable sets $\Omega$; and 
\item $P \{ \| z \|_0  \leq k \} = 1$.
\eitem
\end{definition}

Suppose that the received message $w= B^T u + z$
where $u$ is arbitrary and $z$ is a symmetric $k$-sparse
random vector stochastically independent of $A,B$.
Define $y = A\mu$ and consider the resulting instance
of $(P_1)$.  Then, conditional on each 
fixed realization of $z$, put $x_0 := z$
and apply  Theorem \ref{l1l0equiv} to get
that
\[
    E \{ f_{k-1}(AC)  | z \} \geq (1-\eps) f_{k-1}(C).
\]
implies
\[
    P \{ x_1 = z | z \} \geq 1-\eps. 
\]
By independence of $z$ and $A$,
\[
  E \{ f_{k-1}(AC)  | z \} = E f_{k-1}(AC).
\]
Apply now Corollary \ref{lem:expfaces} to infer
(1.3) \qed

\subsection{Empirical Results}\label{int:empirical}

The phenomena uncovered by Theorems \ref{thm:weakp} and \ref{thm:weakpm}
can be observed empirically.  For a given $(\delta,\rho)$ pair,
pick a large $N$, generate a random $A$ of dimensions $n = \lfloor \delta N \rfloor $
by $N$,
and check whether for $k = \lfloor \rho \cdot n\rfloor $, a randomly chosen
$k$-face $F$ of $Q$, yields a
projected simplex $AF$ 
that is also a face of $AQ$; here $Q=T^{N-1}$ or $Q=C^{N}$. 
This can be verified by linear programming.

Let $F$ be a $k-1$-face of $Q = T^{N-1}$ or $C^N$.
Then the elements of $F$ have nonzeros in only $k$ coordinates.
If $Q = T^{N-1}$ the nonzeros are nonnegative;
if $Q = C^N$ the nonzeros have a definite sign pattern
particular to the interior of $F$. In the following result,
let $\chi_F$ denote the barycenter of the face $F$,
and let $(\cP)$ denote problem $(P_1)$ if $Q= C^N$ or 
problem $(LP)$ if $Q = T^{N-1}$, where
\[
    (LP)  \quad\quad \min 1'x  \mbox{ subject to }  y = Ax, \; x \geq 0.
\]

\begin{theorem}  \cite{Do05_signal,DoTa05_signal}
Let $Q = T^{N-1}$ or $C^N$.
The following statements about a face $F$ of $Q$
are equivalent.
\bitem
 \item $AF$ is a face of $AQ$.
 \item Let  $y_F = A \chi_F$. 
 Then $\chi_F$ is the unique solution
 of the instance of $(\cP)$ defined by $(y_F,A)$.
\eitem
\end{theorem}

Thus, to check Theorems \ref{thm:weakp} and \ref{thm:weakpm},
one checks that for a randomly-generated
vector $x_0$ with $k$ 
nonzeros, the corresponding
vector $y = Ax_0$ generates an instance of
either $(LP)$ or $(P_1)$ uniquely solved by $x_0$;
$(LP)$ corresponds to $x_0\ge 0$ and $Q= T^{N-1}$
whereas $(P_1)$ corresponds to $x_0$ with entries 
of either sign and $Q=C^N$.
If such uniqueness holds,
we call that experiment a success. Theorems \ref{thm:weakp} 
and \ref{thm:weakpm} imply that for $k$
below a given threshold, success is very likely while above
that threshold, success is very unlikely.

We conducted $44,000$ such
experiments with the common value $N=10,000$, 
exploring the $(\delta,\rho)$ domain as follows.  
We considered $n=10,15,20,\ldots,100$; 
for each value of $n$, eleven values of the sparsity, $k$, 
were chosen near the asymptotic thresholds, $n\cdot |2\log (\delta)|^{-1}$.  
At each combination of $k$ and $n$, two hundred random problem
instances were generated. 

\begin{figure}
\begin{center}
\includegraphics[height=2.5in]{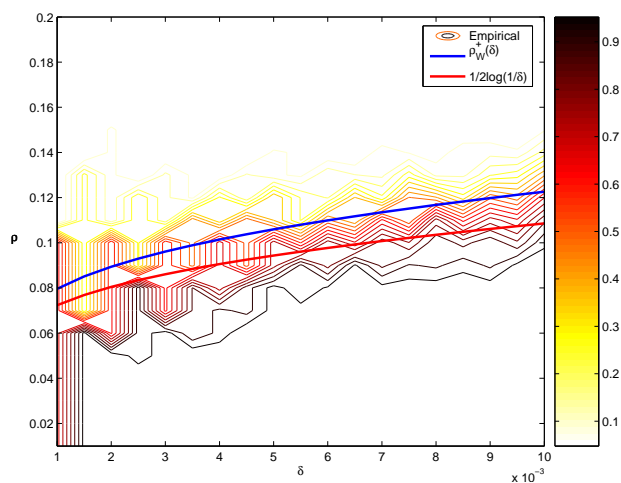} 

(a)

\includegraphics[height=2.5in]{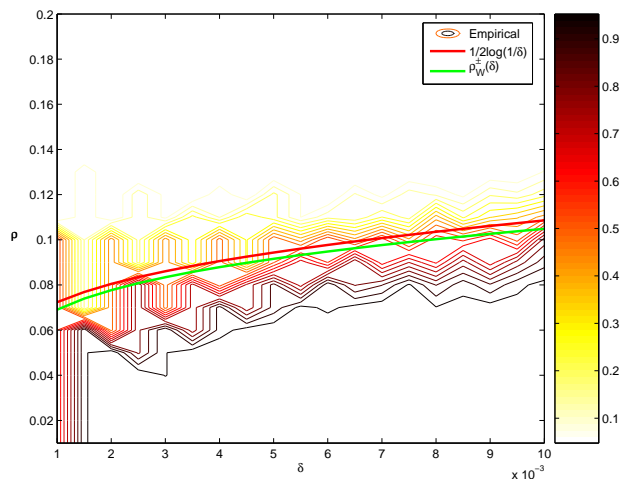} 

(b)

\includegraphics[height=2.5in]{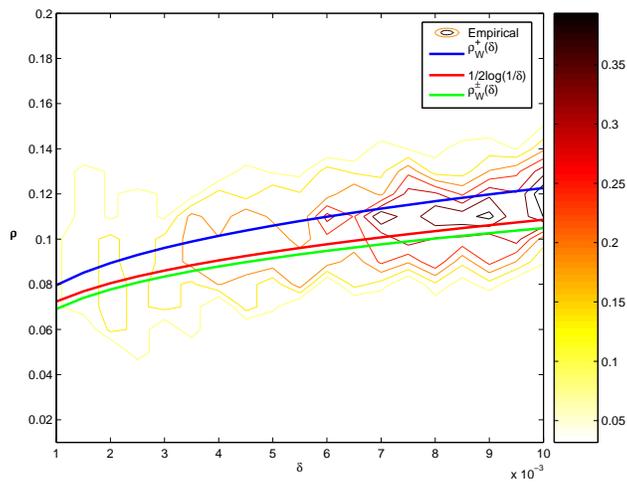} 

(c)

\end{center}
\caption{Panel (a) Success fraction, simplex; 
Panel (b) Success fraction, cross-polytope;
Panel (c) Fraction successful for Simplex but not for Cross-Polytope.   
$N=10,000$.}
\label{fig:emp_region}

\end{figure}

Figures \ref{fig:emp_region} (a)-(c) summarize our results.  
A region of the $(\delta,\rho)$ plane is decorated with
a shaded attribute depicting the fraction of successful  
experiments.
Figure \ref{fig:emp_region}(a) shows the simplex case, 
along with the threshold $\rho_W^+(\delta)$ and its 
asymptotic approximant, $|2\log (\delta)|^{-1}$.  
Figure \ref{fig:emp_region}(b) shows the cross-polytope case, 
with the threshold $\rho_W^{\pm}(\delta)$ and the 
approximant, $|2\log (\delta)|^{-1}$.  
To better highlight the (subtle) difference between the 
simplex and cross-polytope cases, 
Figure \ref{fig:emp_region}(c) shows the fraction of cases where 
the simplex experiments were successful and the 
cross-polytope experiments were not. 

Figure \ref{fig:emp_region} (a) and (b) display a remarkable match 
between the thresholds $\rho_W^+(\delta)$, $\rho_W^\pm(\delta)$ 
and their asymptotic approximations, $|2\log (\delta)|^{-1}$.  
Both curves track  the  observed
empirical phase transition.   This empirical transition
is of course not a true discontinuity, because we are working
with finite problem size $N=10,000$; instead
it is a relatively abrupt change. 
Still, some relatively sharp distinctions
can be made; there is a definite region
where   the simplex experiment is typically successful 
 but the cross-polytope experiment is 
not \-- see Figure \ref{fig:emp_region} (c).  
 
 For $\delta$ near $1/100$ the empirical transitions
 at $N=10,000$ show a clear agreement 
 with the theoretical thresholds $\rho_W(\delta)$ 
and the $|2\log (\delta)|^{-1}$  asymptotic approximant.  
Fixing the region $\delta\in[1/1000,1/100]$ 
explored in Figure \ref{fig:emp_region} and increasing $N$ offers better 
resolution in $k/n$; the sharper empirical transition 
is again in agreement with the theoretical  thresholds. 

\subsection{Exponentiality}\label{subsec:exp}

A key element of our proofs, which we have not emphasized
in the formal statement of our theorems, is that the 
rate of approach to all the limits of interest is exponential
in the proportional growth case.  Thus, we have shown
that, for $\rho < \rho_S^\star(\delta)$, there are positive numbers
$\psi_i(\star,\delta,\rho)$ yielding
\[
    f_{k}(Q)  - Ef_{k}(AQ) \} \leq N^{\psi_1} \exp(- \psi_2 N ).
\]
In all our proofs can be found explicit calculations
of these exponents and remainders.
Similarly,
we have shown
that, for $\rho < \rho_W^\star(\delta)$, there are positive numbers
$\psi_i(\star,\delta,\rho)$ yielding
\[
    \frac{E f_{k}(AQ)}{f_{k}(Q)}  \geq 1 - N^{\psi_1} \exp(- \psi_2 N ).
\]
In the non-proportional growth case, analogous relations
hold, provided we consider triples $(k_n,n,N_n)$
along a trajectory $k_n = \lfloor r_S^\star(n/N_n; \tau) \cdot n \rfloor$,
with $\tau >  2e$ or $k_n = \lfloor r_W^\star(n/N_n; \tau) \cdot n \rfloor$,
with $\tau >  2$, Theorem \ref{thm:5}.  

In fact our results are
strong enough to yield explicit values effective at moderate $(k,n,N)$.
The following two Theorems
follow directly from equations 
(\ref{eq:psicom_precise}), (\ref{eq:psiext_precise_p}), 
(\ref{eq:psiext_precise_pm}), (\ref{IntGoal1}), 
(\ref{eq:psiface_precise_p_upper})
and (\ref{eq:psiface_precise_pm}).
The notations $\tilde{\nu}$, $\hat{\nu}$, etc.
are defined in those equations.

\begin{theorem}\label{thm:strong_precise}
\[
f_{k}(T^{N-1})  - E f_{k}(AT^{N-1}) 
< (N+3)^5 \exp\left(N\left[(\Psi_{com}^+-\Psi_{ext}^+)(\nu,\gamma)
-\Psi_{int}^+(\tilde{\nu},\tilde{\gamma})\right]\right)
\]
and 
\[
f_{k}(C^{N})  - E f_{k}(AC^{N}) 
< (N+3)^5 \exp\left(N\left[\Psi_{com}^{\pm}(\nu,\gamma)
-\Psi_{ext}^{\pm}(\widehat{\nu})
-\Psi_{int}^{\pm}(\tilde{\nu},\tilde{\gamma})\right]\right),
\]
each uniformly over $k=1,2,\ldots n$ and $n=1,2,\ldots N-1$.
\end{theorem}

\begin{theorem}\label{thm:weak_precise}
\[
\frac{E f_{k}(AT^{N-1})}{f_{k}(T^{N-1})}  \geq 1 -
(N+3)^{11/2} \exp\left(N\left[(\Psi_{com}^{+}
-\Psi_{ext}^{+}
-\Psi_{face}^{+})(\nu,\gamma)
-\Psi_{int}^{+}(\tilde{\nu},\tilde{\gamma})
\right]\right)
\]
and
\[
\frac{E f_{k}(AC^{N})}{f_{k}(C^{N})}  \geq 1 -
(N+3)^{11/2} \exp\left(N\left[(\Psi_{com}^{\pm}
-\Psi_{face}^{\pm})(\nu,\gamma)
-\Psi_{ext}^{\pm}(\widehat{\nu})
-\Psi_{int}^{\pm}(\tilde{\nu},\tilde{\gamma})
\right]\right),
\]
each uniformly over $k=1,2,\ldots n$ and $n=1,2,\ldots N-1$.
\end{theorem}

Equipped with these explicit bounds,
we can make nonasymptotic bounds
answering a variety of interesting
questions:
\begin{description}
\item[Q1 (Setting \ref{int:Hull})] At a particular choice of  $N$ and $n$,
for what values of $k$ is there a positive
chance that a standard Gaussian
point cloud has a $k$-neighborly convex hull?
Similarly, consider
the symmetrized  Gaussian point cloud
with $N$ points $\{ a_1,-a_1,a_2,-a_2,\dots, a_{N/2}, -a_{N/2}\}$
(and the $a_i $iid standard normal). For what values
of $k$ is the resulting convex hull $k$-centrally neighborly?

\item[Q2 (Setting 1.1.2)] At particular values of $N$ and $k$,
what values of $n$ are associated
with at least a 99\% success rate in
recovering a $k$-sparse object from
$n$ random questions?
Similarly, what values of $n$ are associated
with at least a 99\% success rate in
recovering a $k$-sparse object from
$n$ random questions supposing we know 
that the $k$-sparse object is nonnegative?
\end{description}

These questions can be answered
by establishing the bounds
\[
f_k(T^{N-1}) - E f_k(AT^{N-1}) < 1, \;\;\mbox{ or }\;\;
f_k(C^{N}) - E f_k(AC^{N}) < 1, \qquad \mbox{ for } (Q1)
\]
and 
\[
\frac{E f_k(AT^{N-1})}{f_k(T^{N-1})} \geq 0.99, \;\;\mbox{ or }\;\;
     \frac{E f_k(AC^N)}{f_k(C^N)} \geq 0.99, \qquad \mbox{ for } (Q2);
\]
for the given $(k,n,N)$ of interest.
Simply plugging in the expressions in 
Theorems \ref{thm:strong_precise}
and \ref{thm:weak_precise}, we 
immediately get bounds of the required form.

Since we have developed a series of computational tools
to evaluate the $\Psi_{net}^\star$ and related quantities,
it is rather easy for us to
numerically compute nonasymptotic bounds
answering Q1-Q2. 

\begin{figure}[h]
\begin{center}
\begin{tabular}{cc} 
\includegraphics[width=2.35in]{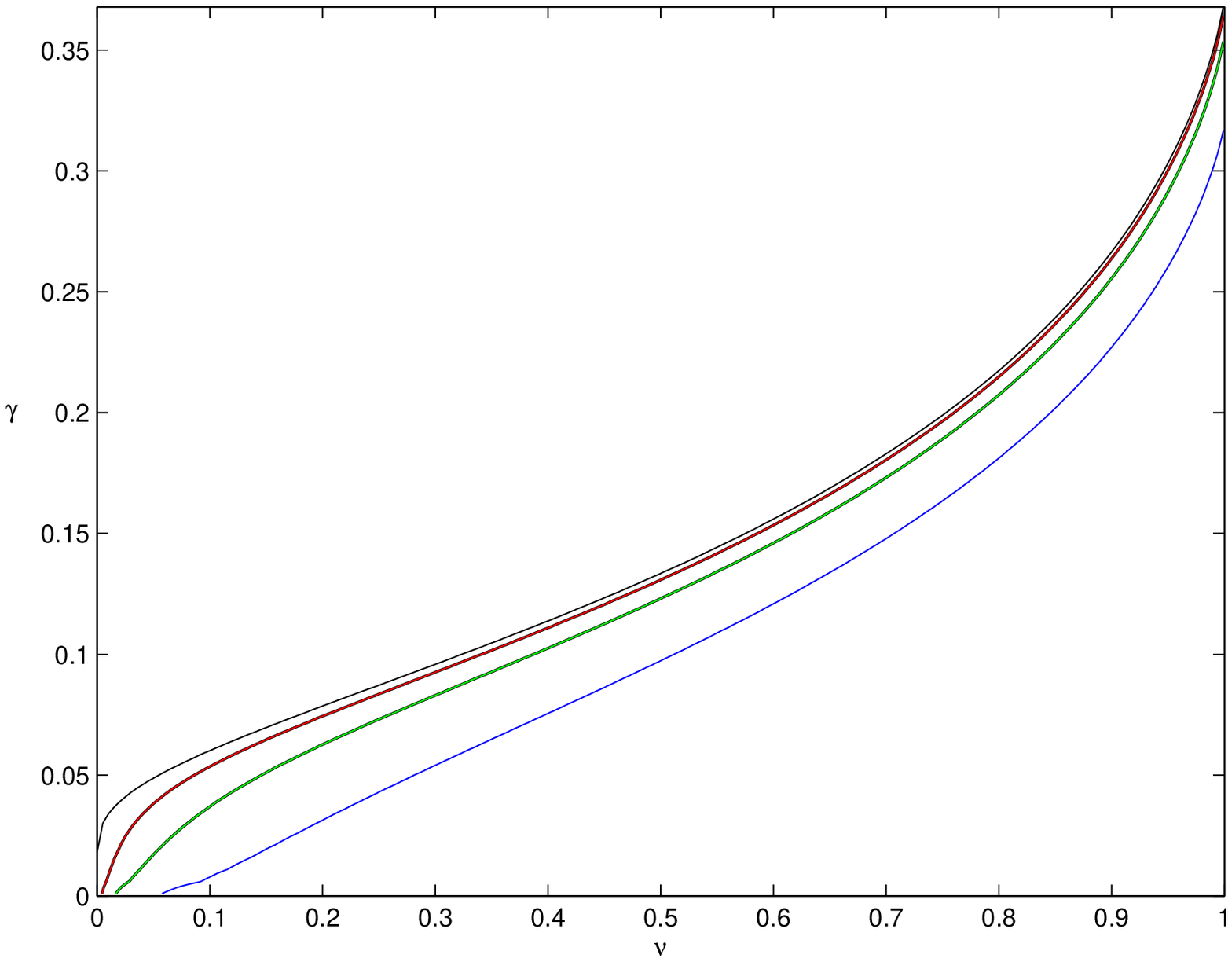} &
\includegraphics[width=2.35in]{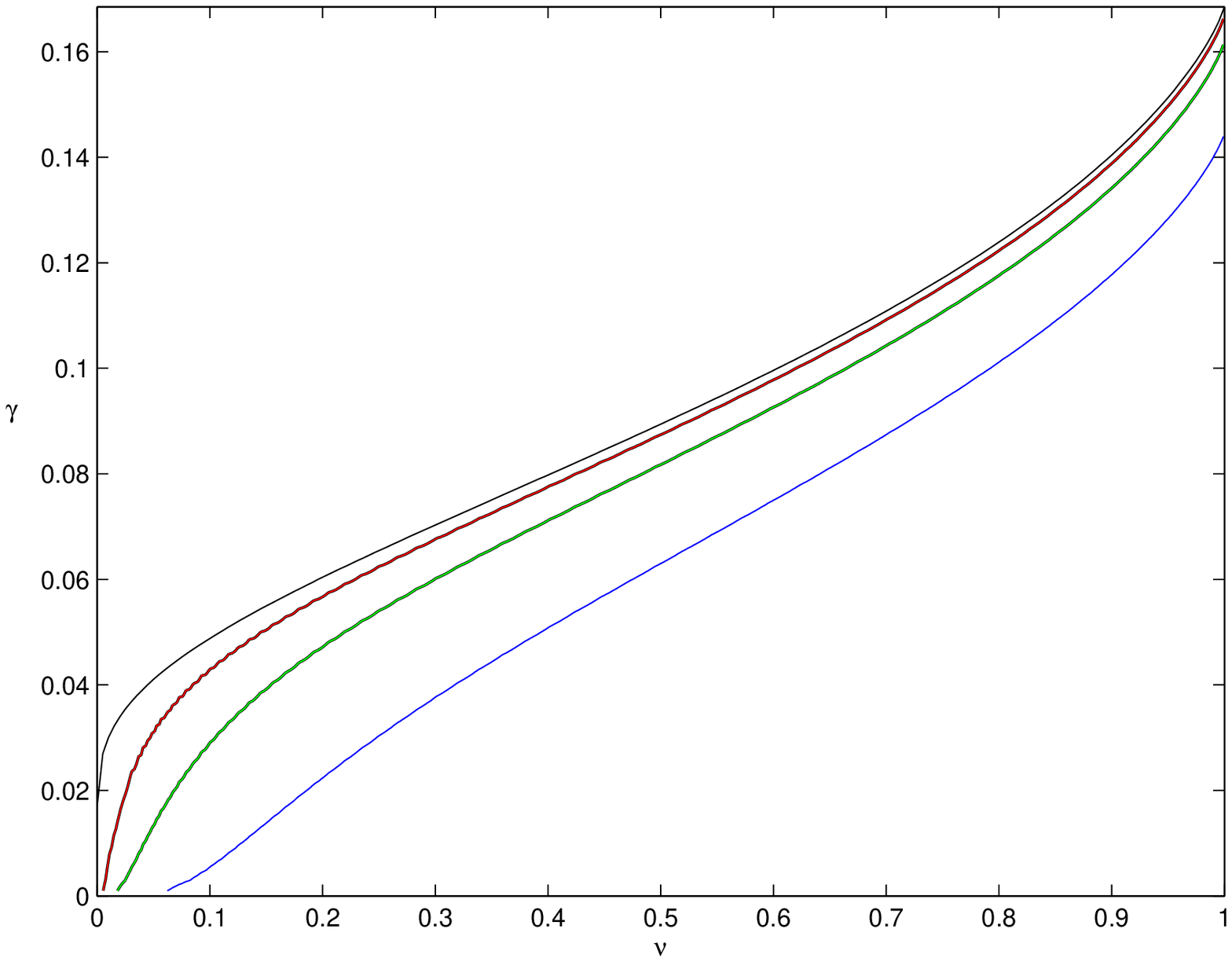} \\
(a) & (b) \\
\includegraphics[width=2.35in]{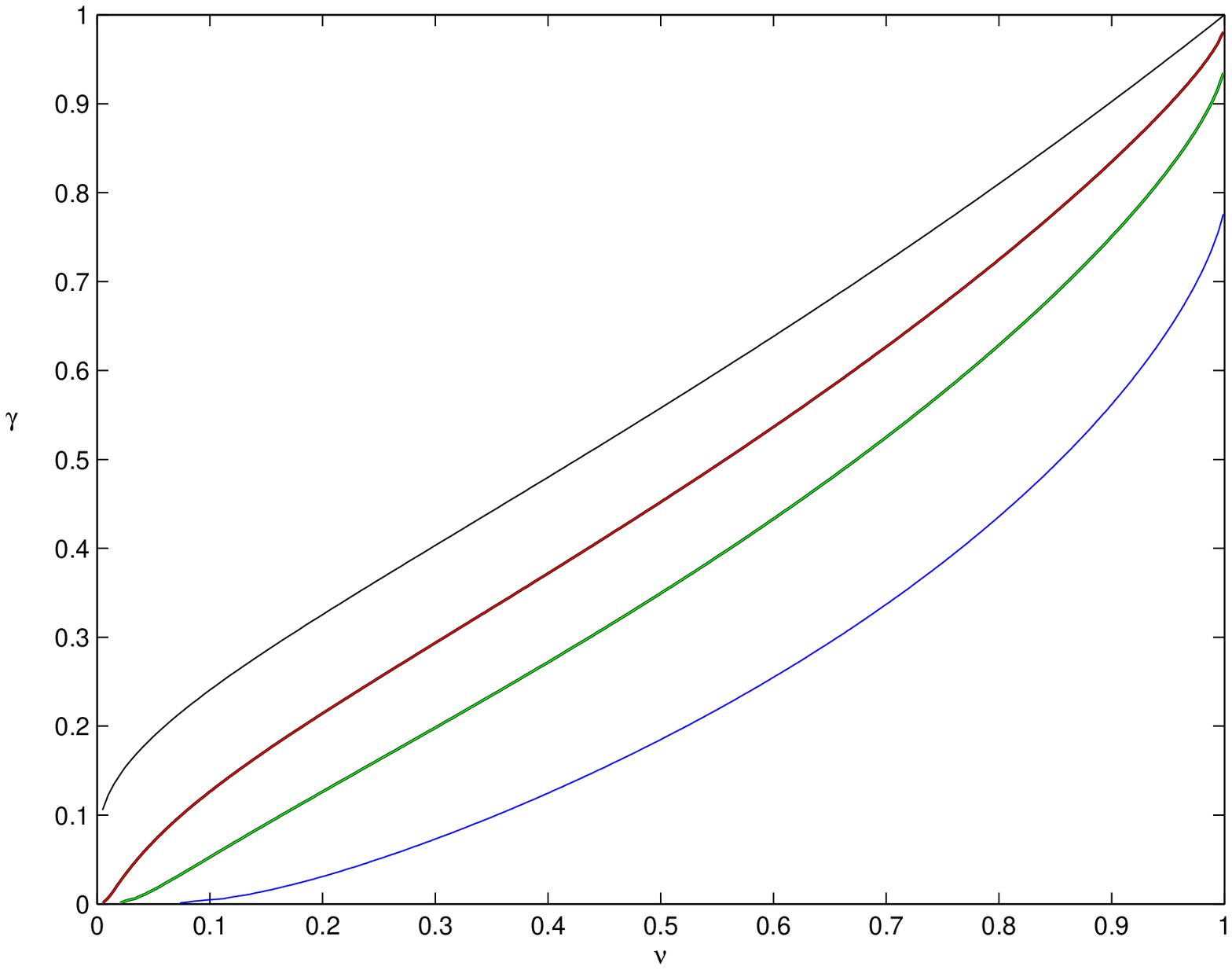} &
\includegraphics[width=2.35in]{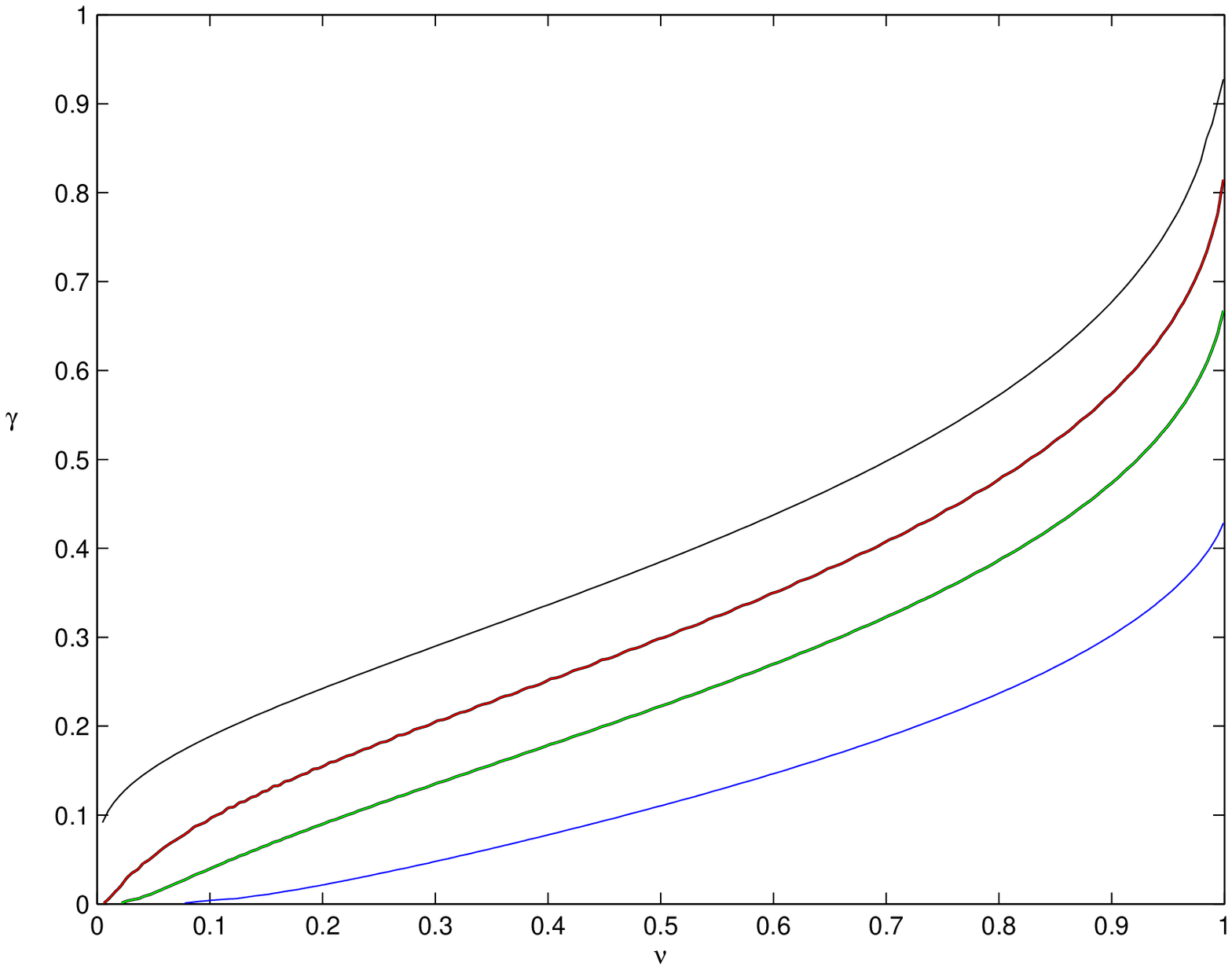} \\
(c) & (d)
\end{tabular}
\end{center}
\caption{Panels (a) and (b): Unit level curves ($|Bound(k,n,N)|=1$) 
for the upper bounds in 
Theorem \ref{thm:strong_precise} for N=200 (blue), 1000 (green), 
and 5000 (red); $Q=T^{N-1}$ (left) and $Q=C^N$ (right).  
The asymptotic, $N\goto\infty$, 
limits $\rho_S^{\star}$ are shown in black.
Panels (c) and (d): The $1/100$ level curves for the lower bounds 
($|Bound(k,n,N)|=1/100$) in Theorem 
\ref{thm:weak_precise} again for N=200 (blue), 1000 (green), 
and 5000 (red); $Q=T^{N-1}$ (left) and $Q=C^N$ (right).  
The asymptotic, $N\goto\infty$, limits 
$\rho_W^{\star}$ are shown in black.}
\label{fig:probability_N}
\end{figure}

Figures \ref{fig:probability_N}(a-b) are relevant to Q1.
They show the unit level set $ Bound(k,n,N) =1$
for the bounds in Theorem \ref{thm:strong_precise} 
for $N=200,1000,$ and $5000$.

\begin{corollary}
For a given $N \in \{200,1000,5000\}$,
consider values of $k$ and $n$ such that $(n/N,k/n)$ lies
strictly beneath the curve corresponding to that $N$ 
depicted in Figures \ref{fig:probability_N}(a).
There exist $n \times N$ matrices $A$ such that $AT^{N-1}$ is 
$k$-neighborly. They can be obtained with
positive probability by random sampling from the standard Gaussian
distribution.

Consider values of $k$ and $n$ such that $(n/N,k/n)$ lies
strictly beneath the curve corresponding to that $N$ 
depicted in Figures \ref{fig:probability_N}(b).
There exist $n \times N$ matrices $A$ so that $AC^N$ is
centrally $k$-neighborly. Such matrices can be obtained with
positive probability by random sampling from the standard Gaussian
distribution on $\bR^{n\times N}$.

\qed
\end{corollary}
  
Figures \ref{fig:probability_N}(c-d) are relevant to Q2.
They show the domain in the phase diagram in which, on average,
at least $99\%$ of faces survive the prescribed dimension 
reduction.  

\begin{corollary}
For a given $N \in \{200,1000,5000\}$,
consider values of $k$ and $n$ such that $(n/N,k/n)$ lies
strictly beneath the curve for that $N$ 
depicted in Figures \ref{fig:probability_N}(c).
Fix a given face $F$ of $T^{N-1}$ independently of $A$.
There is at least a 99\% chance  that $AF$ is a face of $AT^{N-1}$.
Again, we refer to $A$ generated
by random sampling from the standard Gaussian
distribution.

Consider values of $k$ and $n$ such that $(n/N,k/n)$ lies
strictly beneath the curve for that $N$ 
depicted in Figures \ref{fig:probability_N}(d).
Fix a given face $F$ of $C^N$ independently of $A$.
There is at least a 99\% chance  that $AF$ is a face of $AC^N$.
Here probability refers to random sampling from the standard Gaussian
distribution on $\bR^{n\times N}$.

\end{corollary}

Due to the exponentiality of the bounds in 
Theorem \ref{thm:weak_precise}, there are no perceptible 
changes in Figure \ref{fig:probability_N} when the specified levels 
used in calculating those figures are changed,
ie. if we  changed to 
50\% success from $99\%$ success rate in panels (c-d), the figures
would not change substantially.  

It should also be noted from Figure \ref{fig:probability_N} that even 
for small $N$, say 200, when $\delta = n/N$ is relatively large there 
is already 
a large region below the level curves.  However, for $N$ and 
$n/N$ simultaneously small, our bounds
become weak or useless.  For instance, the $N=200$ contour in 
Figure \ref{fig:probability_N}(a) reaches zero at about $n/N=1/20$, 
corresponding to $n=10$.  

In fact, the bounds we presented in 
Theorems \ref{thm:strong_precise} and 
\ref{thm:weak_precise} do not indicate
the full power of our approach.  Those bounds,
in fact, are presented here because they
follow immediately from what has been done above,
and they seem easy for readers to digest. 
For projections to very low dimensional spaces, the polynomial 
factors of the bounds provided in Theorems \ref{thm:strong_precise} and 
\ref{thm:weak_precise} become important; to go beyond
the work reported here, care must be taken 
in combining equations (\ref{eq:psicom_precise}), (\ref{eq:psiext_precise_p}), 
(\ref{eq:psiext_precise_pm}), (\ref{IntGoal1}), 
(\ref{eq:psiface_precise_p_upper}), and (\ref{eq:psiface_precise_pm}) 
to arrive at (\ref{eq:lost_faces}); 
and also perhaps in sharpening the underlying remainder estimates.

\subsection{Relation to Other Work}

We discussed face-counting related work in the body of the text
as the opportunity arose.  We now mention several categories of 
related literature.

\subsubsection{How Neighborly can a Polytope Be?}
Theorems \ref{thm:strongp}, \ref{thm:strongpm}, and \ref{thm:5} 
imply the following.  
{\it For  $N \gg n$, both large, and $N$ subexponential in $n$, there exist
polytopes $P$ which are:
\bitem
 \item    $k$-neighborly with 
 \[
     k \sim \frac{n}{2e \log(N/(n\cdot \sqrt{\pi}))};
 \]
indeed, simply take $P = AT^{N-1}$
 where $A$ has Gaussian i.i.d. entries.
 \item    centrally $k$-neighborly with
  \[
     k \sim \frac{n}{2e \log(N/(n\cdot 2\sqrt{\pi}))};
 \]
 indeed, simply take $P = AC^{N}$
 where $A$ has Gaussian i.i.d. entries.
\eitem
}
Recently, the problem of showing the existence 
of high-dimensional neighborly polytopes has
attracted a resurgence of interest. After fundamental work
in the 1950's-1970's starting with D. Gale \cite{Gale,Gale63} and
extending through P. McMullen and
G.C. Shephard \cite{McMuShep}  and R. Schneider \cite{Schneider},
the subject was very quiet. 
Now, as Schneider
wrote one of us, ``the subject has come to life again''.
Our own work \cite{DoTa05_polytope,DoTa05_signal,Do05_polytope,Do05_signal}
carefully studied the questions of neighborliness and
central neighborliness of projections of random polytopes
in the proportional growth setting.  Our attempt was
 to characterize the
 exact location of the 
asymptotic phase transitions associated with strong 
and weak neighborliness.
Linial and Novik \cite{LinialNovik} gave exponential bounds 
on the probability that $AC^N$ is centrally neighborly; 
note that Rudelson and Vershynin's work \cite{Vershynin}
came earlier and implies similar bounds by duality.  Both 
\cite{LinialNovik,Vershynin} use  
a geometric functional analysis approach which gave inequalities
akin to,
\[
    P \{  f_{k}(AC)  \neq f_{k}(C) \} \leq \psi_1 \exp(- \psi_2 n ), \qquad n > n_0,
\]
valid for $k < cn/ \log(N/n)$ with unspecified constants. 
As our paper was nearing completion, 
we learned that Rudelson and Vershynin \cite{VershyninCISS} had 
been able to supply specific constants.  Their result is as follows.

\begin{theorem}\label{thm:RV}[Rudelson and Vershynin \cite{VershyninCISS}[Theorem 4.1]]
Fix $(k,n,N)$.  Let $A$ be an $n\times N$ random matrix 
from the standard Gaussian distribution.  Let $x_0$ be the 
vector with $k$ nonzeros and let $y=A x_0$.
Let $\Omega$ be the event $\{$The instance of $(P_1)$ defined by 
$(y,A)$ has $x_0$ for its unique solution$\}$.  Then
\begin{equation}\label{eq:RV_probability}
P(\Omega)\ge 1-3.5\exp\left(-\left[\sqrt{n} -\sqrt{m(k,N)}\right]^2/18\right),
\end{equation}
where
\begin{equation}\label{o1}
m(k,N)\le c_1 k\log(c_3 N/k) (1+o(1))
\end{equation}
with $c_1:=6+4\sqrt{2}\approx 11.66$ and $c_3:=e^{3/2}\approx 4.48$.
\end{theorem}

This striking result illustrates the ability to obtain explicit 
constants using existing approaches from geometric functional 
analysis; the proof is admirably short.  Moreover, 
it opens the important question of 
getting explicit results in the finite-sample non-asymptotic case.
When we learned of this result, we decided to include here a
quantitative comparison, illustrating the relative
strengths of our different results.
Theorem \ref{thm:RV} can be recast in a form
similar to the cross-polytope portion of Theorem  
\ref{thm:strong_precise}:
\begin{equation}\label{eq:RV_recast}
P \{  f_{k}(AC^{N})  \neq f_{k}(C^{N}) \} 
\le 3.5\exp\left(-\left[\sqrt{n} -\sqrt{m(k,N)}\right]^2/18\right).
\end{equation}

While this result has the appearance
of an explicit, finite-$N$ result, note the term $o(1)$
in (\ref{o1}) which needs to be converted into
an explicit numerical term to enable concrete
finite-$N$ comparisons. It seems that one can bound this
by $\approx 1/10$ for $N > 1000$.  Rudelson and Vershynin have 
informed us of upcoming explicit bounds on the $o(1)$ term.
It appears that replacing the $o(1)$ term by zero
in the definition of $m(k,N)$ gives a lower bound
on the actual expression developed by Rudelson and Vershynin.
Placing this lower bound for $m$ in the right side
of (\ref{eq:RV_recast}) yields a lower bound on the right hand side of
Rudelson and Vershynin's upper bound.
We now make some comparisons
between our actual upper bounds
and this lower bound on Rudelson's and Vershynin's
upper bound.  To be fair, Rudelson and Vershynin's 
interest was in simply obtaining reasonable bounds by 
geometric functional analysis, which is rather different 
from our focus here.

\begin{figure}[h]
\begin{center}
\begin{tabular}{cc} 
\includegraphics[width=2.35in]{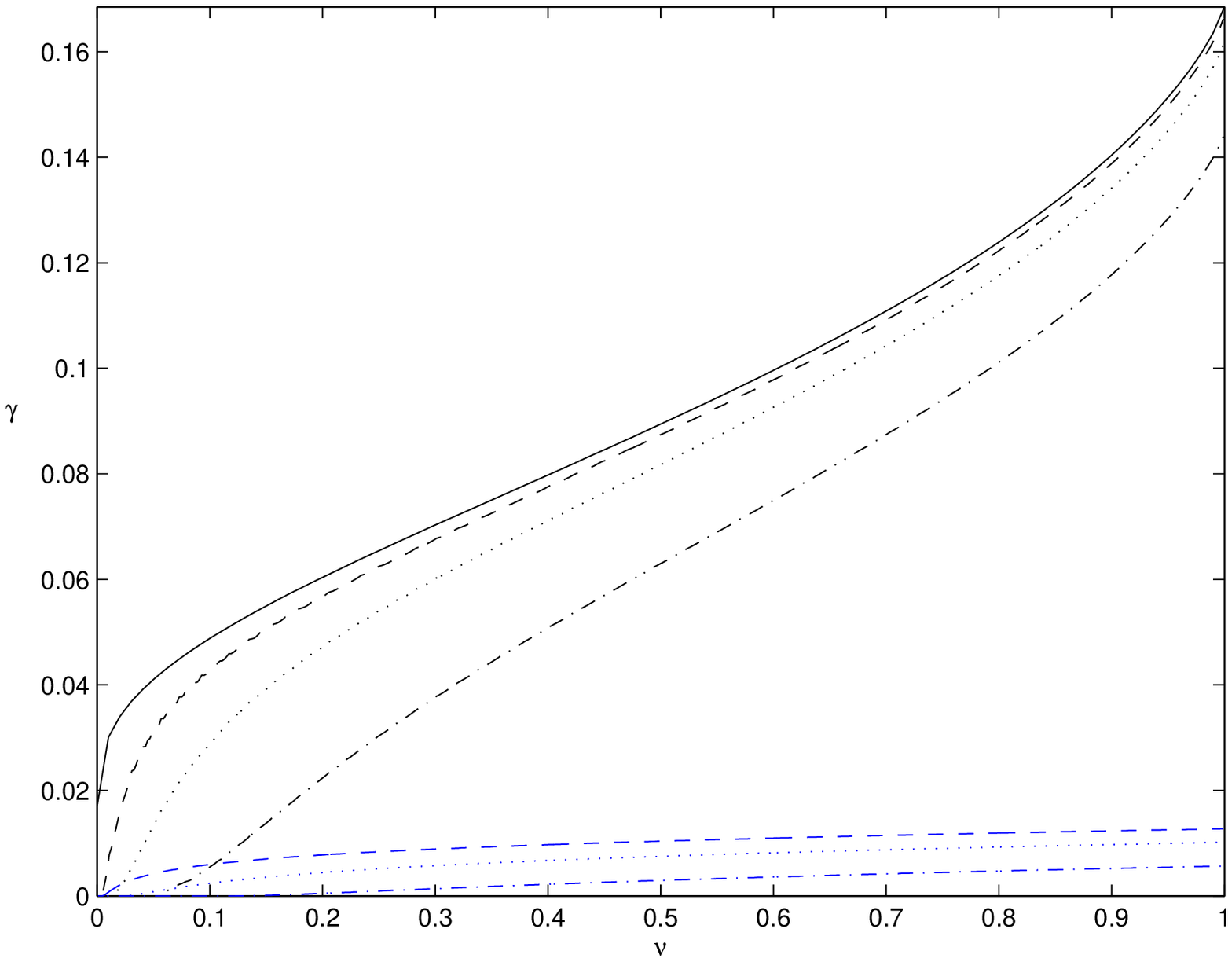} &
\includegraphics[width=2.35in]{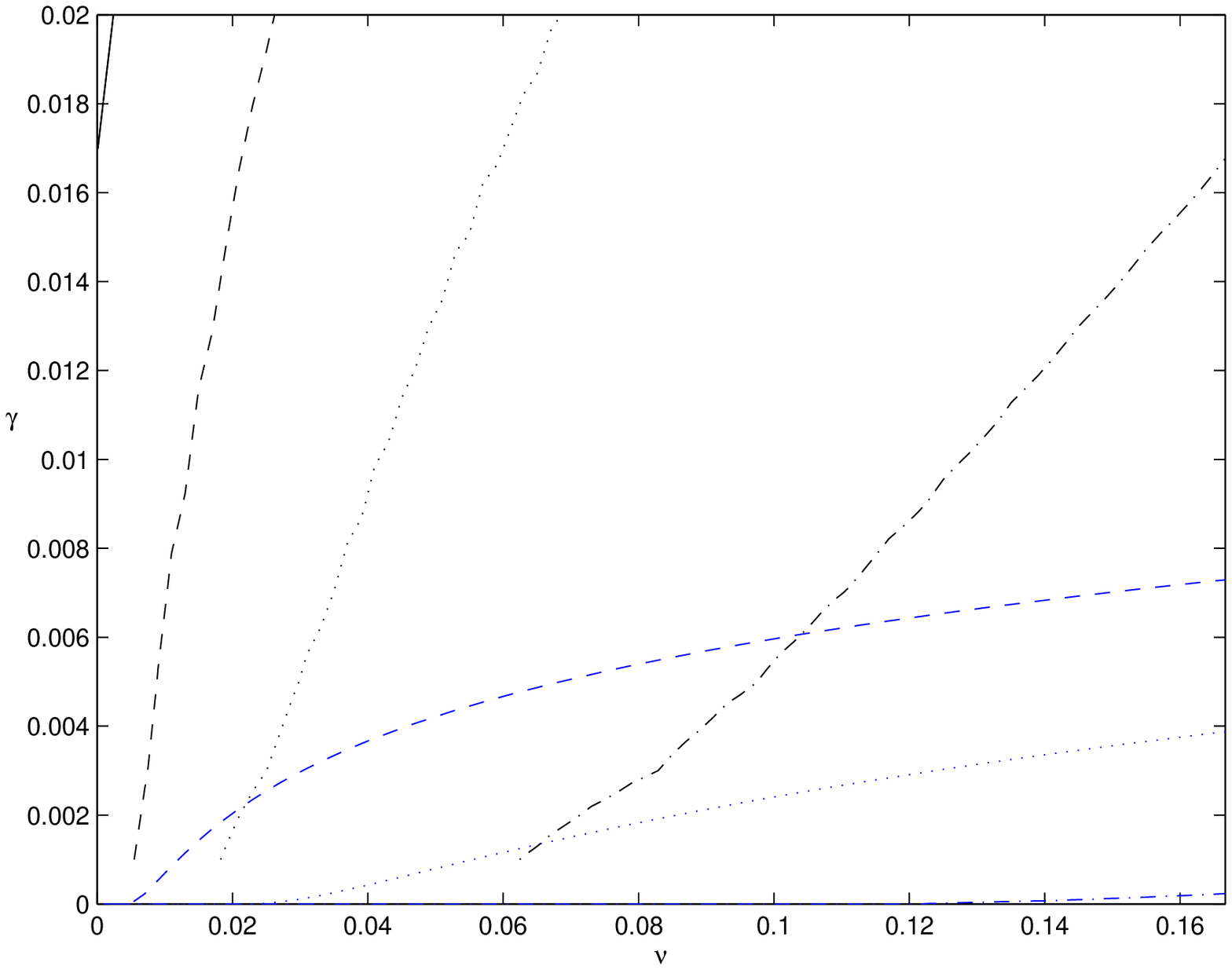} 
\end{tabular}
\end{center}
\caption{Left: Unit level curves $(|Bound(k,n,N)|=1)$ 
for the upper bounds in equation 
(\ref{eq:RV_recast}) (blue) associated with Theorem \ref{thm:RV}, 
and the cross-polytope, $C^N$, portion of 
Theorem \ref{thm:strong_precise} (black).
Curves present the cases $N$=200 (dot-dash), 1000 (dot), 
and 5000 (dash).
The asymptotic, $N\goto\infty$, limit 
$\rho_W^{\pm}$ is shown in solid.  Right: Enhanced bottom left portion 
of left panel.  The upper 4 curves come from our approach; the lower 
three curves come from the results of 
Rudelson and Vershynin \cite{VershyninCISS}.}
\label{fig:RV}
\end{figure}

Figure \ref{fig:RV} illustrates how (\ref{eq:RV_recast}) compares
to the finite-sample bounds developed
by techniques of this paper in Subsection \ref{subsec:exp}.
The curves hugging the bottom
of the display are those implied
by Rudelson-Vershynin's Theorem \ref{thm:RV},
when we replace the $o(1)$ term by 0;
the much higher curves are
those implied by our Theorem \ref{thm:strong_precise}.
The streamlined appearance of the bounds in
Theorem \ref{thm:RV} come at the cost
of a remarkably small region of effectiveness.
For instance, 
the curves associated with $P \{ f_k(AC^N) = f_k(C^N) \} > 0$ 
which follow from Theorem \ref{thm:RV} do not exceed $0.0151$  
for any $N$; whereas, Figure \ref{fig:RV} illustrates that for 
even a modest $N=200$, this same probability level 
following from Theorem \ref{thm:strong_precise} is exceeded well 
before $n/N=1/6$.


\subsubsection{How Many Projections are needed to recover a $k$-sparse object?}
One reason that study of neighborliness  ``has come to life again''
is the surprising implications for speeding up
key processes in medical imaging and proteomics.
The general idea, often labeled {\it Compressed Sensing}
\cite{CS} is that images, spectra, and other real-world objects
are highly compressible, and that this compressibility makes
it possible to reconstruct such objects accurately
from relatively few carefully chosen
generalized samples. In effect Section 1.1.2 has
described an abstract model of compressed sensing.

In the application scenario, $x_0$ represents the
coefficients of an image to be acquired and the rows of $A$
represent a random set of of linear combinations (measurements) 
which will be used to reconstruct $x_0$.
In effect, we are saying that if $x_0$ has $N$ pixels
but only $k \ll N$ nonzeros in (say) a wavelet basis, 
and if $k$ and $N$ are large, then we only
need $n$ measurements, where
\[
   n \geq 2k \log(N/n)(1 + o(1)).
\] 
In contrast, $N$ is the `standard' number of samples; the point is that
for objects which are $k$-sparse with $k$ small,
we can easily have $n \ll N$ if $x_0$ is highly sparse.
(In fact real objects will not 
exhibit such strict sparsity -- $k$ zeros and $N-k$ nonzeros --
but because the $(P_1)$ has an $\ell_1$ stability property \cite{CS},
we can pretend that this is so without distorting the problem.)

The interested reader may pursue
 the papers of Cand\`es and collaborators,
\cite{Candes1,CandesTao}, other theoretical
work \cite{Vershynin,Nowak,TroppGilbert} and 
much recent applied work \cite{ExtCS,baraniuk}.

The quantitative approach developed
here is precise about how much data 
would be needed.  Most of the 
cited theoretical work is qualitative,
often leaving the constants unspecified.
An important point: in Section 1.1.2 and in (1.2) we are studying
the equivalent of {\it weak central neighborliness}.
We argued in Section \ref{int:empirical},
that this is the empirically relevant notion,
we repeat here that ordinary (strong) central neighborliness is simply
not empirically observable.  Nevertheless,
most authors have effectively studied 
implications of ordinary (strong) central neighborliness.
That notion is hard to analyse, and appears to 
indicate a far more pessimistic view of
what is possible than what one actually observes in practice.

\subsubsection{Fast Decoding of Error-Correcting Codes}

In general, decoding of linear error-correcting codes
is NP-hard \cite{Berlek}. However, fast decoding of 
specific error-correcting codes has
been an object of great practical and theoretical attention over
the last 10 years, with great advances in turbo codes
and in LDPC codes (Gallager Codes). 
We proposed in Section 1.1.3 above a simple scheme for fast decoding of
random linear codes over $\bR$ using $(P_1)$. 
The scheme we proposed is equivalent
to one proposed by Cand\`es and Tao \cite{CandesTaoDecode}
and studied further by Rudelson and Vershynin \cite{Vershynin}.
Using the notation of our Section 1.1.3
their decoder solves the $\ell_1$-minimization problem
\[
     \min_u \| w - B^T u \|_1,
\]
yielding the reconstruction $u_1$, say. 
The equivalence of such minimization with the
one proposed in Section \ref{int:ecc} is shown 
in \cite{Do05_signal}[Section 8].

It is of course crucial to know
how many errors such a scheme can correct.
The theoretical literature  (Cand\`es-Tao/Rudelson-Vershynin)
gives qualitative results, saying that one can correct at least
$c n/ \log(N/n)$ errors, with $c$ left unspecified, or else 
specified as a constant which seems much smaller than what would 
be expected based on a comparison of those papers' results with 
the results obtained here.

The problem solved in those visionary papers  is to show existence of
integer-valued matrix pairs $A,B$ allowing block coding of
messages of length $m$ as blocks of length $N$,
such that all patterns of at most $k$ errors can be corrected,
here $n = N-m$.
Our results here change the problem so that $A$ and $B$ are generated 
by partitioning a
uniformly-distributed random projection matrix (n.b. not with
integer-valued entries); with this change, we get a precise asymptotic formula
$k = n  \rho_S^{\pm}(n/N) (1 + o_p(1))$.  For the case $N \gg n$ we have 
proven the formula 
$\rho_S^{\pm}(n/N) \sim 1/2e \log(N/(n\sqrt{\pi}))$.

If we change the problem again slightly so that 
the goal is to correct {\it nearly all} rather than
{\it all} error patterns, then for the case of long block codes,
we get a precise asymptotic
formula  $k = n \rho_W^{\pm}(n/N) (1 + o_P(1))$. For the case $N \gg n$
we have proven the formula $\rho_W^{\pm}(n/N) \sim 1/2 \log(N/n)$.

Conceivably, such results for the ``changed problems''
we just mentioned may be better than
for the original problem; i.e. the situation for general random matrices
may be more optimistic than for matrices with integer entries.
However, our empirical results with Rademacher random
matrices indicate that our formula $n \rho_W^{\pm}(n/N)$
accurately describes the integer-valued case as well, i.e. 
accurately describes the number of errors
which can typically be corrected by 
such random matrices with integer-valued entries.

\section*{Appendix: Proofs of Key Lemmas}
\setcounter{equation}{0}
\setcounter{table}{0}
\setcounter{figure}{0}

\subsection{Proof of Lemma \ref{lem:ygam}.}

We develop (\ref{eq:lemygam}) 
in two stages.  Initially, we derive the
asymptotic behavior of $\sgam$ as $\gamma \goto 0$; we then
substitute that into equation (\ref{eq:ygamma}).
To motivate our approximation of $\sgam$ we use an 
asymptotic series for $R(s)$ appropriate for the regime of $s$ large,
\[
R(s):=s e^{s^2/2}\int_{s}^{\infty}e^{-y^2/2}dy=1-\frac{1}{s^2}
+\frac{1\cdot 3}{s^4}
-\frac{1\cdot 3\cdot 5}{s^6}
+\frac{1\cdot 3\cdot 5\cdot 7}{s^8}
+\cdots;
\]
This is derived as follows. The ratio $R(s) = s \cdot Mills(s)$ for $s > 0$,
where $Mills(s)$ is the usual Mills' ratio for the normal
distribution. 
The corresponding asymptotic series for Mills' ratio
is developed in \cite[Secs 5.37,5.38]{KendallVol1};
H. Ruben \cite{Ruben} credits this series to Laplace.

In \cite[Eq. (5.106)]{KendallVol1} it is shown that 
the error in truncating the series for $Mills()$ at the $s$-th term 
is at most as large as the $s$-th term itself.  $R()$ inherits this property.

It is now convenient to define $L(s,\gamma) := R(s) - 1 + \gamma$
and note that $s_\gamma$ is defined by $L(s_\gamma,\gamma) = 0$.

Keeping the first two terms in the series expansion for $R(s)$ 
and applying the bounds from  \cite[Eq. (5.106)]{KendallVol1}
yields $L(s,\gamma)=\gamma-s^{-2} +3s^{-4}+{\cal O}(s^{-6})$,
which suggests the approximation,
\begin{equation}\label{eq:tsgam}
\sgam \approx \tsgam:= \gamma^{-1/2}-\frac{3}{2}\gamma^{1/2} .
\end{equation}

To quantify the error in this approximation,
invoke the mean value theorem; 
given a smooth function $F(x)$, there is always
a point $w\in[\min (x,y),\max(x,y)]$
satisfying
\begin{equation}\label{taylor}
F(y)=F(x)+(y-x)\frac{d }{d y}F(y)_{|y=w}.
\end{equation}
Hence we can bound $|y-x|$ if we have 
suitable bounds on $|F(y)-F(x)|$ and 
$\frac{d}{dy} F(y) |_{y=w}$.
Apply this principle to $F(s) = L(s,\gamma)$ about $\sgam$,
getting
\begin{equation}\label{sprebound}
|\tsgam -\sgam|\le \left|L(\tsgam,\gamma)\left/ \frac{\partial }{\partial s}L(s,\gamma)_{|s=s_{mid}}\right.\right|,
\end{equation}
for some point $s_{mid}\in[\min(\sgam,\tsgam),\max (\sgam,\tsgam)]$.

The following bounds follow from
 \cite[Eq. (5.106)]{KendallVol1}
\begin{equation}\label{eq:mills}
1-s^{-2}+\frac{5}{2}s^{-4}<R(s)<1-s^{-2}+3s^{-4}\quad\quad\mbox{for}\quad s>\sqrt{30},
\end{equation}
yielding in turn
\begin{equation}\label{eq:Ltsgam}
|L(\tsgam,\gamma)|<\frac{1}{2}\gamma^2, \quad\quad\mbox{for}\quad\gamma<1/30.
\end{equation}

To bound the denominator, note that
\[
\frac{\partial }{\partial s} L(s,\gamma)=[s+s^{-1}]\cdot R(s)-s,
\]
which is a positive decreasing function of $s$; 
this attains its lower bound on the interval $s \in [\min(\sgam,\tsgam),\max (\sgam,\tsgam)]$
at one of the endpoints $\{ \sgam, \tsgam \}$.  At $\tsgam$ we again make use 
of the lower bound on Mills' ratio in equation (\ref{eq:mills})
\begin{eqnarray}
\frac{\partial }{\partial s} L(s,\gamma)_{|s=\tsgam} & = & 
[\tsgam+\tsgam^{-1}] R(\tsgam)-\tsgam \nonumber \\
& > & [\tsgam+\tsgam^{-1}](1-\tsgam^{-2}+\frac{5}{2}\tsgam^{-4}) -\tsgam \nonumber \\
& = & \frac{3}{2}\tsgam^{-3}+\frac{5}{2}\tsgam^{-5}>\frac{3}{2}\gamma^{3/2}.
\end{eqnarray}
For the lower bound at $\sgam$ we assume $|\tsgam-\sgam|\le\frac{1}{2}\gamma^{1/2}$ (which we will verify momentarily), which gives the upper bound  $\sgam\le\gamma^{-1/2}-\gamma^{1/2}$.  From this we have the lower bound,
\begin{eqnarray}
\frac{\partial }{\partial s} L(s,\gamma)_{|s=\sgam} & = &
[\sgam+\sgam^{-1}]\cdot (1-\gamma)-\sgam = (1-\gamma)\sgam^{-1}-\gamma\sgam 
\nonumber \\
& \ge & (1-\gamma)\cdot\frac{1}{\gamma^{-1/2}-\gamma^{1/2}}-\gamma(\gamma^{-1/2}-\gamma^{1/2}) =\gamma^{3/2}. \nonumber
\end{eqnarray}
Using these bounds in equation (\ref{sprebound}) we have:
\begin{equation}\label{eq:tsgam_one}
|\sgam-\tsgam|\le \frac{1}{2}\gamma^{1/2},
\quad\mbox{ for} \quad \gamma\le 1/30, 
\end{equation}
which justifies the earlier claim that $|\sgam-\tsgam|\le \frac{1}{2}\gamma^{1/2}$.  For the following calculations the following estimate suffices,
\begin{equation}\label{sgamma}
\sgam=\gamma^{-1/2}+r_1(\gamma),\quad\mbox{with}\quad |r_1(\gamma)|\le 2\gamma^{1/2},\quad\mbox{for}\quad \gamma\le 1/30.
\end{equation}

%
Combined with (\ref{eq:ygamma}),  this gives 
(\ref{eq:lemygam}) and hence,  Lemma \ref{lem:ygam}.
\qed
\subsection{Proof of Lemma \ref{lem:xnu}}

We first motivate our approximation for $\xnu$,
which solves $2xQ(x)/q(x) = 1 - \nu^{-1}$.
The truncated asymptotic series
\[
\frac{2x Q(x)}{q(x)}=2\pi^{1/2}x e^{x^2}-1+{\cal O}(x^{-2}), \quad x \goto \infty,
\]
suggests approximating $x_\nu$ as the solution to 
\begin{equation}\label{eq:qtrunk}
x e^{x^2}-\frac{1}{2}\pi^{-1/2}\nu^{-1}=0;
\end{equation}
this is exactly of the form
 (\ref{eq:fund}) with $z = \zp:=(2\nu \sqrt{\pi})^{-1}$. 
Our approach for approximate solution
of   (\ref{eq:fund}), carried out to two stages,
 yields the approximant $\xnu$, obeying:
\begin{equation}
\txnu^2:=\log \zp -\frac{1}{2}\log\log \zp.
\end{equation}
Our claim that $\txnu$ accurately approximates 
$\xnu$ as $\nu\rightarrow 0$, 
as stated in Lemma \ref{lem:xnu}, 
will be supported by arguments
similar to those used in proving Lemma \ref{lem:ygam}.  

Let 
\begin{equation}\label{h}
J(x,\nu):=2xe^{x^2}\int_{-\infty}^x e^{-y^2} dy +1-\nu^{-1},
\end{equation}
whose level curve $J(\xnu,\nu)=0$ defines $\xnu$.  
To bound the error in the approximation, $\txnu$, 
we again use the mean value approach (\ref{taylor}), getting
\begin{equation}\label{preboundh}
|\txnu -\xnu|\le \left|J(\txnu,\nu)\left/ 
\frac{\partial }{\partial x} J(x,\nu)_{|x=x_{mid}} \right. \right|,
\end{equation}
for some point $x_{mid}\in[\min(\xnu,\txnu),\max(\xnu,\txnu)]$.
The magnitude of $J(\txnu,\nu)$ can be bounded by
\begin{eqnarray}
|J(\txnu,\nu)| & = & -J(\txnu,\nu) 
= -2\pi^{1/2}\txnu e^{\txnu^2}+\nu^{-1}+r_6(\nu) \nonumber \\
 & \le & -2\pi^{1/2}\txnu e^{\txnu^2}+\nu^{-1}, \qquad  \nu<1/10 \nonumber \\
& = & \nu^{-1}\left[ 1-\left(1-\frac{1}{2}\frac{\log\log\zp}{\log\zp}\right)^{1/2}\right]  \label{eq:identsurp} \\
& \le & \nu^{-1}\frac{3}{8}\frac{\log\log\zp}{\log\zp} 
 \label{hbound}
\end{eqnarray}
where the transition from the first to second line 
utilizes $r_6(\nu):=2\txnu e^{\txnu^2}\int_{x}^{\infty}e^{-y^2}dy-1\le 0$ for $\nu<1/10$.

Turning to the denominator in (\ref{preboundh}), we observe that
on the half-line  $x\ge 0$
the derivative is a positive increasing function,
\begin{eqnarray}\label{hderiv}
\frac{\partial }{\partial x}J(x,\nu) & = & 
2x+2(1+2x^2)e^{x^2}\int_{-\infty}^x e^{-y^2}dy \nonumber \\
& = & \frac{1+2x^2}{x}[J(x,\nu)+\nu^{-1}]-x^{-1};
\end{eqnarray}
a lower bound for $\partial J/\partial x$ over  $[\min(\xnu,\txnu),\max(\xnu,\txnu)]$
is attained at one of the endpoints $\xnu$ or $\txnu$.
At $\xnu$ a simple lower bound is
\begin{equation}\label{hderiv_xnu}
\frac{\partial }{\partial x}J(x,\nu)_{|x=\xnu}
=\frac{1+2\xnu^2}{\xnu}\nu^{-1}-\xnu^{-1}
\ge 2\xnu \nu^{-1} \;\;\mbox{ for } \nu\le 1.
\end{equation}
A similar lower bound holds at $\txnu$,
\begin{eqnarray}
\frac{\partial }{\partial x}J(x,\nu)_{|x=\txnu} & = & 
\frac{1+2\txnu^2}{\txnu}[J(\txnu,\nu)+\nu^{-1}]-\txnu^{-1} \nonumber \\
& \ge & \frac{1+2\txnu^2}{\txnu}\nu^{-1}\left[1-\frac{1}{2}\frac{\log\log\zp}{\log\zp}\right] -\txnu
\qquad [\mbox{by (\ref{eq:identsurp})}] \nonumber \\
& = & \txnu\nu^{-1}+\txnu\nu^{-1}
\left[1+(\log\zp)^{-1}-\frac{\log\log\zp}{\log\zp}-\nu\right] \nonumber \\
& \ge & \txnu\nu^{-1} \;\;\mbox{ for } \nu\le 1/4. \label{hderiv_tilde}
\end{eqnarray}
Combining 
(\ref{hderiv_xnu}) and (\ref{hderiv_tilde}),
\begin{equation}\label{hderiv_bound}
\frac{\partial }{\partial x} J(x,\nu)_{|x=x_{mid}}
\ge \nu^{-1} \mbox{min}(\xnu,\tilde{x}_\nu);
\end{equation}
although crude, this is sufficient for later purposes.  

Shortly we will prove
there is $\nu_0 > 0$ such that 
\begin{equation} \label{eventually}
\mbox{min}(\xnu,\txnu)\ge \frac{3}{4}\txnu, \qquad 0 < \nu < \nu_0.
\end{equation}
Substituting (\ref{hbound}) and (\ref{hderiv_bound}) 
into equation (\ref{preboundh}) gives
\begin{equation}
|\xnu-\txnu|\le \frac{1}{2}\txnu^{-1}\frac{\log\log\zp}{\log\zp}.
\end{equation}
Lemma \ref{lem:xnu} follows by simple substitution of terms.

We now show (\ref{eventually}).
Recall that $\int_{-\infty}^\infty e^{-y^2} dy = \sqrt{\pi}$. Hence on $x \geq 0$,
$J(x,\nu) \leq \tilde{J}(x,\nu) :=  2\sqrt{\pi} x e^{x^2} + 1 - \nu^{-1}$.
As $J(x,\nu)$ is monotone increasing on $(0,\infty)$ it follows
that  $\tilde{J}(x',\nu) < 0$ implies $x_\nu > x'$.

We now show that if  $0 < a < 1$, then 
\begin{equation} \label{h-negative}
\tilde{J}(a \tilde{x}_\nu, \nu) < 0
\end{equation}
for $\nu$ sufficiently small. Setting $a =3/4$, this will imply
$J( \frac{3}{4} \tilde{x}_\nu, \nu) < 0$ for all sufficiently small $\nu$,
and so, for such $\nu$, $\mbox{min}(\xnu,\txnu)\ge \frac{3}{4}\txnu$;
(\ref{eventually}) follows. 

Proceed thus:
\begin{eqnarray*}
\tilde{J}(a \tilde{x}_\nu, \nu) &=&     2 \sqrt{\pi}  \cdot \frac{ a \sqrt{ \log(z^+) - 1/2 \log \log z^+ } }{ \log(z^+)^{a^2/2}} \cdot (z^+)^{a^2}  + 1 - \nu^{-1}\\
                                                 &=&     2 \sqrt{\pi} a  \cdot \log( \nu^{-1} )^{(1- a^2)/2}(1 + o(1) ) \cdot \nu^{-a^2}  + 1 - \nu^{-1} \\
                                                 &=&    o( \nu^{-1})  + 1 - \nu^{-1} , \quad \nu \goto 0.
\end{eqnarray*}
(\ref{h-negative}) follows.
\qed

\subsection{Proof of Lemma \ref{lem:psinetp_nu}}

We will show that 
\begin{equation} \label{GoalLemma44}
M[\frac{\partial }{\partial \nu}\psinetp] (\delta,\rps)\leq\frac{1}{2}\log\left(\frac{2e}{\cs}\right) 
+ o(1) ,\quad\delta\rightarrow 0.
\end{equation}
Because $\cs>2e$, the leading term on the RHS is a negative constant,
showing that for small enough $\delta$ the function $\Psi_{net}^+$
is monotone decreasing in $\nu$ on 
the admissible domain, implying the assertions
of the Lemma.
Now
\begin{eqnarray}
\frac{\partial }{\partial \nu}\psinetp(\nu,\gamma)
& = & \frac{\partial}{\partial \nu}(\psicomp-\psiintp-\psiextp) \label{eq:psinet_nu_deriv} \nonumber \\
& = & \log\xnu+\frac{1}{2}\log(4\pi)+H(\gamma)-(1-\gamma)
\left[\log\left(\frac{\ygam}{\gamma}\right)+\frac{1}{2}\log(2\pi)
+\frac{\gamma-1}{2\gamma}\ygam^2\right]. \nonumber
\end{eqnarray}

Over the interval  $\nu\in [\delta,1)$, 
the first component, $\log\xnu$, is largest 
at $\nu=\delta$.  Applying Lemma \ref{lem:xnu} we have
\begin{eqnarray}
\log\xnu\le \log\xdelta & = & \frac{1}{2}\log\log\zpd
+\log \left( \left[ 1-\frac{1}{2}\frac{\log\log\zpd}{\log\zpd}\right]^{1/2}
+r_3(\delta)(\log\zpd)^{-1/2}\right) \nonumber \\
& < & \frac{1}{2}\log\log\zpd \quad\quad \mbox{for}\quad\delta<1/50.
\label{eq:Boundlogx}
\end{eqnarray}

The RHS of (\ref{eq:psinet_nu_deriv}) is an increasing function of $\gamma$, 
maximized at $\gamma=\rps$.  Using Lemma \ref{lem:ygam},
gives, for all $\gamma$ small enough:
\begin{eqnarray}
& & H(\gamma)-(1-\gamma)
\left[\log\left(\frac{\ygam}{\gamma}\right)+\frac{1}{2}\log(2\pi)
+\frac{\gamma-1}{2\gamma}\ygam^2\right] \nonumber \\
& & < \frac{1}{2}\log{\gamma}+\frac{1}{2}\log(e/2\pi)+6\gamma-\frac{1}{2}\gamma\log\gamma  \label{eq:Boundrest} \\
& & = \frac{1}{2}\log\left[\frac{e}{2\pi \cs \log\zpd}\right]
+{\cal O}\left(\frac{\log\log\zpd}{\log\zpd}\right),\quad\delta\rightarrow 0. \nonumber
\end{eqnarray}

Combining (\ref{eq:psinet_nu_deriv})-(\ref{eq:Boundrest})
  yields (\ref{GoalLemma44}).
\qed

\newcommand{\psinets}{\Psi_{net}^{\star}}
\subsection{Proof of Lemma \ref{lem:psinetpgamma}}

It is sufficient to show that for some $\gamma_0 > 0$
\begin{equation}
\label{Goal45}
\frac{\partial}{\partial\gamma} \Psi_{net}^\star(\nu,\gamma) > \nu/2, \quad \nu \in [\delta,1), \quad 0 <  \gamma < \gamma_0.
\end{equation}
Now
\begin{eqnarray}
\frac{\partial}{\partial\gamma}\psinets(\nu,\gamma) & = & \nu\left[
\gamma^{-1}-2\log\gamma+\log(\ygam)
+\frac{\gamma^2 -1}{2\gamma^2}\ygam^2+\log(1-\gamma)\right.
\nonumber \\
& & \left. +\frac{1}{2}\log(2\pi)+1
+(1-\gamma)^2\left(\frac{\ygam}{\gamma}-\frac{1}{\ygam (1-\gamma)}\right)
\frac{d}{d\gamma}\ygam\right]. \label{eq:psinetgamma_def}
\end{eqnarray}
Lower bounds for each but the last term follow 
either directly or from Lemma \ref{lem:ygam}.  For $\gamma<1/10$, $\ygam$ 
satisfies
\begin{equation}\label{eq:ygamtmp}
\ygam\ge\frac{\gamma^{1/2}}{1-\gamma}-4\gamma^{3/2}>\gamma^{1/2}(1-3\gamma)
\end{equation}
and
\[
\ygam\le\frac{\gamma^{1/2}}{1-\gamma}+4\gamma^{3/2}<\gamma^{1/2}(1+6\gamma),
\]
from which follow both
\[
\log\ygam>\frac{1}{2}\log\gamma-4\gamma,
\]
and
\[
\frac{\gamma^2 -1}{2\gamma^2}\ygam^2>\frac{-1}{2\gamma}(1+6\gamma)^2(1-\gamma^2)
\frac{-1}{2\gamma}(1+16\gamma),
\]
respectively.  The last term in (\ref{eq:psinetgamma_def}) 
requires estimating  
\[
\frac{d}{d\gamma}\ygam=\frac{\sgam}{(1-\gamma)^2}\left[
1-\frac{\gamma (\gamma-1)}{\gamma\sgam^2+\gamma-1}\right].
\]
From (\ref{eq:ygamtmp}) 
\[
\gamma (\gamma-1)>\gamma\sgam^2+\gamma-1>4\gamma^2-3\gamma>4\gamma (\gamma-1)
\]
yielding 

\begin{equation}\label{eq:ygamprime}
0<\frac{d}{d\gamma}\ygam<\frac{3}{4}\frac{\sgam}{(1-\gamma)^2}
<\frac{3}{4}\gamma^{-1/2}+\gamma^{1/2} ,
\end{equation}
for $\gamma\le 1/30$.  As the above quantity is positive, a lower bound 
for the last term in (\ref{eq:psinetgamma_def}) is obtained with a lower 
bound on its multiplicative factor,
\[
(1-\gamma)^2\left(\frac{\ygam}{\gamma}-\frac{1}{\ygam (1-\gamma)}\right) 
 > -8\gamma^{1/2}(1-\gamma)^2
\]
which is obtained from (\ref{eq:ygamtmp}).  With 
(\ref{eq:ygamprime}) we arrive at,
\[
(1-\gamma)^2\left(\frac{\ygam}{\gamma}-\frac{1}{\ygam (1-\gamma)}\right)
\frac{d}{d\gamma}\ygam>-6-8\gamma\quad\mbox{for}\quad\gamma<1/30.
\]
Combining these bounds we have that
\begin{equation}\label{eq:psinet_gamma}
\frac{\partial}{\partial\gamma}\psinets(\nu,\gamma) > \nu\left[
\frac{1}{2}\gamma^{-1}+\frac{3}{2}\log(1/\gamma)-13-14\gamma\right];
\end{equation}
for $\gamma<1/30$, the term in brackets exceeds $1/2$.
(\ref{Goal45}) follows.
\qed

\subsection{Proof of Lemma \ref{lem:psiweakp_nu}}

We will show that over the admissible domain,
\begin{equation} \label{Goal46}
\frac{\partial}{\partial\nu}(\psinetp-\psiweakp) < \left[ \frac{1}{2}\log\left(\frac{2e}{\cw}\right)-\frac{1}{\cw}+ o(1) \right], \quad  \delta \goto 0.
\end{equation}
As $\cw > 2$, this proves the Lemma.
For sufficiently small $\delta$,
we have the inequality
\begin{eqnarray}
\frac{\partial}{\partial\nu}(\psinetp-\psiweakp)
& = & \frac{\partial}{\partial \nu}(\psinetp)
+\gamma\log\nu+\gamma\log\gamma-\gamma\log(1-\nu\gamma)
\nonumber  \\
& < & \log\xnu+\gamma\log\nu+\frac{1}{2}\log\gamma
+\frac{1}{2}\log(2e)+8\gamma+\frac{1}{2}\gamma\log\gamma \nonumber \\
&:=& \Omega(\nu,\gamma), \label{eq:psiweak_nu_deriv}
\end{eqnarray}
say.
We will show that 
\[
  M[\Omega](\delta, \rwp) < \left[ \frac{1}{2}\log\left(\frac{2e}{\cw}\right)-\frac{1}{\cw} \right], \quad 0 <\delta <\delta_0; 
\]
this implies (\ref{Goal46}).
We first note that
\begin{eqnarray}
\frac{\partial}{\partial\gamma}\Omega (\nu,\gamma) & = & \log\nu+\frac{17}{2}+\frac{1}{2}\log\gamma+\frac{1}{2}\gamma^{-1} \nonumber \\
& \ge & \log\delta+\frac{17}{2}+\frac{1}{2}\log(\gamma)+\frac{1}{2}\gamma^{-1}\nonumber \\
& \ge & \left[\frac{\cw}{2}-1\right]\log(1/\delta)-\frac{1}{2}\log\log(1/\delta)+\frac{17}{2}-\frac{1}{2}\log (\cw)
\end{eqnarray}
which for any $\cw>2$ becomes arbitrarily large as $\delta$ approaches zero.
As a result, $\Omega(\nu,\gamma)$ obtains its maximum where $\gamma$ is largest
within the admissible domain,
i.e. at $\gamma=\rwp$.  To find the overall maximum we now 
examine the $\nu$ direction along $\gamma=\rwp$:
\begin{equation}\label{eq:Bnu}
\frac{\partial}{\partial\nu}\Omega(\nu,\gamma) =
\frac{\frac{\partial}{\partial\nu}\xnu}{\xnu}+\frac{\gamma}{\nu}
 =  \nu^{-1}\left[\frac{1}{\cw\log{1/\delta}}-\frac{1}{1+2\xnu^2-\nu}\right].
\end{equation}
>From Lemma \ref{lem:xnu} it follows that for any $\cw>2$, 
for $\delta$ sufficiently small,
\[
\cw\log(1/\delta)>2\log(1/\delta)>2\xdelta^2+1-\delta>2\xnu^2+1-\nu
\]
for $\nu\in[\delta,1)$.  As a result (\ref{eq:Bnu}) is negative for
$\delta$ sufficiently small, indicating that the maximum of $\Omega(\nu,\gamma)$
over the domain of interest is obtained at $(\delta,\rwp)$.  Moreover,
\[
\Omega(\nu,\gamma)\le \Omega(\delta,\rho)  < \frac{1}{2}\log\left(\frac{2e}{\cw}\right)-\frac{1}{\cw}+{\cal O}\left(\frac{\log\log\zp}{\log\zp}\right),\quad\delta\rightarrow 0
\]
giving  (\ref{Goal46}).
\qed

\subsection{Proof of Lemma \ref{lem:psiweakp_gamma}}

>From (\ref{eq:psinet_gamma}) and 
\[
\frac{\partial}{\partial\gamma}\psiweakp(\nu,\gamma)=-\nu\left[\log\gamma+\log\nu-\log(1-\nu\gamma)\right]
\]
we have the lower bound,
\begin{eqnarray}
\frac{\partial}{\partial\gamma}(\psinetp-\psiweakp)(\nu,\gamma)
& > & \nu\left[\frac{1}{2}\gamma^{-1}+\log\nu+\frac{1}{2}\log(1/\gamma)-13-14\gamma\right]
\nonumber \\
& > & \delta\left[\left(\frac{\cw}{2}-1\right)\log(1/\delta)+\frac{1}{2}\log\log(1/\delta)-13 \right. \nonumber \\
&& \qquad \left. +\frac{1}{2}\log \cw-14 [\cw\log(1/\delta)]^{-1} \right] ; \nonumber 
\end{eqnarray}
with the last inequality due to $\gamma\le\rwp$.
For any $\cw > 2$ the above bound is positive for $\delta$ sufficiently small.

\qed

\bibliography{HowManyRP_20060917}

\bibliographystyle{amsplain}
\end{document}